\documentclass[letterpaper, 10 pt, conference]{ieeeconf}%
\usepackage{graphics}
\usepackage{epsfig}
\usepackage{mathptmx}
\usepackage{amsmath}
\usepackage{amssymb}
\usepackage{setspace}
\usepackage[colorinlistoftodos]{todonotes}
\usepackage[colorlinks=true, allcolors=blue]{hyperref}
\usepackage{caption}
\usepackage{subcaption}
\usepackage{hyperref}
\usepackage{xcolor,colortbl}
\usepackage{setspace}
\usepackage{float}
\usepackage{cite}
\usepackage{graphicx}
\usepackage[space]{grffile}
\usepackage{amsfonts}
\usepackage{algorithm}
\usepackage{algpseudocode}
\usepackage{mathtools}
\usepackage{dsfont}
\usepackage{tikz}%
\providecommand{\U}[1]{\protect\rule{.1in}{.1in}}
\usetikzlibrary{calc,trees,positioning,arrows,chains,shapes.geometric,    decorations.pathreplacing,decorations.pathmorphing,shapes,    matrix,shapes.symbols}
\tikzset{
>=stealth',
  punktchain/.style={
    rectangle, 
    rounded corners, 
        draw=black, very thick,
    text width=15em, 
    minimum height=3em, 
    text centered, 
    on chain},
  line/.style={draw, thick, <-},
  element/.style={
    tape,
    top color=white,
    bottom color=blue!50!black!60!,
    minimum width=8em,
    draw=blue!40!black!90, very thick,
    text width=10em, 
    minimum height=3.5em, 
    text centered, 
    on chain},
  every join/.style={->, thick,shorten >=1pt},
  decoration={brace},
  tuborg/.style={decorate},
  tubnode/.style={midway, right=2pt},
}
\IEEEoverridecommandlockouts
\begin{document}

\title{{\LARGE \textbf{Combined Eco-Routing and Power-Train Control of Plug-In Hybrid
Electric Vehicles in Transportation Networks
}}}
\author{Arian Houshmand$^{1}$, Christos G. Cassandras$^{1}$, Nan Zhou$^{1}$, Nasser Hashemi$^{1}$, Boqi Li$^{2}$, and
Huei Peng$^{2}$ \thanks{*This work was
supported in part by NSF under grants ECCS-1509084, DMS-1664644, CNS-1645681,
by AFOSR under grant FA9550-19-1-0158, by ARPA-E's NEXTCAR program under grant
DEAR0000796 and by the MathWorks.}\thanks{$^{1}$The authors are with Division
of Systems Engineering, Boston University, Brookline, MA 02446 USA
\texttt{{\small arianh@bu.edu; cgc@bu.edu; nanzhou@bu.edu; nhashemi@bu.edu}}} \thanks{$^{2}$The authors are
with Department of Mechanical Engineering, University of Michigan, Ann Arbor,
MI 48109 USA \texttt{{\small boqili@umich.edu; hpeng@umich.edu }}}}
\maketitle

\begin{abstract}
We study the problem of eco-routing for Plug-In Hybrid Electric Vehicles
(PHEVs) to minimize the overall energy consumption cost. We propose an
algorithm which can simultaneously calculate an energy-optimal route
(eco-route) for a PHEV and an optimal power-train control strategy over this
route. In order to show the effectiveness of our method in practice, we use a
HERE Maps API to apply our algorithms based on traffic data in the city of
Boston with more than 110,000 links. Moreover, we validate the performance of
our eco-routing algorithm using speed profiles collected from a traffic
simulator (SUMO) as input to a high-fidelity energy model to calculate energy
consumption costs. Our results show significant energy savings (around 12\%)
for PHEVs with a near real-time execution time for the algorithm.


\end{abstract}

\bstctlcite{MyBSTcontrol}

\thispagestyle{empty} \pagestyle{empty}



\section{INTRODUCTION}

\label{sec: Intro} Due to environmental concerns and the high cost of gas,
there has been an increasing interest in vehicles using alternative energy
sources such as Electric Vehicles (EV). However, given battery capacity levels
in current EVs, their adoption is limited by the All-Electric Range (AER). In
this respect, Plug-In Hybrid Electric Vehicles (PHEVs) offer a suitable
alternative, as they can overcome range limitations by using both gas and
electricity. Depending on battery size, PHEVs can be driven 10-40 miles on
electricity, which is roughly the average daily commuting distance in the US
\cite{collia_2001_2003}. Moreover, it is possible to decrease the energy
consumption cost and the carbon footprint of PHEVs using smart eco-routing and
power-train control strategies.

Traditional vehicle routing algorithms seek to find the minimum
time (fastest) or shortest path routes
\cite{bertsekas_dynamic_1995,braekers_vehicle_2016,toth_vehicle_2002}, whereas
eco-routing algorithms find the paths that minimize the total energy
consumption cost. Several eco-routing algorithms have been studied in the
literature for conventional vehicles that are capable of finding the
energy-optimal routes using historical and online traffic data
\cite{barth_environmentally-friendly_2007,boriboonsomsin_eco-routing_2012,andersen_ecotour:_2013,yao_study_2013,yang_stochastic_2014}%
. Kubicka et al \cite{kubicka_performance_2016} performed a case study to compare
the objective values proposed in the eco-routing literature and showed that
the performance of eco-routing algorithms is highly dependent on the energy model
used to calculate the traveling cost of each link. Pourazarm et al studied optimal routing of electric vehicles considering recharging at charging stations \cite{pourazarm_optimal_2016}. De Nunzio et al \cite{de_nunzio_model-based_2017} studied the eco-routing problem for EVs considering road grade and speed changes on each road link. Although eco-routing of
conventional vehicles is well studied, there is little research that addresses
the case of PHEVs \cite{guanetti_control_2018}.
Jurik et al \cite{cela_energy_2014} studied the problem of eco-routing for HEVs considering the vehicle longitudinal
dynamics. Sun et al \cite{sun_save_2016} and Qiao et al
\cite{qiao_vehicle_2016} proposed the Charge Depleting First (CDF) approach to
address eco-routing for PHEVs. Furthermore, it is shown in
\cite{sun_save_2016} that energy-optimal paths typically take more time
compared to the time-optimal routes. More recently, De Nunzio et al.
\cite{de2018constrained} proposed a semi-analytical solution of the
power-train energy management based on Pontryagin's minimum principle to
address the eco-routing of HEVs, and in \cite{salazar2019optimal} the
eco-routing problem for PHEVs is solved by minimizing a combination of time
and energy.

The contributions of this paper are summarized as follows. Based on the work
introduced in \cite{houshmand_eco-routing_2018}, we first review a CDF
strategy for finding the energy optimal route for PHEVs and propose two
methods for solving this problem: a modified version of Dikjstra's algorithm
\cite{dijkstra_note_1959}, and a Hybrid-LP Relaxation algorithm. We then
propose a Combined Routing and Power-Train Control (CRPTC) eco-routing
algorithm for PHEVs that can simultaneously find an energy optimal route as
well as an optimal power-train control strategy along the route. In contrast
to existing methods in the literature where the power-train control strategy
is considered fixed \cite{sun_save_2016,qiao_vehicle_2016}, we allow the
optimizer to find the optimal PT control strategy. We formulate the problem as
a Mixed Integer Linear Program (MILP) and later relax it into a bi-level
optimization problem where the upper level problem finds the eco-route and the
lower level problem determines the optimal PT switching control strategy
between electricity and gas using a Linear Programming (LP) problem
formulation. We show that the bi-level eco-routing algorithm is
computationally more efficient than the CRPTC approach and its results are
very close to the optimal values calculated using the CRPTC algorithm. Using a
HERE Maps API \cite{heremaps}, we developed
a publicly available web-based tool in which we can request and download the
geographical map of a region alongside its traffic information. Using this platform, we applied our eco-routing
algorithms to large urban traffic networks, including the city of Boston
(110,000 links, 50,000 nodes).
As an alternative to such traffic
data, we also use the Simulation of Urban MObility (SUMO)
\cite{krajzewicz_recent_2012} to investigate traffic outcomes and also collect
speed traces of vehicles following eco-routes and fastest routes. We then use
the Vehicle-Engine SIMulation (VESIM) model
\cite{malikopoulos_simulation_2006}, a high fidelity power-train energy
modelling software package, to calculate the actual energy consumption of
travelling through an eco-route and fastest route for each individual
origin-destination (O-D) pair and compare them against each other. This
approach is used to validate the performance of our eco-routing algorithm with
results suggesting energy savings of about 12\% compared to the fastest route.
We show the trade-off between saving energy and
time in Section \ref{sec:numerical results}. 

The remainder of the paper is organized as follows. The PHEV energy
consumption model is presented in Section \ref{sec: energy model}. A modified
Dijkstra's algorithm as well as MILP problem formulation are proposed in
Section \ref{sec: eco-routing problem} to solve the eco-routing problem. In
Section \ref{sec:numerical results}, we explain our traffic data platform and
by using it, we apply our eco-routing algorithms to the urban area of Boston.
In Section \ref{sec:validation}, using SUMO and VESIM we introduce a framework
to validate the performance of our eco-routing algorithm in real-world
scenarios. Finally, conclusions and further research directions are outlined
in Section \ref{sec: conclusions}.


\section{PHEV Energy Consumption Modeling}

\label{sec: energy model} The first step in developing an efficient
eco-routing algorithm is to understand how the PHEV power-train works, and how
one can model its energy consumption cost. Unlike conventional vehicles where
it is possible to analytically estimate fuel consumption costs as functions of
the velocity and acceleration of the vehicle \cite{kamal_model_2013},
estimating a PHEV's fuel consumption is a more involved process. This is
mainly due to the complexity of the PHEV power-train's architecture. A PHEV
can run on electricity, gas or as a hybrid. Moreover, when PHEVs use
electricity, the battery can be recharged using the regenerative brake and/or
other mechanisms \cite{gao_investigation_1999}. As such, we need a
comprehensive model which takes into account the effect of motor/generator
units and the Internal Combustion Engine (ICE) to calculate the fuel rate and
the electrical power demand from the battery pack.

A PHEV power-train has several different components that work together to
drive the vehicle including the engine, motor/generator, inverter, etc. The
interactions between these components should be considered to estimate the
vehicle's energy consumption. The energy consumption of a PHEV over a finite
time horizon can be expressed as follows:
\begin{equation}
\int_{t_{0}}^{t_{f}}(C_{gas}\dot{m}_{gas}(t)+C_{ele}P_{batt}(t))dt
\label{eqn:cost function}%
\end{equation}
where $\dot{m}_{gas}$ is the instantaneous fuel consumption rate, and
$P_{batt}$ is the total electrical power used/generated by the motor/generator
units. Moreover, $C_{gas}$ ($\$/gallon$) and $C_{ele}$ ($\$/kWh$) are the cost
of gas and electricity, respectively. We discuss two possible approaches to
calculate $\dot{m}_{gas}(t)$ and $P_{batt}(t)$ at any operating condition: a
direct method and an indirect method.

\subsection{Direct Method}

\label{sec: Direct Method} One can calculate the vehicle's energy consumption
at any given time, knowing the details of a vehicle's power-train
architecture, efficiency maps of engine and motor/generator units, and
physical parameters of the vehicle \cite{liu_modeling_2008}. To do so, we need
to have the torque and speed demand from the engine and motor/generators at
any given time. Considering the vehicle's specifications, we can translate its
speed and acceleration to torque and rotational speed demand from the engine
and motor/generator units \cite{moura_tradeoffs_2010}. We can then use these
values to extract $\dot{m}_{gas}(t)$ and $P_{batt}(t)$ from the efficiency
maps. We can either use commercially available software such as
\textit{Autonomie/PSAT }\cite{autonomie}, or develop our own
functions by knowing the details of a specific vehicle. There are two
functions through which we can calculate the fuel and electrical energy
consumptions for any given vehicle as follows:
\begin{align*}
\dot{m}_{gas}(t)  &  =f(v(t),a(t))\\
P_{batt}(t)  &  =g(v(t),a(t))
\end{align*}
where $v(t)$ and $a(t)$ are the speed and acceleration of the vehicle at any
given time respectively. Even though this approach may lead to accurate
estimates of energy consumption values, it is an elaborate method which is not
suitable for the purpose of our higher-level eco-routing framework. Hence, we
use a computationally more efficient approach which we call \textquotedblleft
indirect method\textquotedblright.

\subsection{Indirect Method}

We use a simplified energy model which was first proposed by Qiao et al
\cite{qiao_vehicle_2016} to calculate the energy consumption cost of PHEVs.
Instead of using real time driving data for a targeted vehicle, we calculate
the average $\dot{m}_{gas}$ and $P_{batt}$ per mile for different drive cycles
(Table \ref{tab: Drive cycle assignment}) using a modified version of the
Vehicle-Engine SIMulation (VESIM) model reported in
\cite{malikopoulos_simulation_2006} and references therein. In this method, we
consider two driving modes for a PHEV: Charge-Depleting (CD) and Charge
Sustaining (CS). The CD mode refers to the phase where the PHEV acts like an
EV and consumes all of its propulsion energy from the battery pack. Once the
State Of Charge (SOC) of the battery reaches a target value, it switches to
the CS mode in which the vehicle starts using the internal combustion engine
as the main propulsion system and the battery and electric motors are only
used to improve the fuel economy as in HEVs \cite{karabasoglu_influence_2013}.
\begin{figure}[t]
\centering
\includegraphics[width=0.48\textwidth]{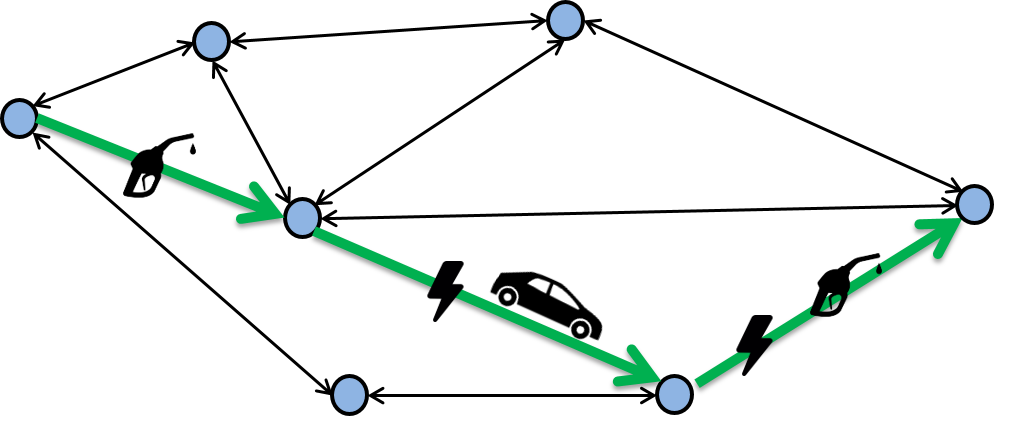}
\caption{Traffic network representation as a directed graph (blue dots represent intersections
and black arcs denote road-links). In green the optimal eco-route is shown
together with the optimal power-train control strategy for switching between
charge depleting and charge sustaining modes.}
\label{fig:eco_route schenatic}%
\end{figure}

Let us consider the traffic network as a directed graph (Fig.
\ref{fig:eco_route schenatic}). Based on their traffic intensity we can
categorize the links into three modes: low, medium, and high traffic intensity
links, and assign different standard drive cycles to them
\cite{qiao_vehicle_2016} (Table \ref{tab: Drive cycle assignment}). For any target vehicle, we can
then use a high-fidelity energy model to calculate
the average electrical energy ($\mu_{CD}$) and gas ($\mu_{CS}$) used to drive
one mile under CD and CS modes respectively under each of theses drive cycles
(Table \ref{tab: conversion factors}):%

\[%
\begin{array}
[c]{lr}%
\mu_{CD_{ij}}=\frac{d_{ij}}{P_{batt_{ij}}}\quad\text{,} & \mu_{CS_{ij}}%
=\frac{d_{ij}}{m_{gas_{ij}}}\quad\text{,}%
\end{array}
\]
where $d_{ij}$ is the length of link $(i,j)$. By knowing $\mu_{CD_{ij}}$ and
$\mu_{CS_{ij}}$ on each link, as well as the network topology (length of each
link), we can determine the average fuel consumption ($m_{gas_{ij}}$) and
electrical power demand from the battery ($P_{batt_{ij}}$) on each link
$(i,j)$. We can then use (\ref{eqn:cost function}) to calculate the total
energy cost for each trip. In this paper we use VESIM as the high-fidelity energy model, and it is calibrated for a PHEV Audi A3 e-tron. Note that this method can easily be extended to conventional vehicles, HEVs and EVs by just considering one of the operational modes. 

\begin{table}[h]
\caption{Drive cycle assignment of each link }%
\label{tab: Drive cycle assignment}
\centering
\resizebox{1\columnwidth}{!}{\renewcommand{\arraystretch}{0.9}
\begin{tabular}{|c|c|}
\hline
\begin{tabular}[c]{@{}c@{}}Traffic Intensity \\ on the link\end{tabular} & Assigned Drive Cycle                     \\ \hline
Low Traffic                                                              & HWFET: Highway Fuel Economy Test         \\ \hline
Medium Traffic                                                           & UDDS: Urban Dynamometer Driving Schedule \\ \hline
High Traffic                                                             & NYC: New York City                       \\ \hline
\end{tabular}
}\end{table}



\begin{table}[h]
\caption{Energy consumption of various driving cycles}%
\label{tab: conversion factors}
\centering
\resizebox{1\columnwidth}{!}{\renewcommand{\arraystretch}{0.9}
\begin{tabular}{|l|l|l|l|l|l|}
\hline
Vehicle Type & Symbol     & Unit     & HWFET & UDDS & NYC  \\ \hline
Audi A3       & $\mu_{CD}$ & $mi/kWh$ & 4.14   & 4.39  & 3.14  \\
& $\mu_{CS}$ & $mi/gal$ & 47.11  & 49.03 & 28.88 \\ \hline
\end{tabular}
}\end{table}

\begin{figure}[h]
\centering
\begin{subfigure}{.16\textwidth}
\centering
\includegraphics[width=0.98\linewidth]{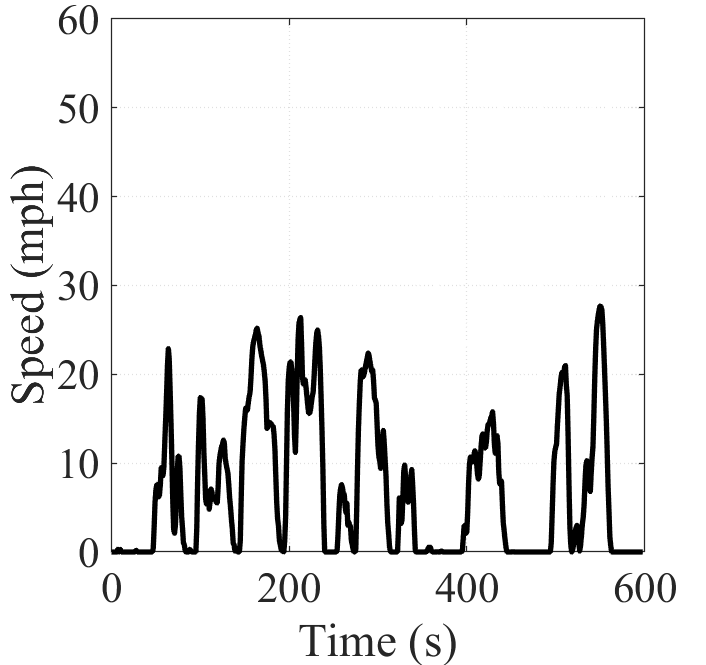}
\caption{}
\label{fig:NYC}
\end{subfigure}\begin{subfigure}{.16\textwidth}
\centering
\includegraphics[width=0.98\linewidth]{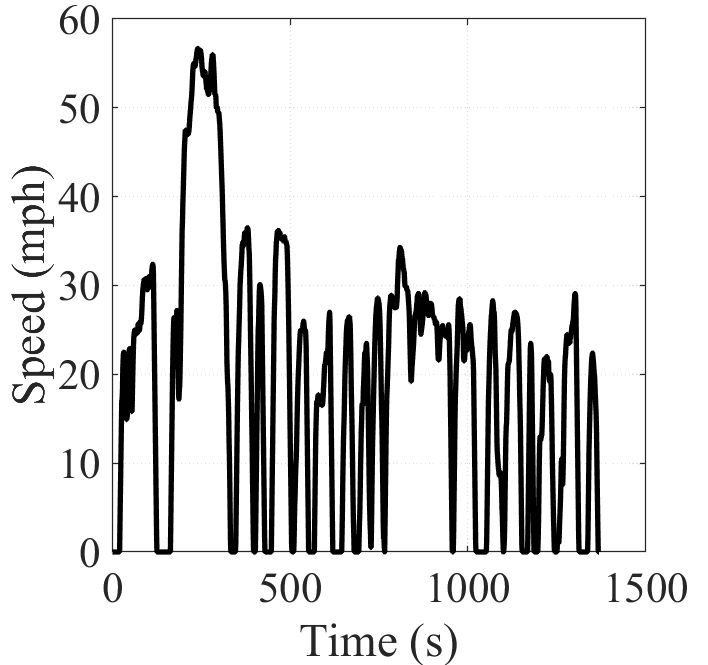}
\caption{}
\label{fig:UDDS}
\end{subfigure}\begin{subfigure}{.16\textwidth}
\centering
\includegraphics[width=0.98\linewidth]{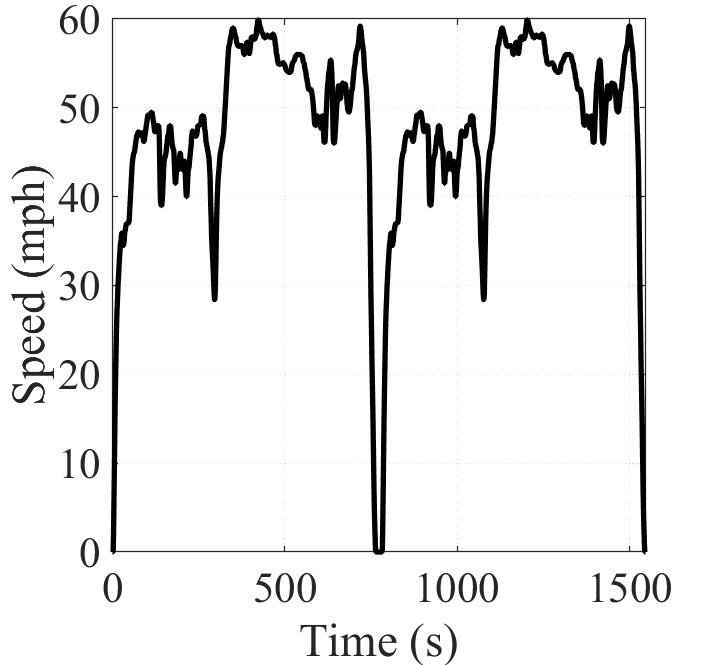}
\caption{}
\label{fig:HWFET}
\end{subfigure}
\caption{ (a) NYC drive cycle: high traffic links; (b) UDDS drive cycle:
medium traffic links; (c) HWFET drive cycle: low traffic links}%
\label{fig:Drive Cycles}%
\end{figure}


\section{Single Vehicle Eco-Routing}

\label{sec: eco-routing problem}

\label{sec: single vehicle routing} In view of the two driving modes (CD and
CS), coming up with an eco-routing algorithm requires knowledge of how the PT
controller switches between the two modes on each link a priori;
alternatively, we can let the algorithm decide the PT control strategy while
finding the optimal route.

In this section, we first review the \textit{Charge Depleting First} (CDF)
eco-routing approach \cite{sun_save_2016,qiao_vehicle_2016} and propose two
methods to solve it. Next, we propose a \textit{Combined Routing and
Power-Train Control} (CRPTC) algorithm to solve the eco-routing problem.

\subsection{Problem Formulation}

We model the traffic network as a directed graph $G=(\mathcal{N},\mathcal{A})$
with $\mathcal{N}={1,...,n}$ and $|\mathcal{A}|=m$ with the arc (link)
connecting node \textit{i} to \textit{j} denoted by $(i,j)\in\mathcal{A}$. The
set of nodes that are incoming/outgoing to node \textit{i} are defined as:
$\mathcal{I}(i)=\{j\in\mathcal{N}|(j,i)\in\mathcal{A}\}$ and $\mathcal{O}%
(i)=\{j\in\mathcal{N}|(i,j)\in\mathcal{A}\}$, respectively. We consider the
single-origin-single-destination eco-routing problem where origin and
destination nodes are denoted by $o$ and $d$ respectively. The energy cost
associated with the vehicle on link $(i,j)$ is denoted by $c_{ij}$. We use
$E_{i}$ to represent the vehicle's residual battery energy at node $i$.
Moreover, we denote the selection of arc $(i,j)$ by $x_{ij}\in\{0,1\}$. The
problem objective is to determine a path from node \textit{o} to \textit{d} so as to minimize the total
energy cost consumed by the vehicle to reach the destination. We consider two
approaches to solve this problem as follows.

\subsection{Charge Depleting First (CDF)}

In this approach \cite{sun_save_2016,qiao_vehicle_2016}, we assume that the
PHEV always starts every trip in the CD mode and uses electricity to drive the
vehicle until it drains all the energy out of the battery pack. Afterwards, it
switches to the CS mode and starts using gas to drive the vehicle. Even though
this approach is generally sub-optimal, it is motivated by the fact that it
eliminates the need for complicated PHEV power-train control strategies to
switch between ICE and electric motors \cite{qiao_vehicle_2016}. As a result,
we can formulate the eco-routing problem using a Mixed-Integer Nonlinear
Programming (MINLP) framework as follows:%

\begin{subequations}
\label{prob: CDF MINLP}
\begin{gather}
\min_{x_{ij},(i,j)\in\mathcal{A}}\sum_{(i,j)\in\mathcal{A}}c_{ij}%
x_{ij}\\%
\label{eqn: CDF cost}
\begin{array}
[c]{lc}%
s.t. & c_{ij}=%
\begin{cases}
C_{gas}\dfrac{d_{ij}}{\mu_{CS_{ij}}}; & E_{i}\leq0\\
C_{ele}\dfrac{d_{ij}}{\mu_{CD_{ij}}}; & E_{i}\geq\frac{d_{ij}}{\mu_{CD_{ij}}%
}\\
C_{ele}E_{i}+C_{gas}\dfrac{d_{ij}-\mu_{CD_{ij}}E_{i}}{\mu_{CS_{ij}}}; &
\text{otherwise}%
\end{cases}
\end{array}
\\%
\begin{array}
[c]{lc}%
E_{j}=\sum\limits_{i\in\mathcal{I}(j)}(E_{i}-\dfrac{d_{ij}}{\mu_{CD_{ij}}%
})x_{ij}, & \forall j\in\mathcal{N}%
\end{array}
\\%
\begin{array}
[c]{lr}%
\smashoperator{\sum_{i:(i,j)\in\mathcal{A}}}x_{ij}+\mathds{1}_{j=o}%
=\smashoperator{\sum_{k:(j,k)\in\mathcal{A}}}x_{jk}+\mathds{1}_{j=d} & \forall
j\in\mathcal{N}%
\end{array}
\label{eqn: CDF-cons_flow}\\%
\label{eqn: CDF x}
\begin{array}
[c]{lr}%
x_{ij}\in\{0,1\} & \forall(i,j)\in\mathcal{A}%
\end{array}
\end{gather}
\end{subequations}
where $C_{gas}$ and $C_{ele}$ are the price of gas ($\$/gallon$) and
electricity ($\$/kWh$) respectively, and $E_{i}$ is the remaining electrical
energy at node \textit{i}. The conversion factors $\mu_{CD_{ij}}$ and
$\mu_{CS_{ij}}$, taken from Table \ref{tab: conversion factors}, are functions
of the traffic intensity on each link $(i,j)$. Note that
(\ref{eqn: CDF-cons_flow}) is the flow conservation constraint
\cite{bertsimas_introduction_1997}, and $\mathds{1}_{\{.\}}$ is a Boolean
indicator function. We assume that the vehicle has enough gas and electrical
power to complete the trip and that $E_{o}\geq0$ (initial energy at the origin).

The nonlinearities in problem (\ref{prob: CDF MINLP}) arise from the
dependency of $c_{ij}$ in (\ref{eqn: CDF cost}) on the remaining energy in the
battery, which means $c_{ij}$ is route-dependent. If $c_{ij}$ were a priori
known, the problem would have been reduced to a shortest path problem and we
could solve it using one of the highly studied algorithms such as Dijktra's
algorithm \cite{dijkstra_note_1959}. To address this issue, in what what
follows we propose two methods to solve this problem: CDF-Dijkstra and
Hybrid-LP relaxation. 

\begin{figure}[ptb]
\centering
\includegraphics[width=0.45\textwidth]{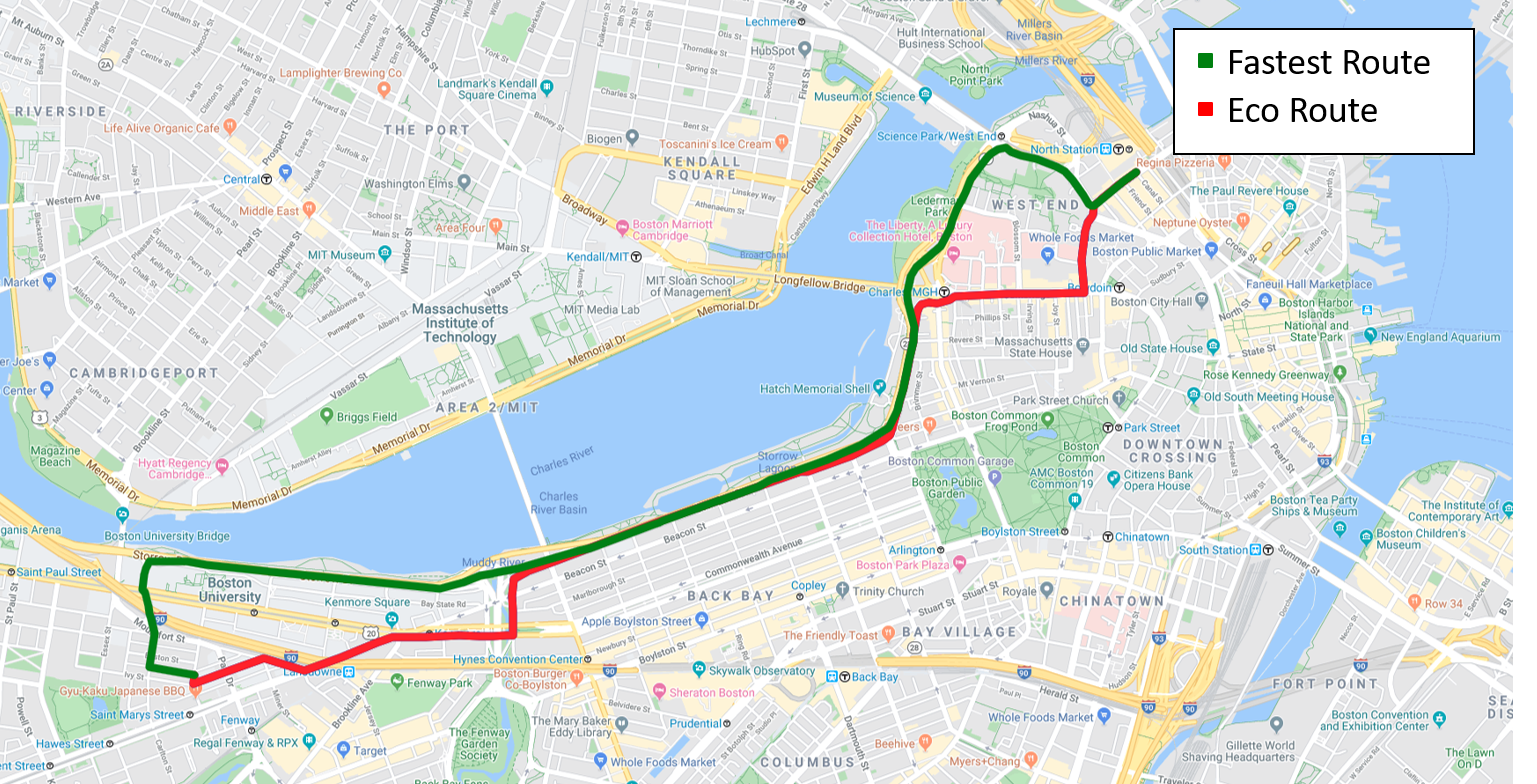}\caption{{Visualization
of eco-route (CRPTC) and fastest route on a map}}%
\label{fig:sample route}%
\end{figure}

\subsubsection{Charge Depleting First (CDF) - CDF- Dijkstra}

Dijkstra's algorithm finds the shortest path between two points if the weight
of each link is a known value. In (\ref{prob: CDF MINLP}) the energy cost of
each link (\ref{eqn: CDF cost}) is a function of remaining energy in the
battery pack. As a result, if we keep track of the remaining energy in the
battery at the end node of each candidate path while searching the graph using
Dijkstra's algorithm, the energy cost of links would be known to us and we can
find the eco-route using a modified version of Dijkstra's algorithm. This
modified version, which can find the CDF eco-route by keeping track of
remaining energy in the battery, is given in Algorithm
\ref{algo: CDF-Dijkstra} (see Appendix).

\subsubsection{Charge Depleting First (CDF) - A Hybrid-LP Relaxation Approach}

We propose an alternative solution to this problem by reducing the MINLP
problem (\ref{prob: CDF MINLP}) to a simpler one which can be solved using a
combination of Linear Programming (LP) and a simple dynamic programming-like
algorithm, in order to guarantee global convergence. The nonlinearities of the
problem arise in (\ref{eqn: CDF cost}) where $c_{ij}$ is a function of
$x_{ij}$. We show that we can reduce this piece-wise constant function to a
constant function, and the MINLP can be converted to a LP by using the
properties of the minimum cost flow problem \cite{hillier_introduction_2012}.
The proposed algorithm is as follows:

\begin{enumerate}
\item Find the shortest path and calculate the energy cost on this path using (\ref{eqn: CDF cost}) and set it to
$\rho$. We are going to use $\rho$ as a reference in the next step.

\item From the origin, construct all paths reaching node $p$ such that
$E_{p}\leq0$ for the
first time and stop constructing the path at this node. Disregard the paths
with a total energy cost greater than $\rho$ and save the remaining paths in a matrix.

\item Since $E_p \leq 0$, we can only use gas (CS mode) to travel through links which belong to the paths outgoing from node $p$.
\item Assuming knowledge of traffic modes on each link, the least energy cost
path from node $p$ in step 2 to the destination node can be found from:
\begin{small}
\begin{subequations}
\label{prob: hybrid LP}%
\begin{gather}
\min_{x_{ij},(i,j)\in\mathcal{A}}\sum_{(i,j)\in\mathcal{A}}c_{ij}%
x_{ij}\label{eqn: hybrid LP}\\%
\begin{array}
[c]{lc}%
s.t. & c_{ij}=C_{gas}\dfrac{d_{ij}}{\mu_{CS_{ij}}}%
\end{array}
\\%
\label{eqn: hybrid LP cons. flow}
\begin{array}
[c]{lr}%
\smashoperator{\sum_{i:(i,j)\in\mathcal{A}}}x_{ij}+\mathds{1}_{j=p}%
=\smashoperator{\sum_{k:(j,k)\in\mathcal{A}}}x_{jk}+\mathds{1}_{j=d} & \forall
j\in\mathcal{N}%
\end{array}
\\%
\begin{array}
[c]{lr}%
x_{ij}\in\{0,1\} & \forall(i,j)\in\mathcal{A}%
\end{array}
\end{gather}
\end{subequations}
\end{small}
Note that constraint (\ref{eqn: hybrid LP cons. flow}) ensures that by solving (\ref{prob: hybrid LP}), we are finding the optimal
path from $p$ to $d$.

\item Using the property of the minimum cost flow problem
\cite{hillier_introduction_2012}, problem (\ref{eqn: hybrid LP}) is equivalent
to an LP problem with the integer restriction of $x_{ij}$ relaxed: $0\leq
x_{ij}\leq1$.

\item Find the path from node $o$ to $p$ with the least energy cost. By the
principle of optimally, the optimal path from \textit{o} to $d$ is the one determined
in this manner followed by the path selected by steps 4 and 5 from node $p$ to
$d$.

\item Find the paths in step 6 for all nodes $p$ such that $E_{p}\leq0$, then
choose the one with the minimum energy cost. The selected path is the minimum
energy cost path.

\item If there are paths without any node such that $E_{p}\leq0$ (generated at
step 2), compare their cost function values with the cost functions in step 6.
The optimal route is the minimum among them.
\end{enumerate}

\subsection{Combined Routing and Power-Train Control (CRPTC)}\label{sec: CRPTC} 

Based on Table \ref{tab: conversion factors}, the CD mode
has the best efficiency on medium traffic links. As such, if we always
consider using the CD mode at the beginning of each trip and then switch to
the CS mode when we run out of battery, we miss the opportunity to harness the
effectiveness of the CD mode on medium traffic links towards the end of the
route. With this motivation, we propose a new algorithm which finds the energy-optimal
routing decisions as well as the PT controller decision to switch between CD
and CS modes. Let $y_{ij}\in[0,1]$ be an additional decision variable on
link $(i,j)$ which represents the fraction of the link's length over which we
use the CD mode (thus, if we only use the CD mode over link $(i,j)$, then
$y_{ij}=1$). Considering the new decision variable, we can formulate the CRPTC
problem as follows:

\begin{small}
\begin{subequations}
\label{prob:CRPTC minlp}%
\begin{gather}
\min_{x_{ij},y_{ij},(i,j)\in\mathcal{A}}\sum_{(i,j)\in\mathcal{A}}%
[c_{gas}\dfrac{d_{ij}}{\mu_{CS_{ij}}}(1-y_{ij})+c_{ele}\dfrac{d_{ij}}%
{\mu_{CD_{ij}}}y_{ij}]x_{ij}\label{eqn:CRPTC Obj}\\%
\begin{array}
[c]{lr}%
s.t.\smashoperator{\sum_{i:(i,j)\in\mathcal{A}}}x_{ij}+\mathds{1}_{j=o}%
=\smashoperator{\sum_{k:(j,k)\in\mathcal{A}}}x_{jk}+\mathds{1}_{j=d} & \forall
j\in\mathcal{N}\\
&
\end{array}
\\
\label{eqn:CRPTC energy cons}%
\sum_{(i,j)\in\mathcal{A}}\dfrac{d_{ij}}{\mu_{CD_{ij}}}y_{ij}x_{ij}\leq
E_{o}\\%
\begin{array}
[c]{lr}%
x_{ij}\in\{0,1\} & \forall(i,j)\in\mathcal{A}%
\end{array}
\\%
\begin{array}
[c]{lr}%
y_{ij}\in\lbrack0,1] & \forall(i,j)\in\mathcal{A}%
\end{array}
\end{gather}
\end{subequations}
\end{small}
 Note that constraint (\ref{eqn:CRPTC energy cons}) ensures that the total
electrical energy used in the CD mode would be less than the initial available
energy in the battery ($E_{o}$). Since we have the term $x_{ij}y_{ij}$ in the
problem formulation, this is a MINLP and we may not be able to determine a
global optimum. Hence, we transform (\ref{prob:CRPTC minlp}) into a MILP by
introducing an intermediate decision variable $z_{ij}=x_{ij}y_{ij}$. We can
then make use of the following inequalities to transform the existing MINLP
problem into a MILP problem:
\begin{equation}
z_{ij}\geq0,\text{ \ }z_{ij}\leq y_{ij},\text{ \ }z_{ij}\geq y_{ij}%
-(1-x_{ij}),\text{ \ }z_{ij}\leq x_{ij} \label{eqn:Zij inequalities}%
\end{equation}
Considering $z_{ij}$ and (\ref{eqn:Zij inequalities}), we can reformulate
problem (\ref{prob:CRPTC minlp}) as follows:
\begin{small}
\begin{subequations}
\label{prob: CRPTC MILP}
\begin{gather}
\min_{x_{ij},y_{ij},z_{ij},(i,j)\in\mathcal{A}}\sum_{(i,j)\in\mathcal{A}%
}(c_{gas}\dfrac{d_{ij}}{\mu_{CS_{ij}}}x_{ij}+(c_{ele}\dfrac{d_{ij}}%
{\mu_{CD_{ij}}}-c_{gas}\dfrac{d_{ij}}{\mu_{CS_{ij}}})z_{ij}%
)\label{eqn:Combined Single Vehicle energy}\\%
\label{eqn: CRPTC flow cons}
\begin{array}
[c]{llr}%
s.t. & \smashoperator{\sum_{i:(i,j)\in\mathcal{A}}}x_{ij}+\mathds{1}_{j=o}%
=\smashoperator{\sum_{k:(j,k)\in\mathcal{A}}}x_{jk}+\mathds{1}_{j=d} & \forall
j\in\mathcal{N}\\
&  &
\end{array}
\\
\label{eqn: CRPTC energy constraint}
\sum_{(i,j)\in\mathcal{A}}\dfrac{d_{ij}}{\mu_{CD_{ij}}}z_{ij}\leq E_{1}\\
z_{ij}\geq0\\
z_{ij}\leq y_{ij}\\
z_{ij}\geq y_{ij}-(1-x_{ij})\\
z_{ij}\leq x_{ij}\\%
\begin{array}
[c]{cc}%
x_{ij}\in\{0,1\} & \forall(i,j)\in\mathcal{N}%
\end{array}
\\%
\label{eqn: CRPTC y}
\begin{array}
[c]{cc}%
y_{ij}\in\lbrack0,1] & \forall(i,j)\in\mathcal{N}%
\end{array}
\end{gather}
\end{subequations}
\end{small}
 This is a MILP problem which can be solved to determine a global optimum.


\subsection{Combined Routing and Power-Train Control: A bi-level approach}

Since MILP problems, such as the one above, are typically NP-hard, we now
investigate a bi-level optimization approach in which the upper-level problem
finds the energy-optimal route considering the CDF approach and the
lower-level problem calculates the optimal PT control strategy by solving an
LP problem. As a result, we can formulate the bi-level optimization problem as follows:
\begin{enumerate}
\item Using Algorithm \ref{algo: CDF-Dijkstra}, solve problem (\ref{prob: CDF MINLP}) and find the optimal routing decision vector
$\mathbf{x}^{\ast}=[x_{ij},(i,j)\in\mathcal{A}]$. Note that for finding the
eco-route, we solve the CDF problem which is the baseline for PT control strategy.

\item Fix the routing decision vector $\mathbf{x}^{\ast}$ calculated in step 1
using CDF and find the optimal switching strategy $\mathbf{y}^{\ast}%
=[y_{ij}^{\ast},(i,j)\in\mathcal{A}]$ by solving the following LP problem:
\begin{subequations}
\label{prob: Bi-Level L2}
\begin{gather}
\min_{y_{ij},(i,j)\in\mathcal{A}}\sum_{(i,j)\in\mathcal{A}}[c_{gas}%
\dfrac{d_{ij}}{\mu_{CS_{ij}}}(1-y_{ij})+c_{ele}\dfrac{d_{ij}}{\mu_{CD_{ij}}%
}y_{ij}]x_{ij}\\
\sum_{(i,j)\in\mathcal{A}}\dfrac{d_{ij}}{\mu_{CD_{ij}}}y_{ij}x_{ij}\leq
E_{1}\\%
\label{eqn: bi-level const}
\begin{array}
[c]{cc}%
x_{ij}=x_{i,j}^{\ast} & \forall(i,j)\in\mathcal{A}%
\end{array}
\\%
\begin{array}
[c]{cc}%
y_{ij}\in\lbrack0,1] & \forall(i,j)\in\mathcal{A}%
\end{array}
\end{gather}
\end{subequations}
 Note that constraint (\ref{eqn: bi-level const}) enforces routing
variables to be equal to the ones calculated in step 1.

\item The optimal route is $\mathbf{x}^{\ast}$ calculated in step 1 and the
optimal PT switching strategy is $\mathbf{y}^{\ast}$ found by solving problem
(\ref{prob: Bi-Level L2}).
\end{enumerate}

Note that this solution is sub-optimal, but the computational time is orders
of magnitudes faster than solving problem (\ref{prob: CRPTC MILP}). This is due
to the fact that the upper-level problem finds the optimal route using the
computationally efficient CDF-Dijkstra algorithm, while the lower-level problem
solves an LP on a small set of decision variables (only links selected by the
router in step 1). We will further discuss the execution time of each
algorithm in subsequent sections.

\begin{figure}[h]
\caption{Bi-level eco-routing structure}%
\centering\begin{tikzpicture}
[node distance=.8cm,
start chain=going below,]
\node[punktchain, join] (level-1) {Find eco-route ($\textbf{x}^*$) by solving problem (\ref{prob: CDF MINLP}) using CDF-Dijkstra Algorithm (\ref{algo: CDF-Dijkstra})};
\node[punktchain, join] (level-2 )      {Solve Problem (\ref{prob: Bi-Level L2}) to find the optimal PT control strategy ($\textbf{y}^*$) over the eco-route };
\end{tikzpicture}\end{figure}
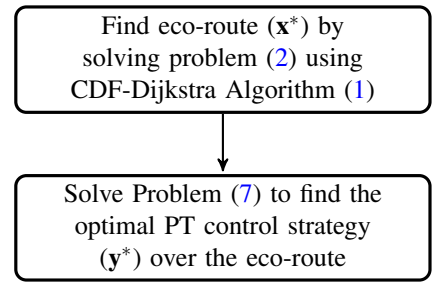

\subsection{Eco-routing with Time Consideration}

There is typically a trade-off between energy savings and time savings in
choosing between eco-routes and fastest routes and, as shown in
\cite{houshmand_eco-routing_2018}, eco-routes can be up to 20\% slower than
fastest routes (Fig. \ref{fig:sample route}). As a result, in order to consider a balance between time and energy
we can introduce a time component into the objective function of our
eco-routing problem. Considering time in the eco-routing problem, we can
re-write problem (\ref{prob: CDF MINLP}) as follows: 
\begin{equation}
\begin{aligned} \min_{x_{ij}, (i,j)\in \mathcal{A}}&\sum_{(i,j)\in\mathcal{A}}(\alpha\frac{t_{ij}}{\beta_1}  + (1-\alpha)\frac{c_{ij}}{\beta_2})x_{ij}\\ \text{s.t. } & (\ref{eqn: CDF cost})-(\ref{eqn: CDF x}) \end{aligned} \label{prob: CDF with time}%
\end{equation}
where $\alpha\in\lbrack0,1]$ is a time-to-energy weighting factor, and
$\beta_{1}>0$ and $\beta_{2}>0$ are normalization factors for time and energy
respectively. Note that if we select $\beta_1 = t_{ij}^{max}$ and $\beta_2 = c_{ij}^{max}$, then the two terms are ensured to be in [0,1]. The two max constants are upper bounds selected based on the topology of the network and the pricing structure.
We can use the same analogy and modify the CRPTC problem
(\ref{prob: CRPTC MILP}) to include time as follows: 
\begin{equation}
\begin{aligned} \min_{x_{ij},y_{ij},z_{ij},(i,j)\in\mathcal{A}}&\sum_{(i,j)\in\mathcal{A}}(\frac{1-\alpha}{\beta_2})(c_{gas}\dfrac{d_{ij}}{\mu_{CS_{ij}}}x_{ij}+ \\ & (c_{ele}\dfrac{d_{ij}}{\mu_{CD_{ij}}}-c_{gas}\dfrac{d_{ij}}{\mu_{CS_{ij}}})z_{ij})+\frac{\alpha}{\beta_1} t_{ij}x_{ij} \\ \text{s.t. } & (\ref{eqn: CRPTC flow cons})-(\ref{eqn: CRPTC y}) \end{aligned} \label{prob: CRPTC with time}%
\end{equation}


\section{Numerical results}

\label{sec:numerical results} An eco-routing case study was presented in
\cite{houshmand_eco-routing_2018} using traffic data from the Eastern
Massachusetts highway sub-network collected by \textit{INRIX}
\cite{zhang_price_2018-1,zhang_data-driven_2017}. In this paper, we extend our
analysis to a larger network which includes the entire Boston urban area.

\subsection{Traffic Data Platform}

\begin{figure}[ptb]
\centering
\begin{subfigure}{.48\textwidth}
\centering
\includegraphics[width=0.98\linewidth]{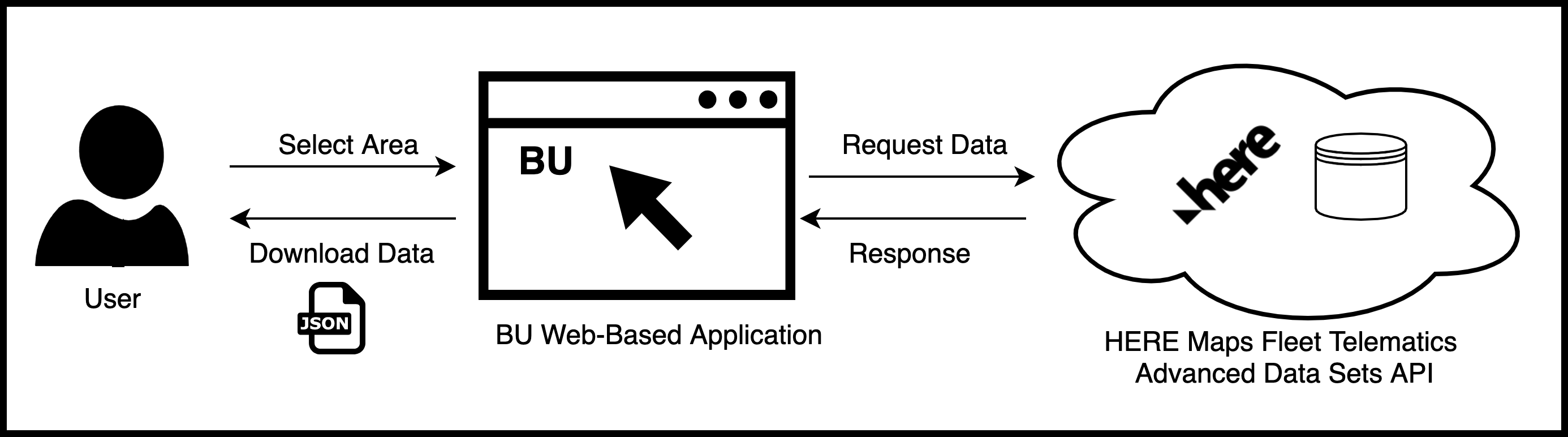}
\caption{}
\label{fig:HERE Architechture}
\end{subfigure}
\begin{subfigure}{.24\textwidth}
\centering
\includegraphics[width=0.98\linewidth]{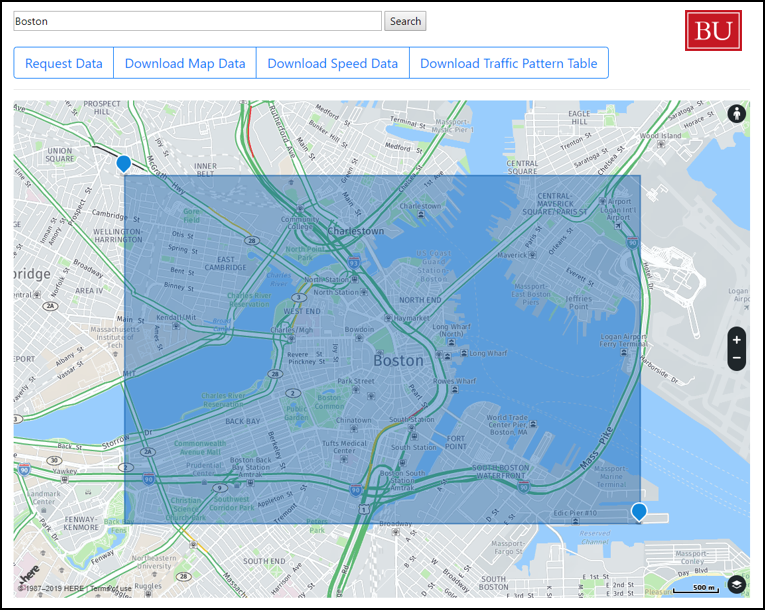}
\caption{}
\label{fig:UI dl}
\end{subfigure}\begin{subfigure}{.24\textwidth}
\centering
\includegraphics[width=0.98\linewidth]{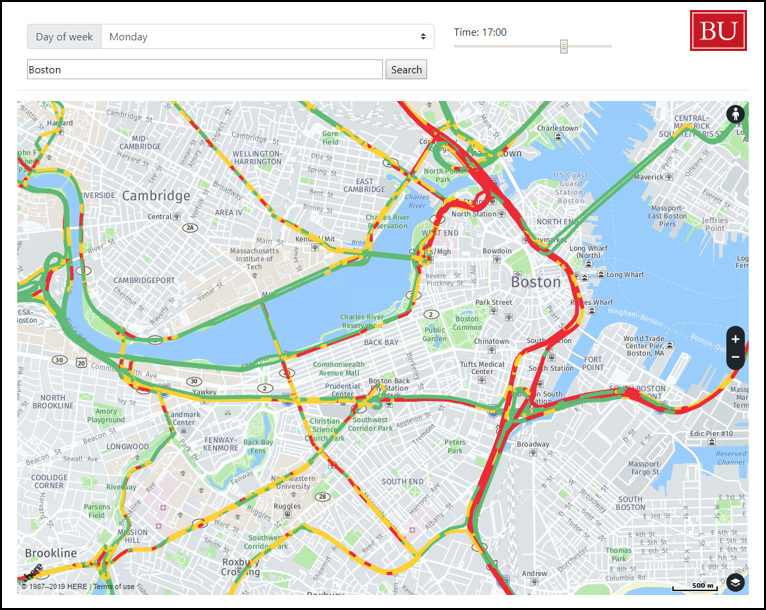}
\caption{}
\label{fig:UI vis}
\end{subfigure}
\caption{ (a) BU traffic data platform architecture; (b) BU traffic data
platform user interface: Downloading tool; (c) Platform user interface:
Average speed visualization}%
\label{fig:BU HERE Platform}%
\end{figure}

Using a HERE Maps API, we developed a web-based tool in which we can request
and download the geographical map of a region alongside its traffic
information (free flow and average speeds)
(\href{http://www.bu.edu/codes/simulations/traffic_downloader/}{http://www.bu.edu/codes/simulations/traffic\_downloader}%
). The traffic data includes average speeds of all roads for every 15 minutes
of a typical week and their free flow speeds. Moreover, we have topological
data of the selected region including how the links are connected to each
other, positions of nodes and links, link lengths, road grade, etc. The
structure of this platform is shown in Fig. \ref{fig:HERE Architechture} .

\begin{figure}[ptb]
\centering
\includegraphics[width=0.48\textwidth]{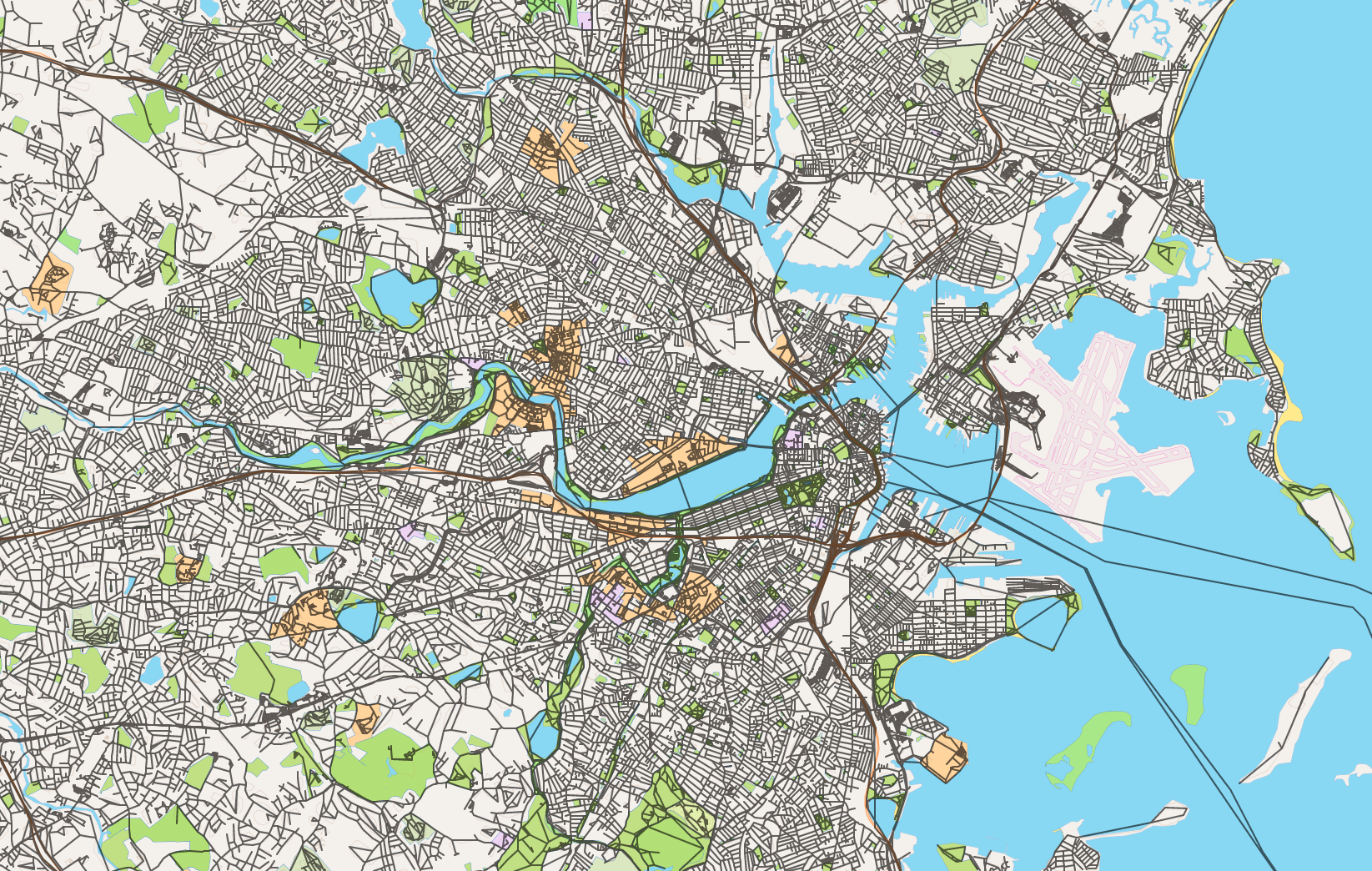}\caption{{Available
links in our case study network of Boston}}%
\label{fig:HERE Boston}%
\end{figure}

\subsection{Data Preprocessing}

Using this platform, we downloaded traffic data of the Boston urban area (Fig.
\ref{fig:HERE Boston}) which includes more than 110,000 links and 50,000
nodes. Since the energy model in our eco-routing algorithm depends on the
traffic intensity of each link (Table \ref{tab: conversion factors}), we need
to first categorize links in terms of traffic intensity: low, medium, high.
Using the average speed of every link and its corresponding free flow speed
data, we introduce a new variable called \emph{speed factor} ($S_{ij}$) for
each link $(i,j)$ as follows:
\[
S_{ij}=\frac{\bar{v}_{ij}}{f_{ij}}%
\]
where $\bar{v}_{ij}$ and $f_{ij}$ are the average speed and free flow speed
respectively of link $(i,j)$. Note that the speed factor is a normalized speed
value on each link indicating how congested the link is. As a result, we can
use $S_{ij}$ to categorize links into three modes: $(i)$ If $S_{ij}\leq0.5$,
link $(i,j)$ is categorized as a high traffic intensity link and values of the
NYC drive cycle are assigned to it (Table \ref{tab: conversion factors}),
$(ii)$ If $0.5<S_{ij}<0.75$, link $(i,j)$ is categorized as a medium traffic
intensity link and values of UDDS are assigned to it, $(iii)$ If $S_{ij}%
\geq0.75$, link $(i,j)$ is categorized as a low traffic intensity link and
values of HWFET are assigned to it.

\subsection{Performance Measurement Baseline}

In order to measure the performance of our eco-routing algorithms, we consider
the time-optimal path (fastest-route) as the baseline. In this respect, for
each O-D pair we find both the eco-route and fastest route and then compare
the energy cost of travelling through both. The fastest route can be
determined by solving the following problem: 
\begin{subequations}
\label{eqn: travel time}%
\begin{gather}
\min_{x_{ij},(i,j)\in\mathcal{A}}\sum_{(i,j)\in\mathcal{A}}t_{ij}x_{ij}\\%
\begin{array}
[c]{lr}%
s.t. & (\ref{eqn: CDF-cons_flow})%
\end{array}
\\
t_{ij}=\dfrac{d_{ij}}{\bar{v}_{ij}}\\
x_{ij}\in\{0,1\}
\end{gather}
where $t_{ij}$ and $d_{ij}$ are the traveling time and length of link $(i,j)$
respectively. In order to calculate the energy cost of travelling through the
fastest route, we consider the CDF policy for the car and calculate its energy
cost using (\ref{eqn: CDF cost}).

\subsection{Energy Cost and Travel Time Comparison Results}

We use the urban area of Boston (Fig. \ref{fig:HERE Boston}) as our case-study
network. In order to show the impact of traffic intensity and distance between
O-D pairs on the performance of the eco-routing algorithms, we randomly select
100 O-D pairs in this network and calculate eco-routes (CDF, CRPTC, and
Bi-level) as well as the fastest route between each of these O-D pairs over
different hours of a day (8:00 am, 12:00 pm, 3:00 pm, 5:00 pm, and 9:00 pm).
As in \cite{qiao_vehicle_2016}, we assume that the initial available energy in
the battery is $E_{ini}=5.57kWh$ and the cost of gas and electricity are
$C_{gas}=2.75\$/gal$ and $C_{ele}=0.114\$/kWh$ respectively
\cite{karabasoglu_influence_2013}. Since the amount of allowable electrical
energy depletion $\Delta E$ affects the efficiency of eco-routing algorithms,
five different values for $\Delta E$ are selected: 0\textit{kWh},
0.5\textit{kWh}, 1\textit{kWh}, 2.5\textit{kWh}, and 5.7\textit{kWh}. The
average energy and time comparison plots over the selected O-D pairs for
different allowable $\Delta E$ values are shown in Figs.
\ref{fig:Energy costs} and \ref{fig:Travel times} respectively.
\begin{figure}[h]
\centering
\begin{subfigure}{.24\textwidth}
\centering
\includegraphics[width=0.98\linewidth]{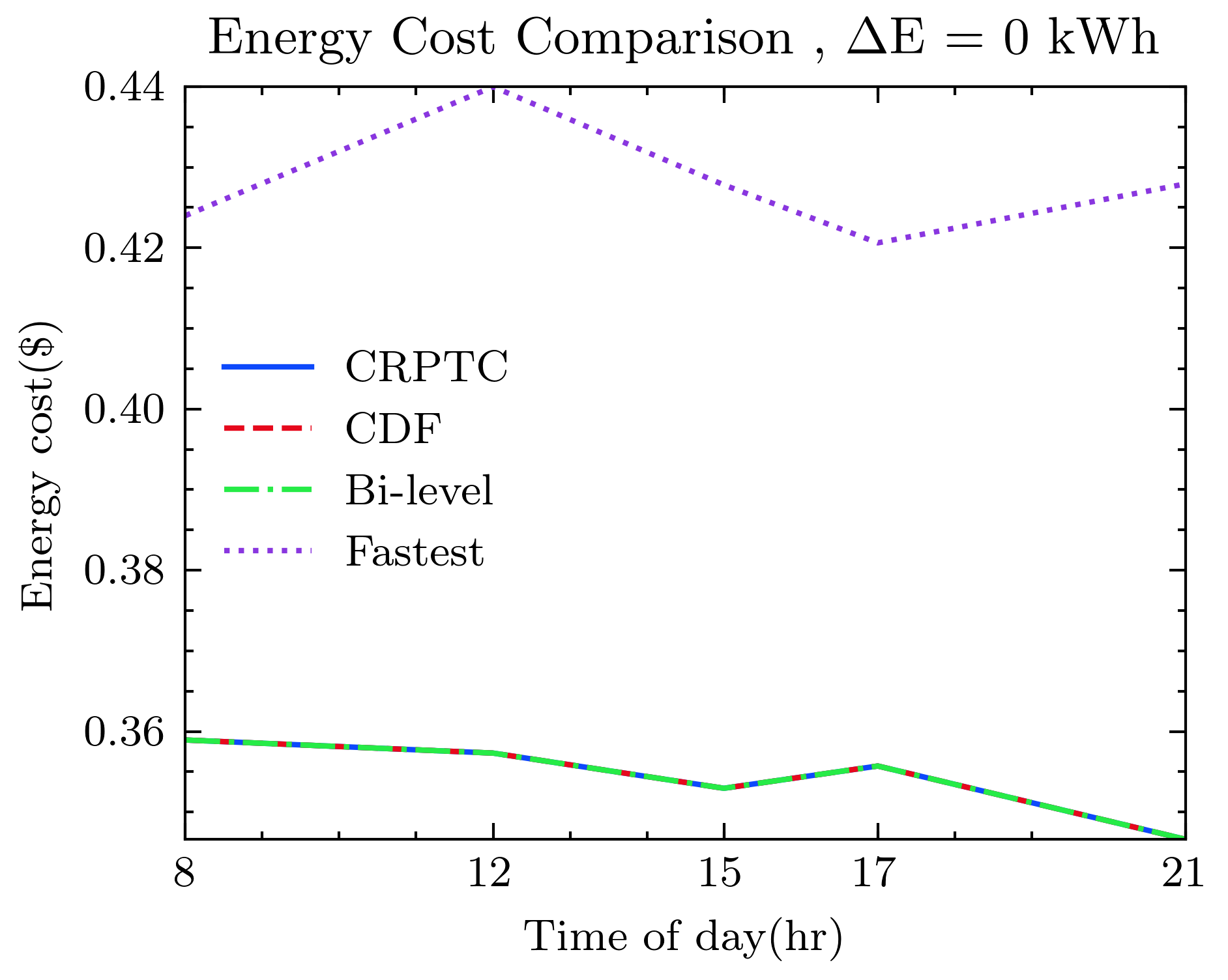}
\caption{}
\label{fig:energy comp E0}
\end{subfigure}\begin{subfigure}{.24\textwidth}
\centering
\includegraphics[width=0.98\linewidth]{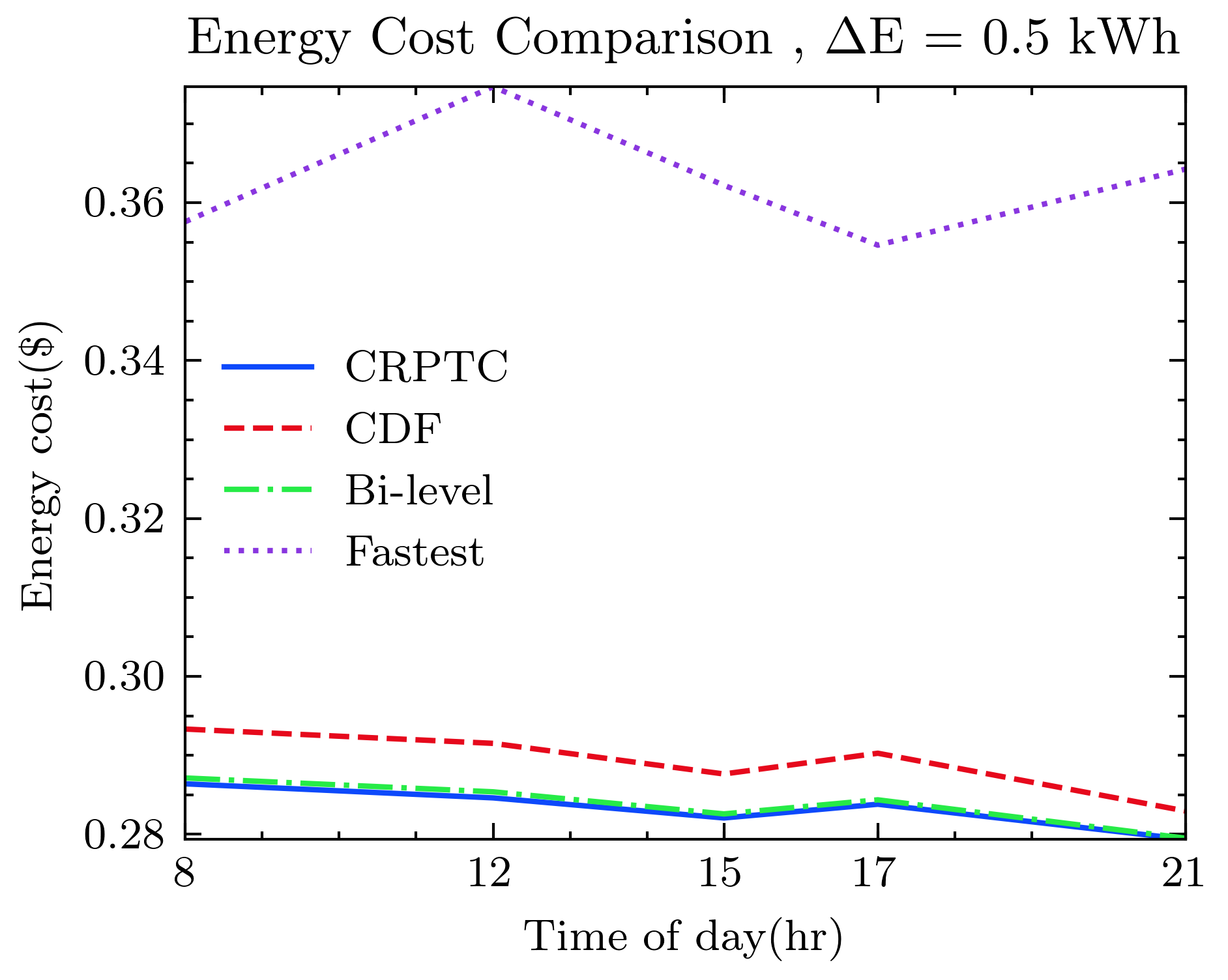}
\caption{}
\label{fig:energy comp E0.5}
\end{subfigure}
\begin{subfigure}{.24\textwidth}
\centering
\includegraphics[width=0.98\linewidth]{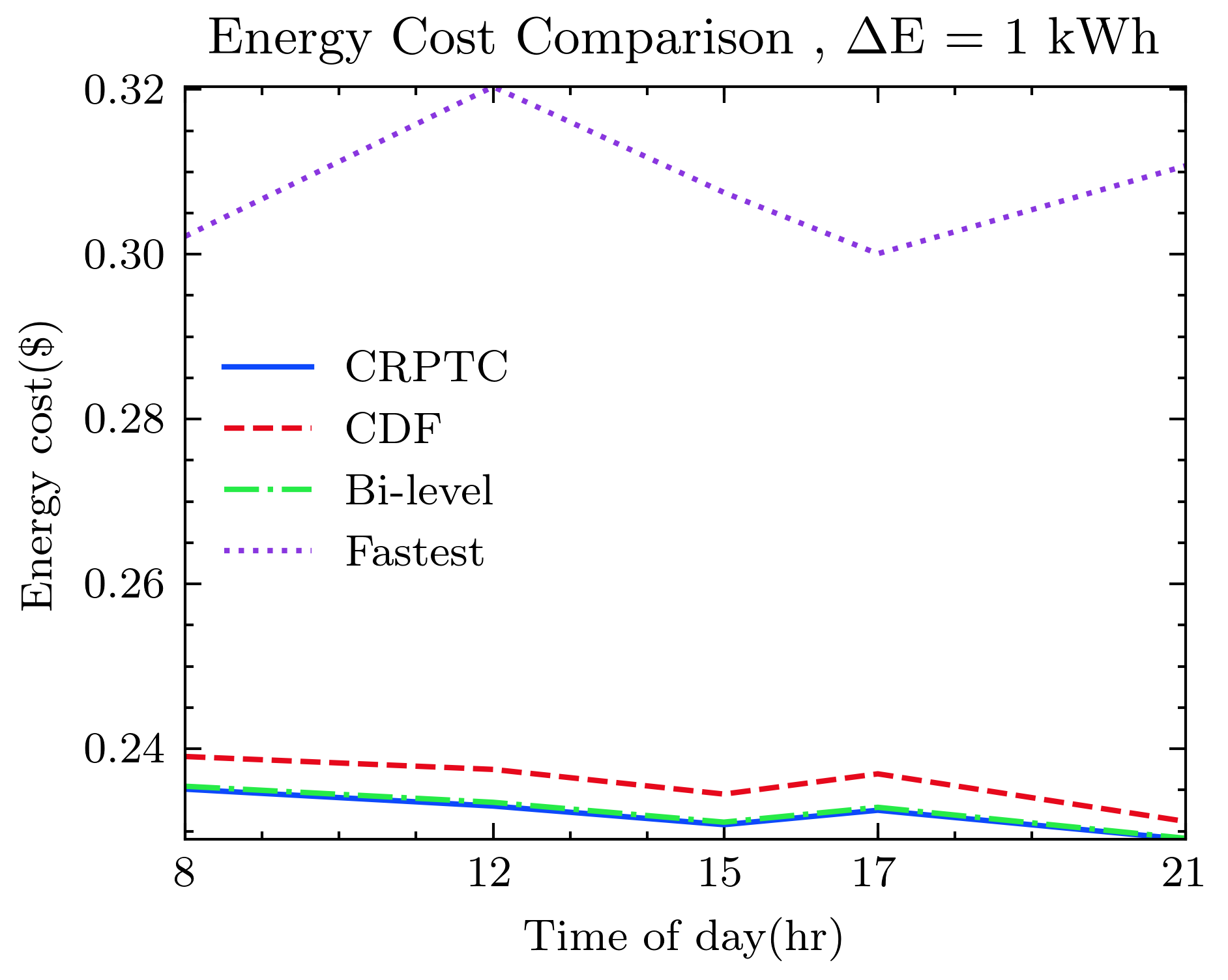}
\caption{}
\label{fig:energy comp E1}
\end{subfigure}%
\begin{subfigure}{.24\textwidth}
\centering
\includegraphics[width=0.98\linewidth]{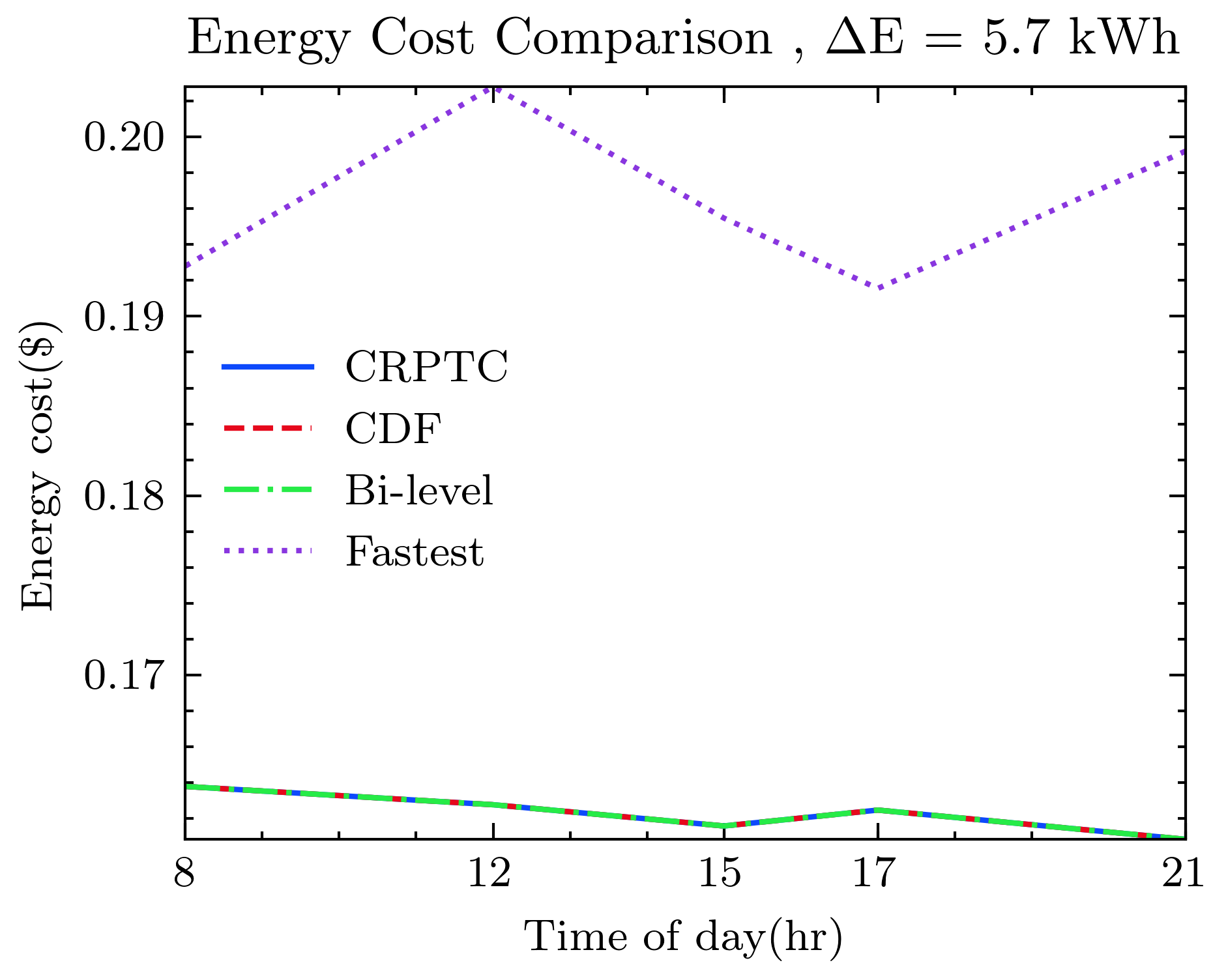}
\caption{}
\label{fig:energy comp E5.7}
\end{subfigure}
\caption{ Average energy cost comparison for the selected O-D pairs with different allowable $\Delta E$ values}
\label{fig:Energy costs}%
\end{figure}

\begin{figure}[h]
\centering
\begin{subfigure}{.24\textwidth}
\centering
\includegraphics[width=0.98\linewidth]{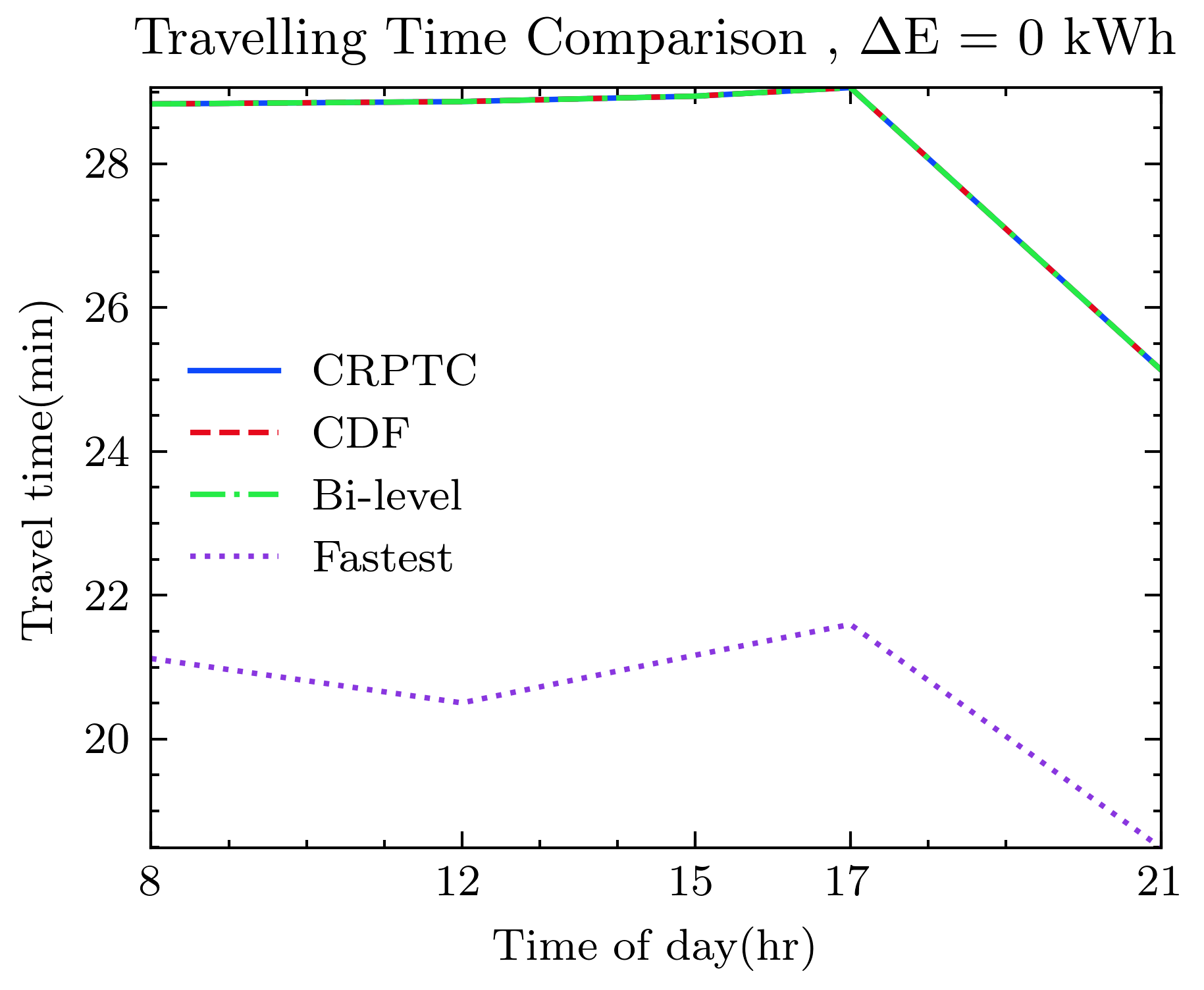}
\caption{}
\label{fig:time comp E0}
\end{subfigure}\begin{subfigure}{.24\textwidth}
\centering
\includegraphics[width=0.98\linewidth]{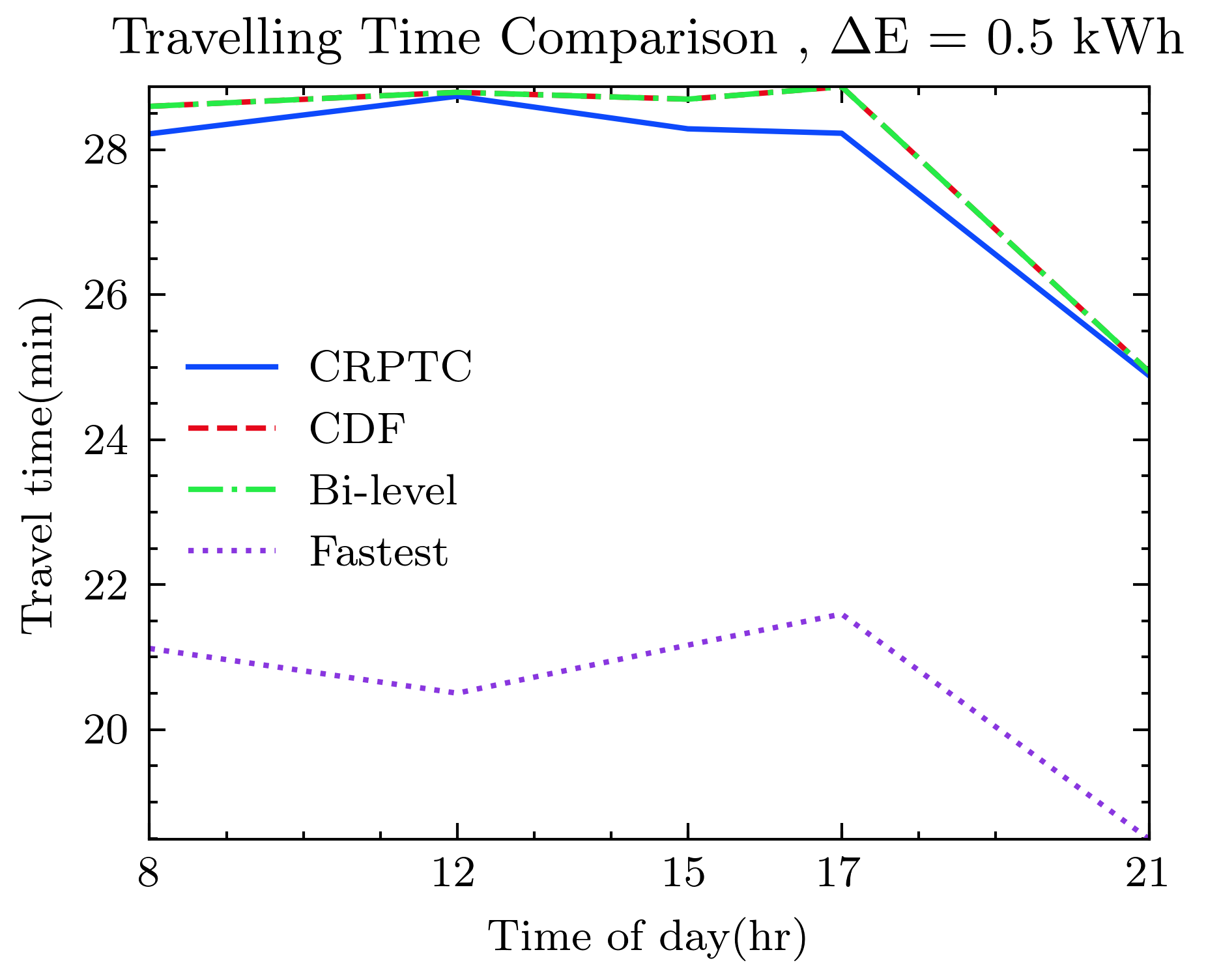}
\caption{}
\label{fig:time comp E0.5}
\end{subfigure}
\begin{subfigure}{.24\textwidth}
\centering
\includegraphics[width=0.98\linewidth]{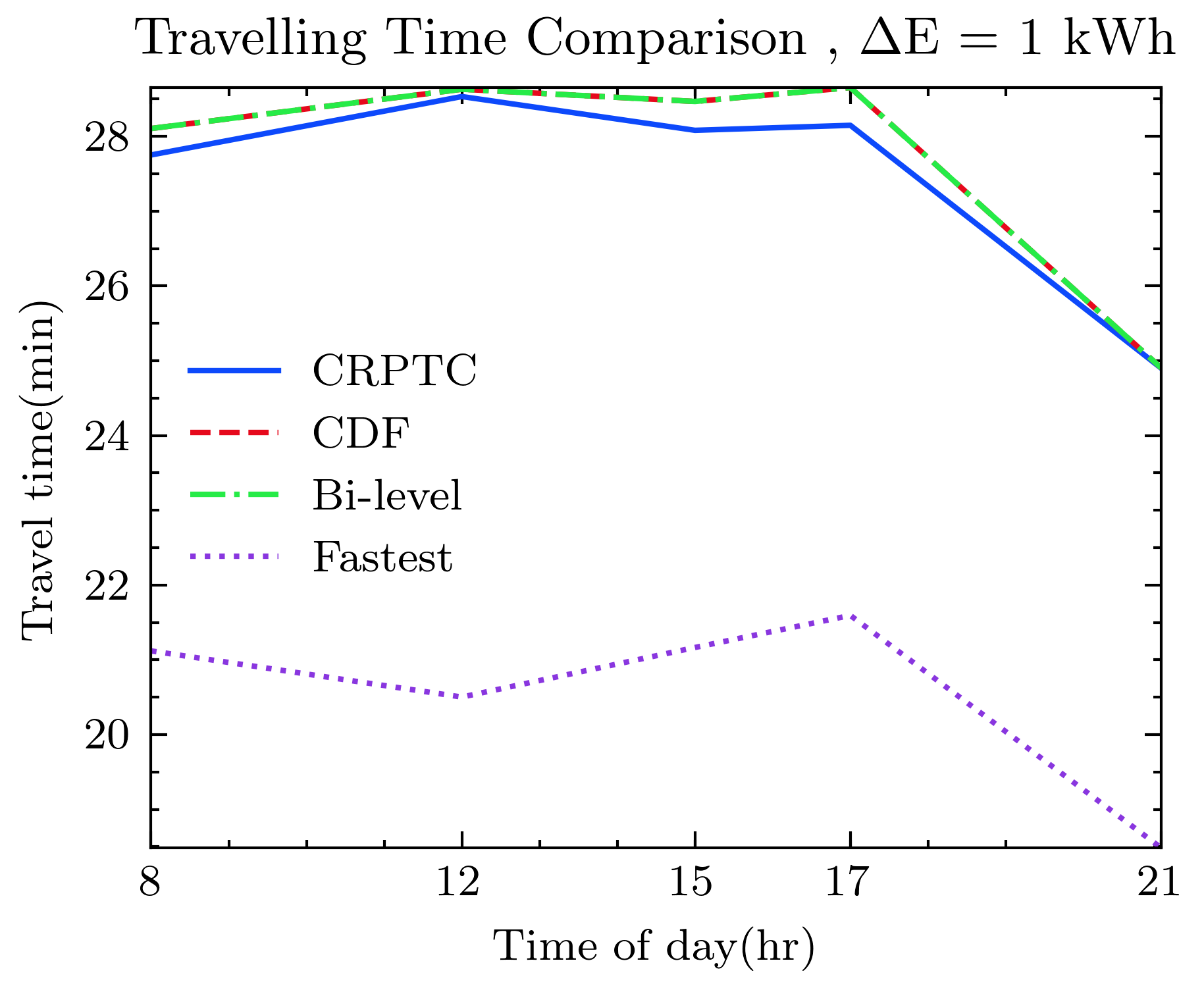}
\caption{}
\label{fig:time comp E1}
\end{subfigure}%
\begin{subfigure}{.24\textwidth}
\centering
\includegraphics[width=0.98\linewidth]{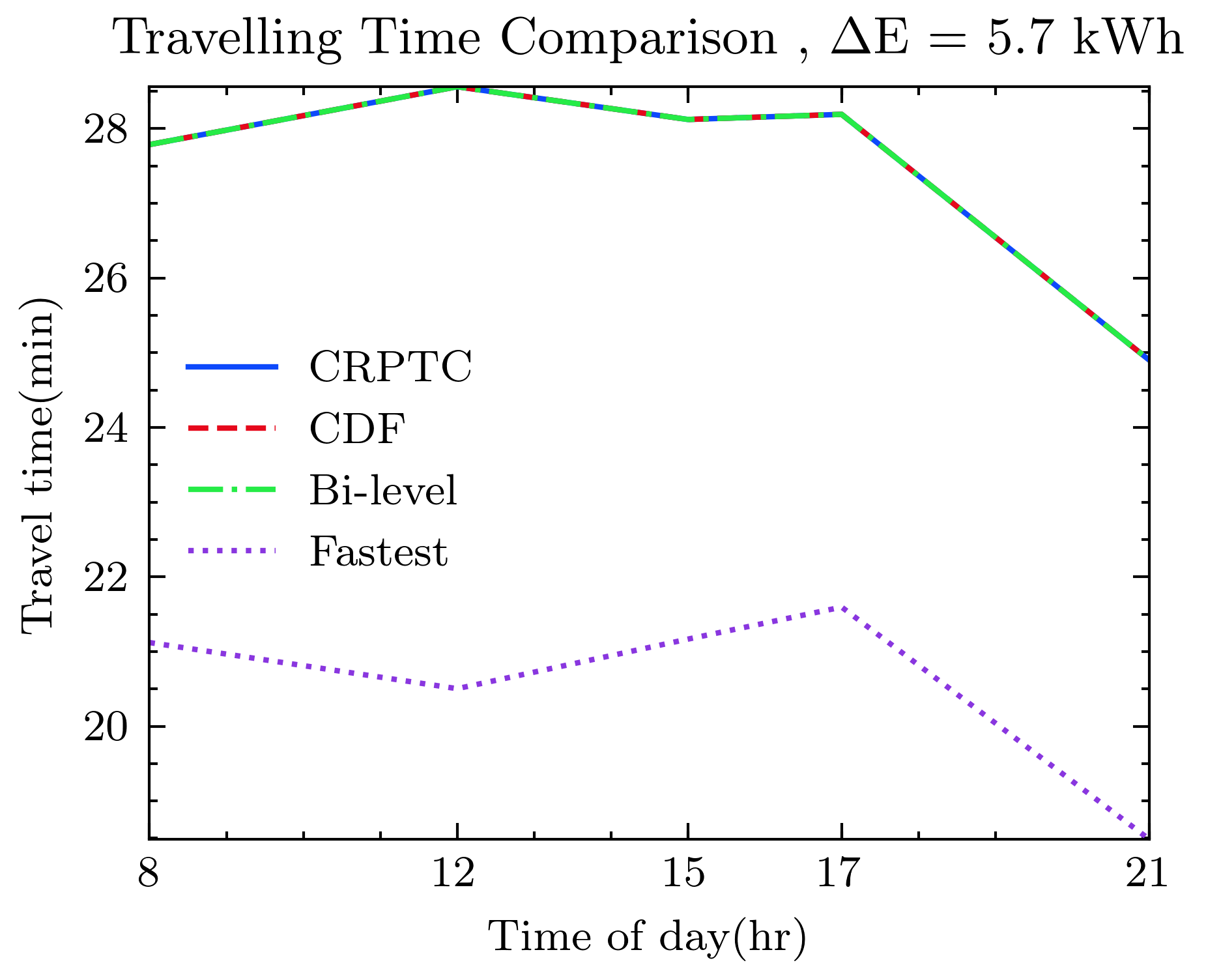}
\caption{}
\label{fig:time comp E5.7}
\end{subfigure}
\caption{ Average travel time comparison for the selected O-D pairs with
different allowable $\Delta E$ values}%
\label{fig:Travel times}%
\end{figure}
We have also compared the energy cost of different eco-routing algorithms
against the fastest route and compared them to each other showing their energy
saving distributions in Fig. \ref{fig:energy savings box} as box-plots where red line is the median and green triangle is the mean. Note that in these plots we are not showing the outlier data points.

\begin{figure}[h]
\centering
\begin{subfigure}{.24\textwidth}
\centering
\includegraphics[width=0.98\linewidth]{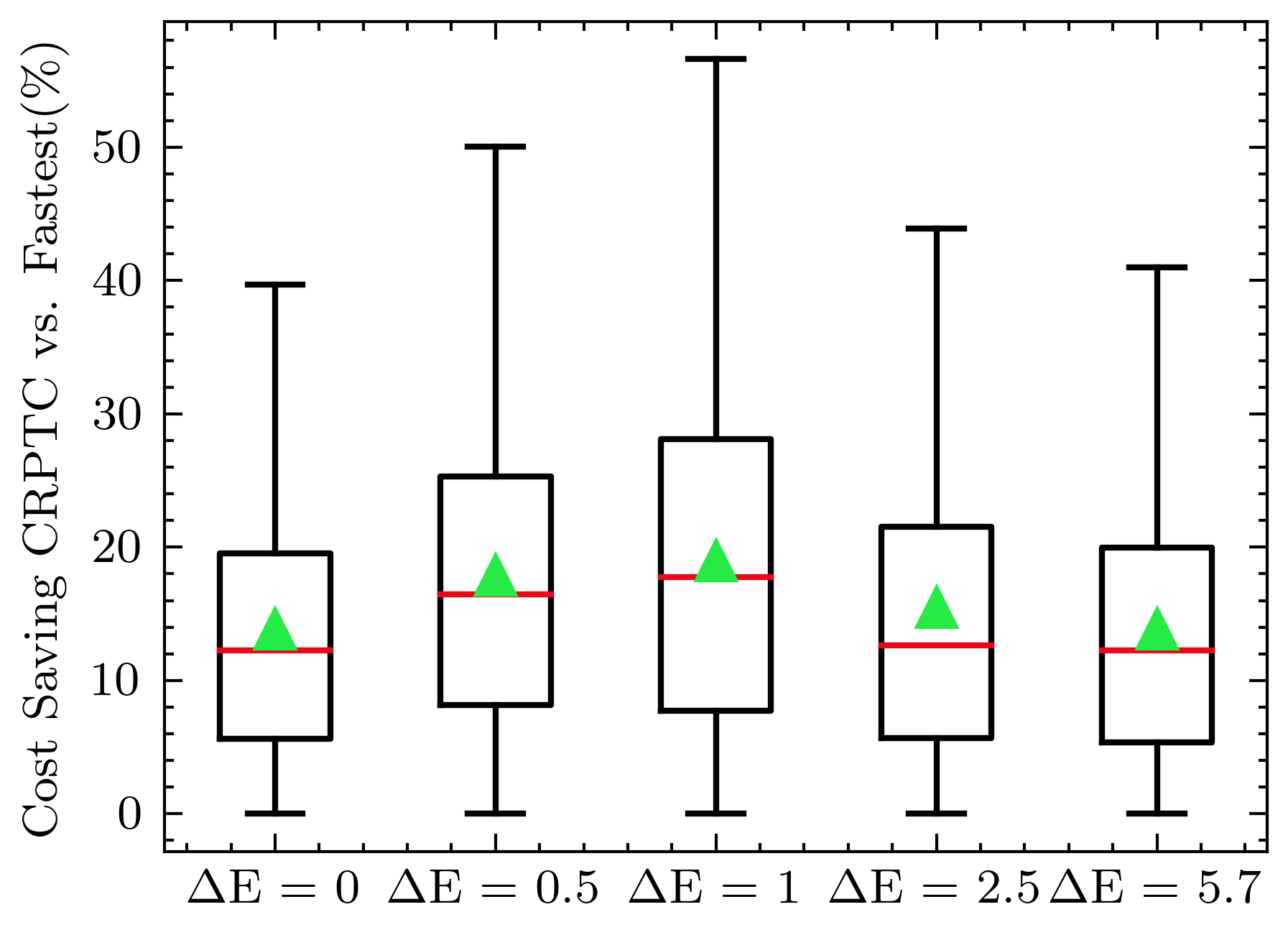}
\caption{}
\label{fig:cost saving box CRPTC vs fastest}
\end{subfigure}\begin{subfigure}{.24\textwidth}
\centering
\includegraphics[width=0.98\linewidth]{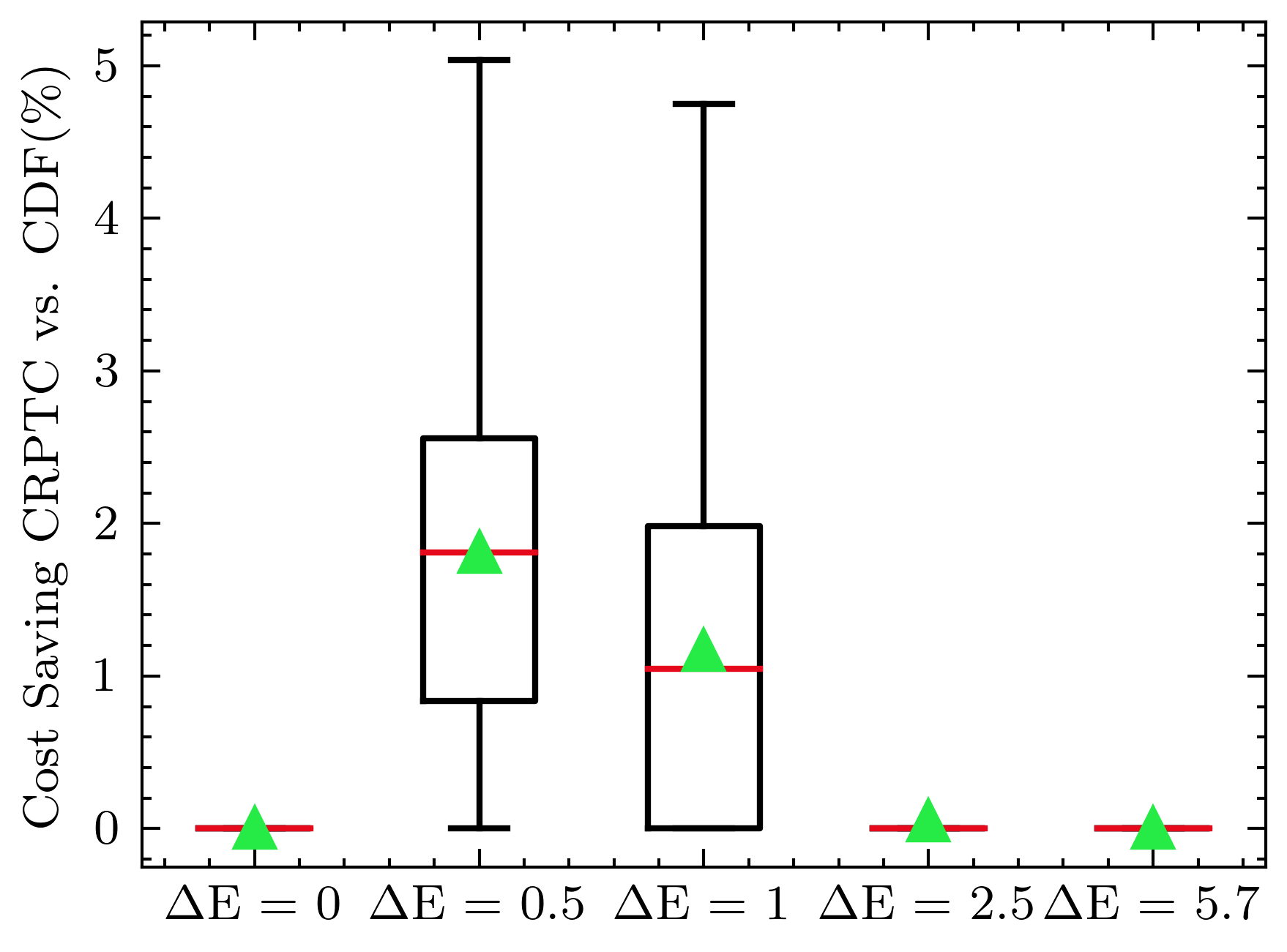}
\caption{}
\label{fig:cost saving box CRPTC vs CDF}
\end{subfigure}
\begin{subfigure}{.24\textwidth}
\centering
\includegraphics[width=0.98\linewidth]{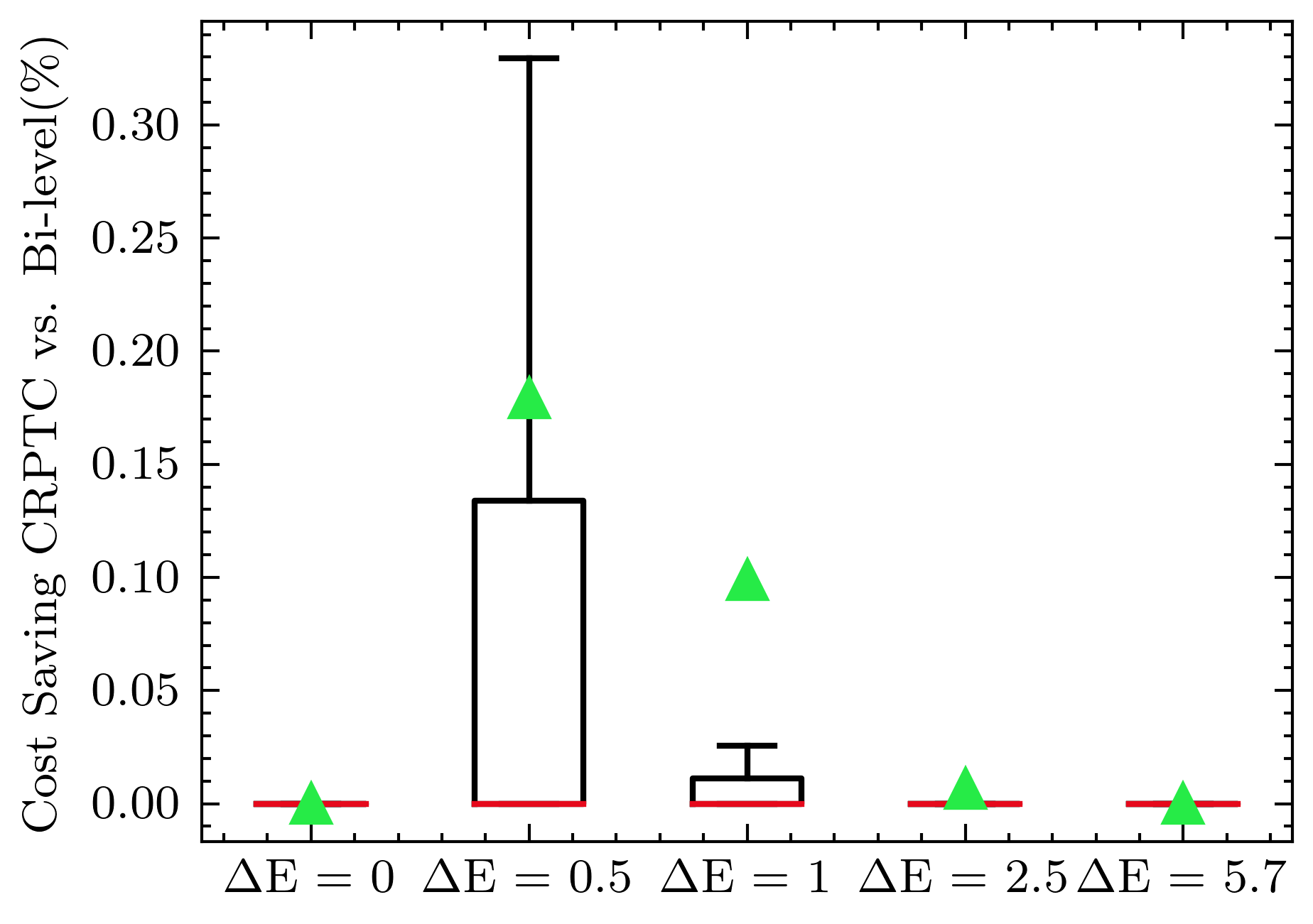}
\caption{}
\label{fig:cost saving box CRPTC vs bilevel}
\end{subfigure}\begin{subfigure}{.24\textwidth}
\centering
\includegraphics[width=0.98\linewidth]{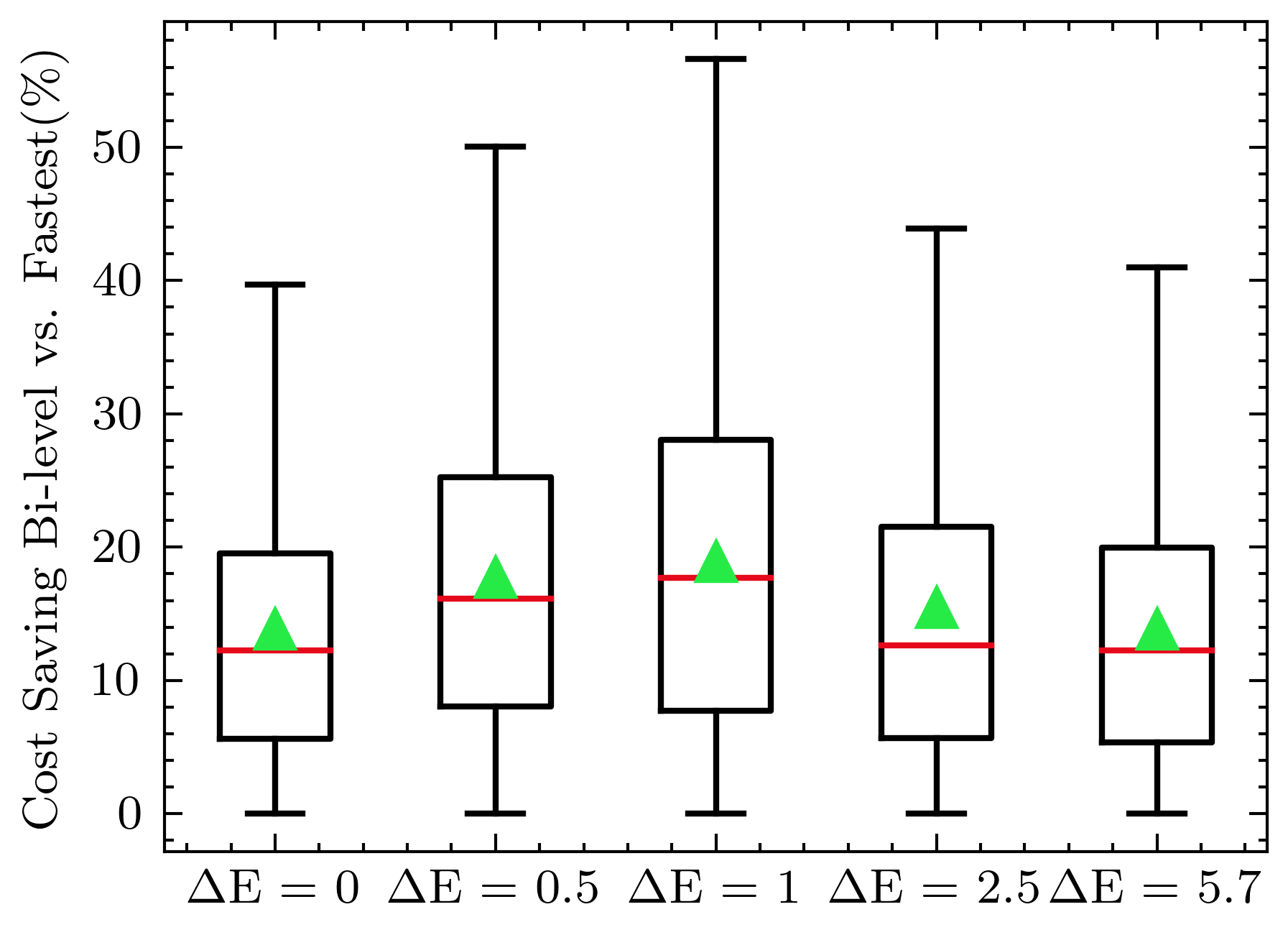}
\caption{}
\label{fig:cost saving box bilevel vs fastest}
\end{subfigure}
\begin{subfigure}
{.24\textwidth}
\centering
\includegraphics[width=0.98\linewidth]{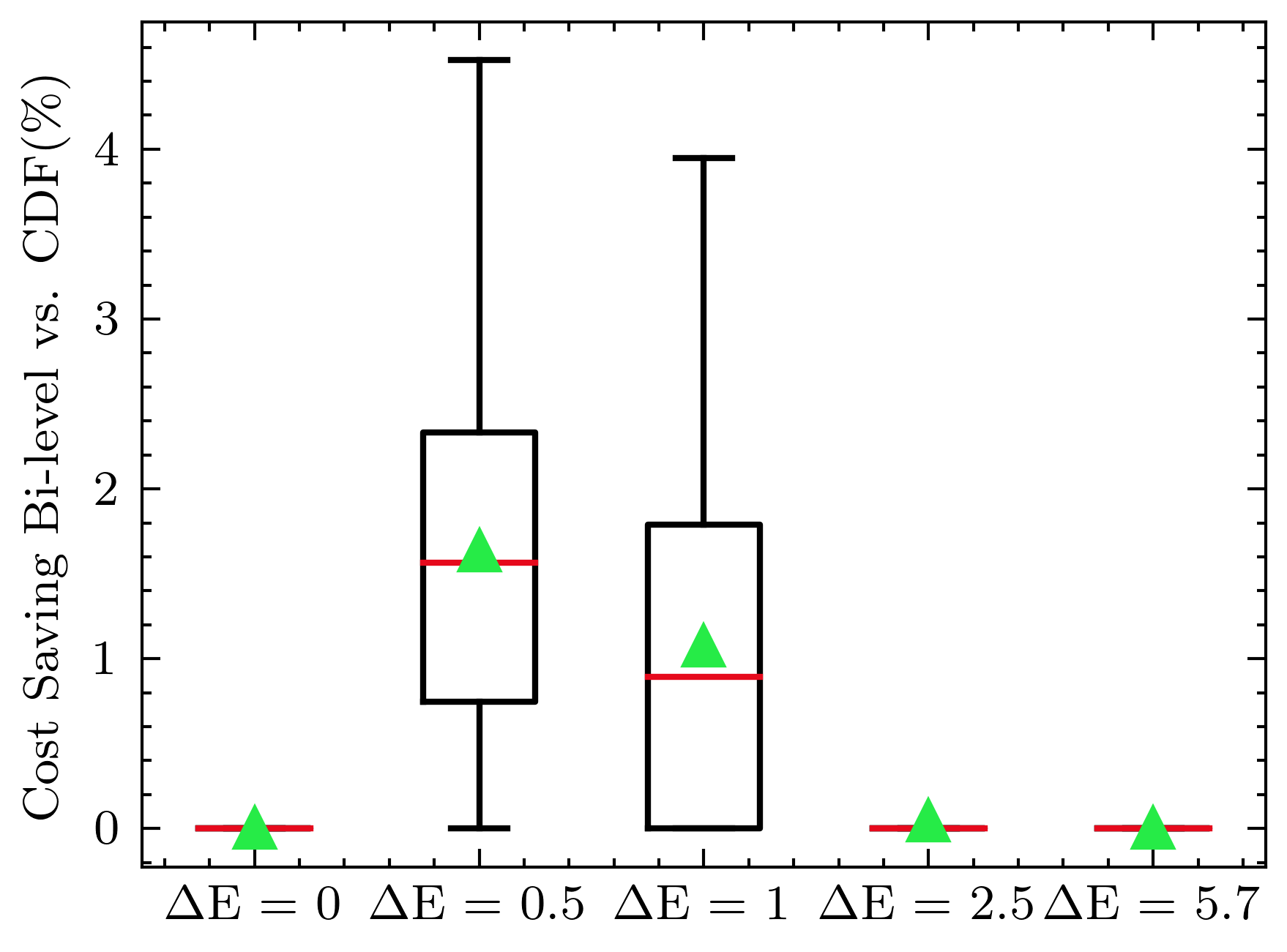}
\caption{}
\label{fig:cost saving box bilevel vs CDF}
\end{subfigure}\begin{subfigure}{.24\textwidth}
\centering
\includegraphics[width=0.98\linewidth]{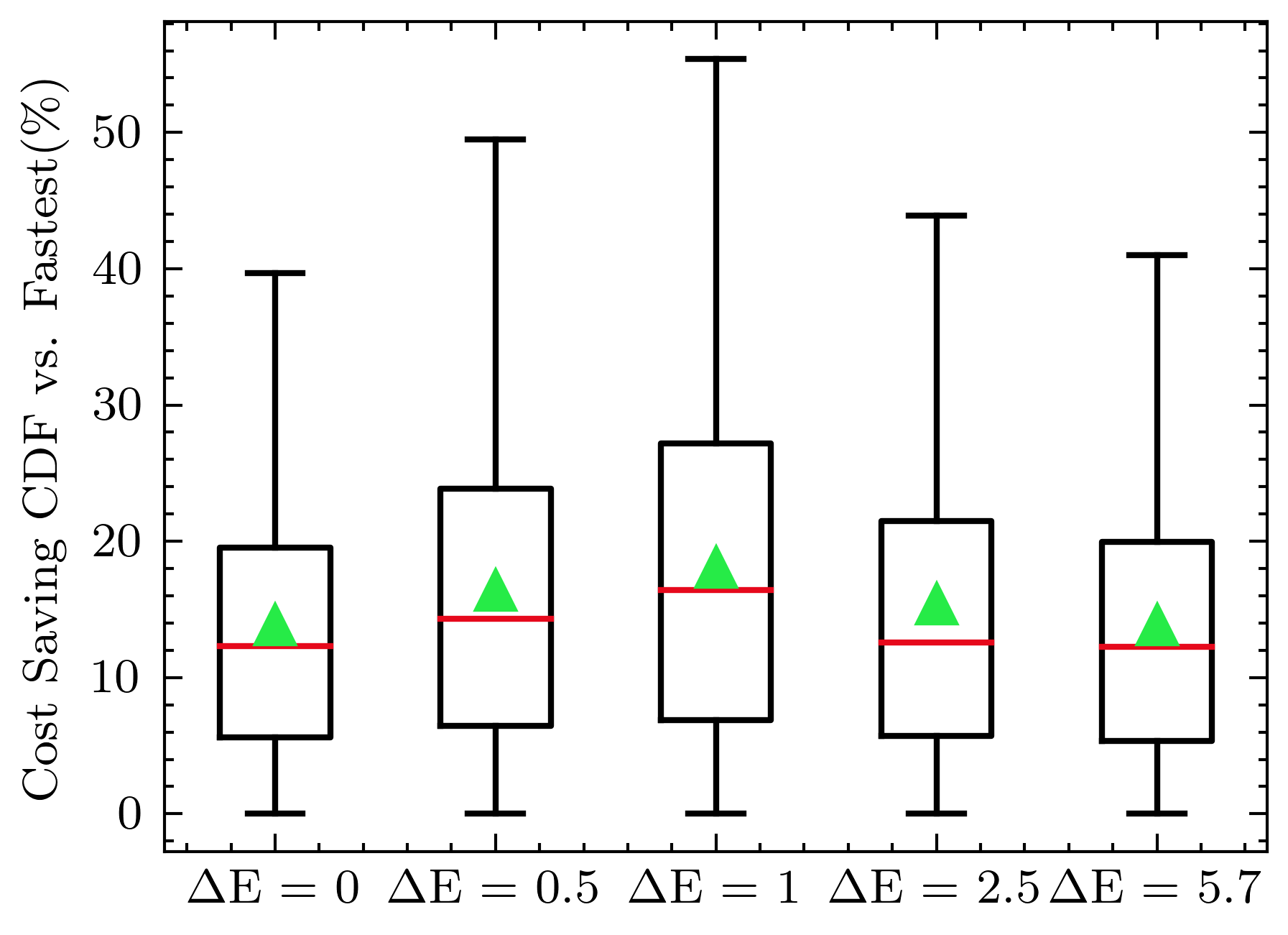}
\caption{}
\label{fig:cost saving CDF vs fastest}
\end{subfigure}
\caption{ Average energy cost saving distribution for the selected O-D pairs}%
\label{fig:energy savings box}%
\end{figure}
As we can see in Figs. \ref{fig:Energy costs} and \ref{fig:Travel times}, there is a trade-off between
energy and time while travelling through the eco-route and fastest route. To
better quantify this trade-off, we compare the travelling time of each
eco-routing algorithm (CRPTC, CDF, and bi-level) against that of the fastest
route and show their box-plot distribution in Fig. \ref{fig:time savings box} . Note
that the eco-route that is chosen by the CDF algorithm is the same as that of
the bi-level approach; as a result the traveling time results of these two
approaches are the same.

\begin{figure}[h]
\centering
\begin{subfigure}{.24\textwidth}
\centering
\includegraphics[width=0.98\linewidth]{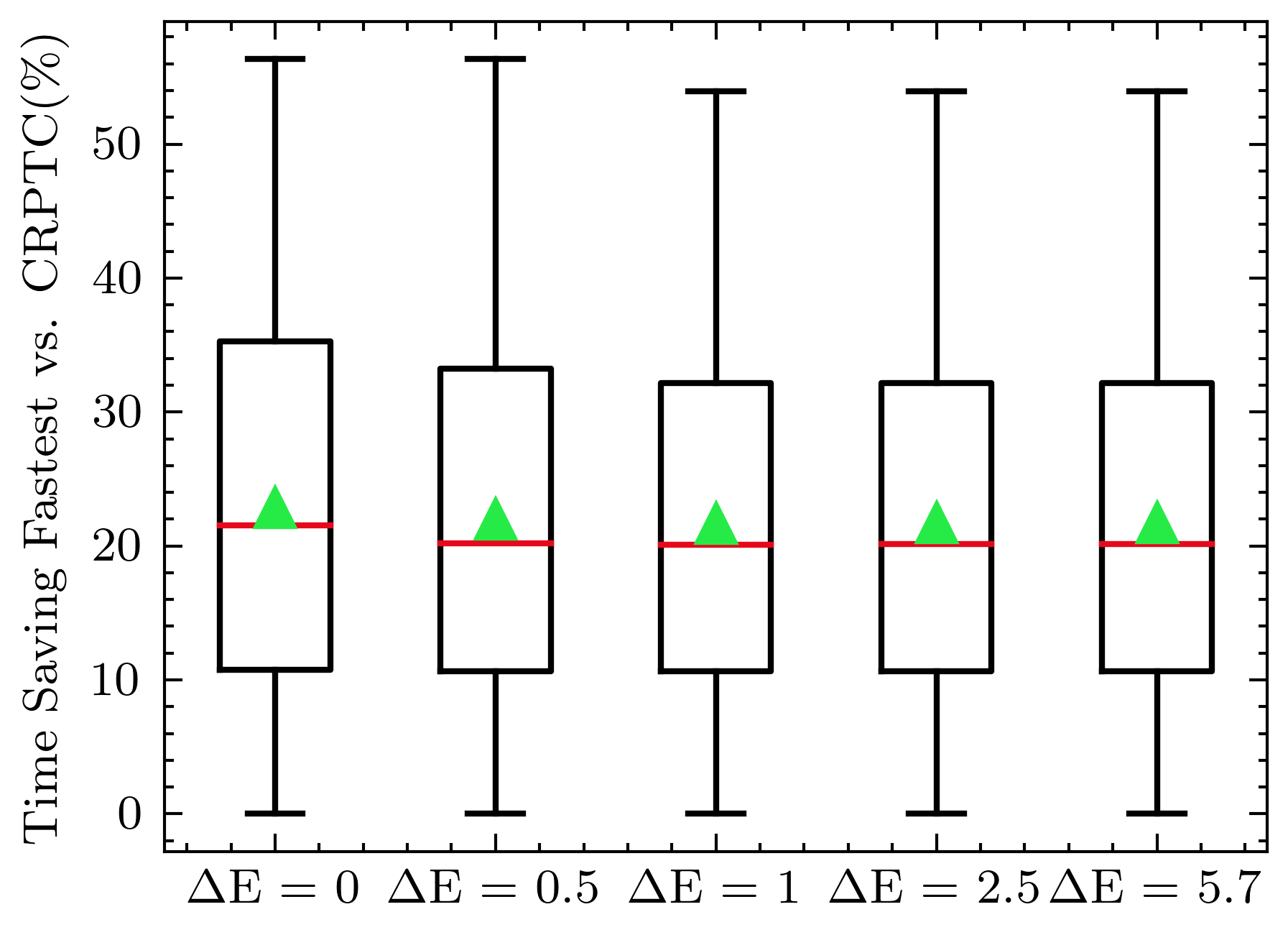}
\caption{}
\label{fig:Time saving box fastest vs CRPTC}
\end{subfigure}\begin{subfigure}{.24\textwidth}
\centering
\includegraphics[width=0.98\linewidth]{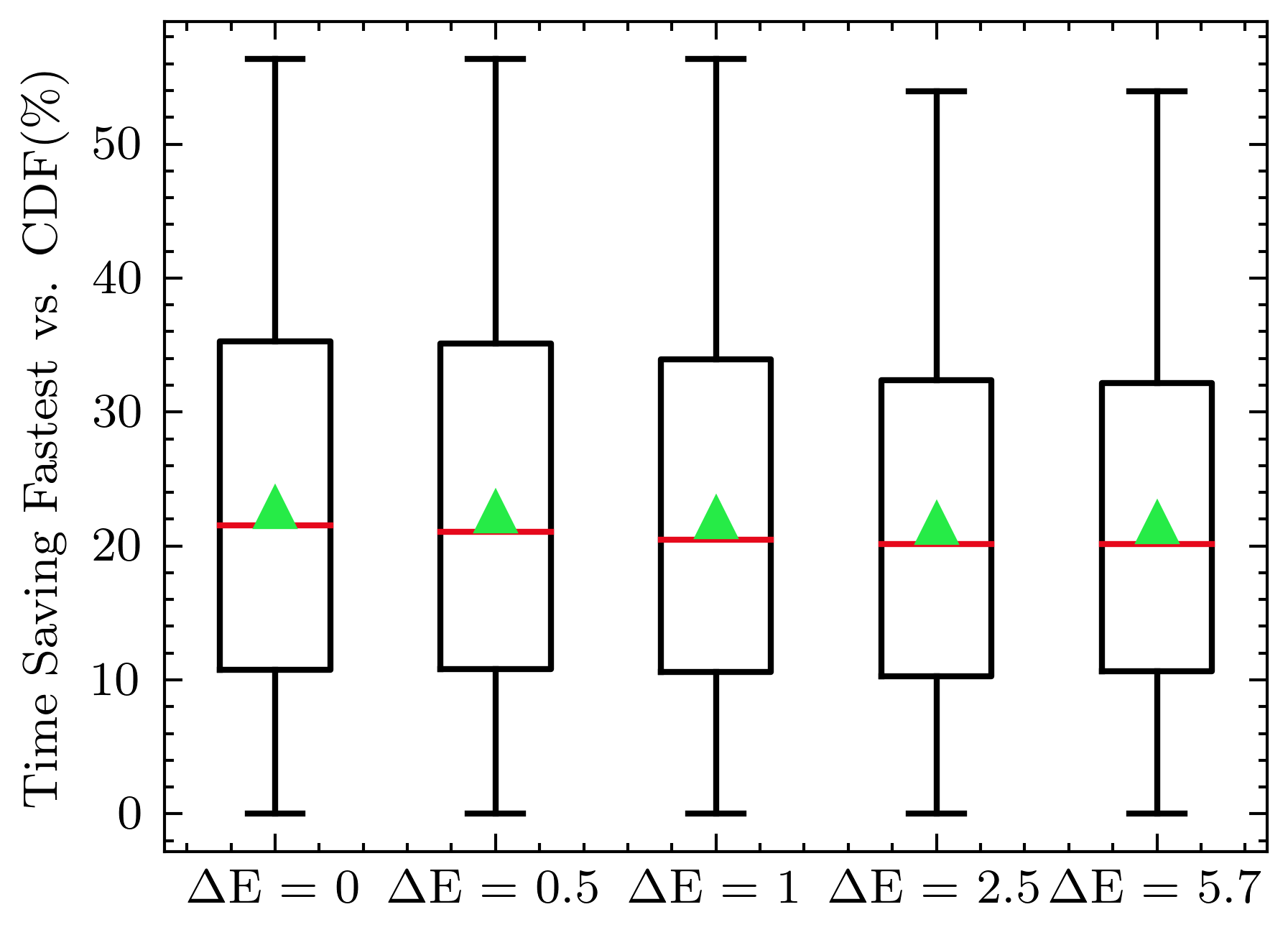}
\caption{}
\label{fig:Time saving box Fastest vs CDF}
\end{subfigure}
\begin{subfigure}{.24\textwidth}
\centering
\includegraphics[width=0.98\linewidth]{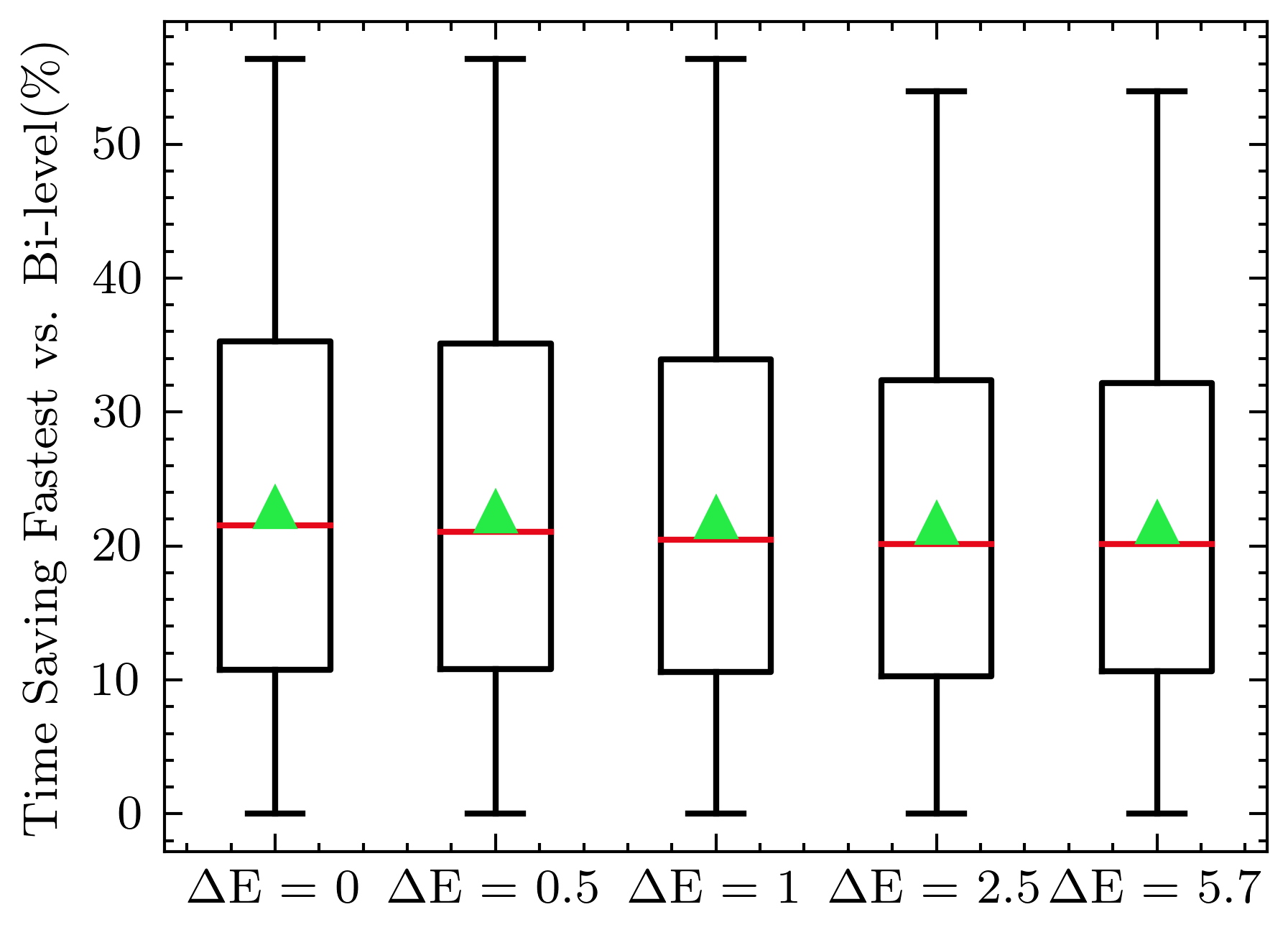}
\caption{}
\label{fig:Time saving box Fastest vs Bi-level}
\end{subfigure}\begin{subfigure}{.24\textwidth}
\centering
\includegraphics[width=0.98\linewidth]{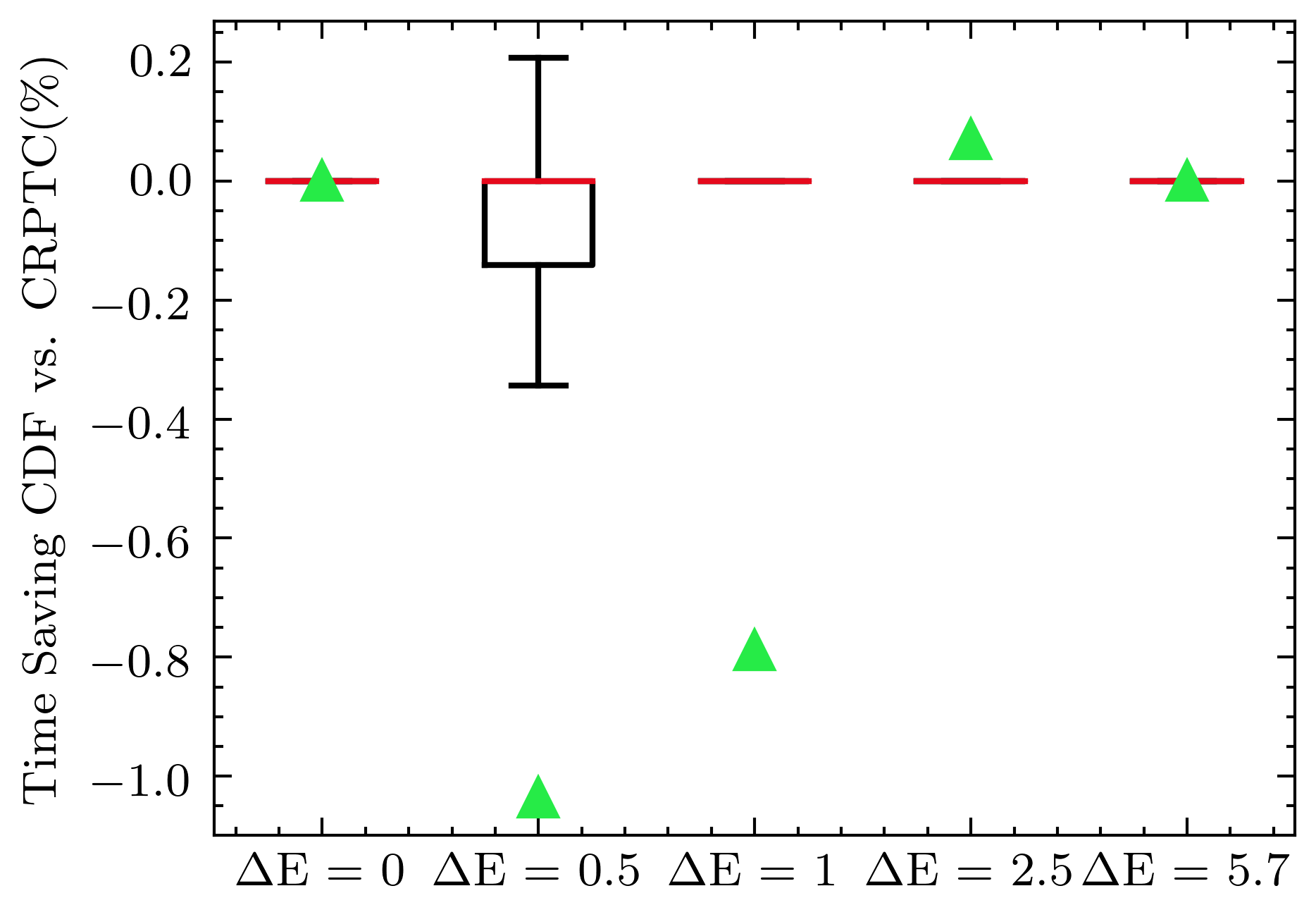}
\caption{}
\label{fig:Time saving box CDF vs CRPTC}
\end{subfigure}
\caption{ Average travel time saving distribution for the selected O-D pairs}%
\label{fig:time savings box}%
\end{figure}

To show the interdependence of eco-routing performance and the distance
between O-D pairs, we organized O-D pairs based on the shortest distance
between them and reported the average energy saving and travel time saving
values based on the distance between origin and destination in Figs.
\ref{fig:Energy savings based on distance} and
\ref{fig:Time savings based on distance}.
\begin{figure*}[ptb]
\centering
\begin{subfigure}{.33\textwidth}
\centering
\includegraphics[width=0.98\linewidth]{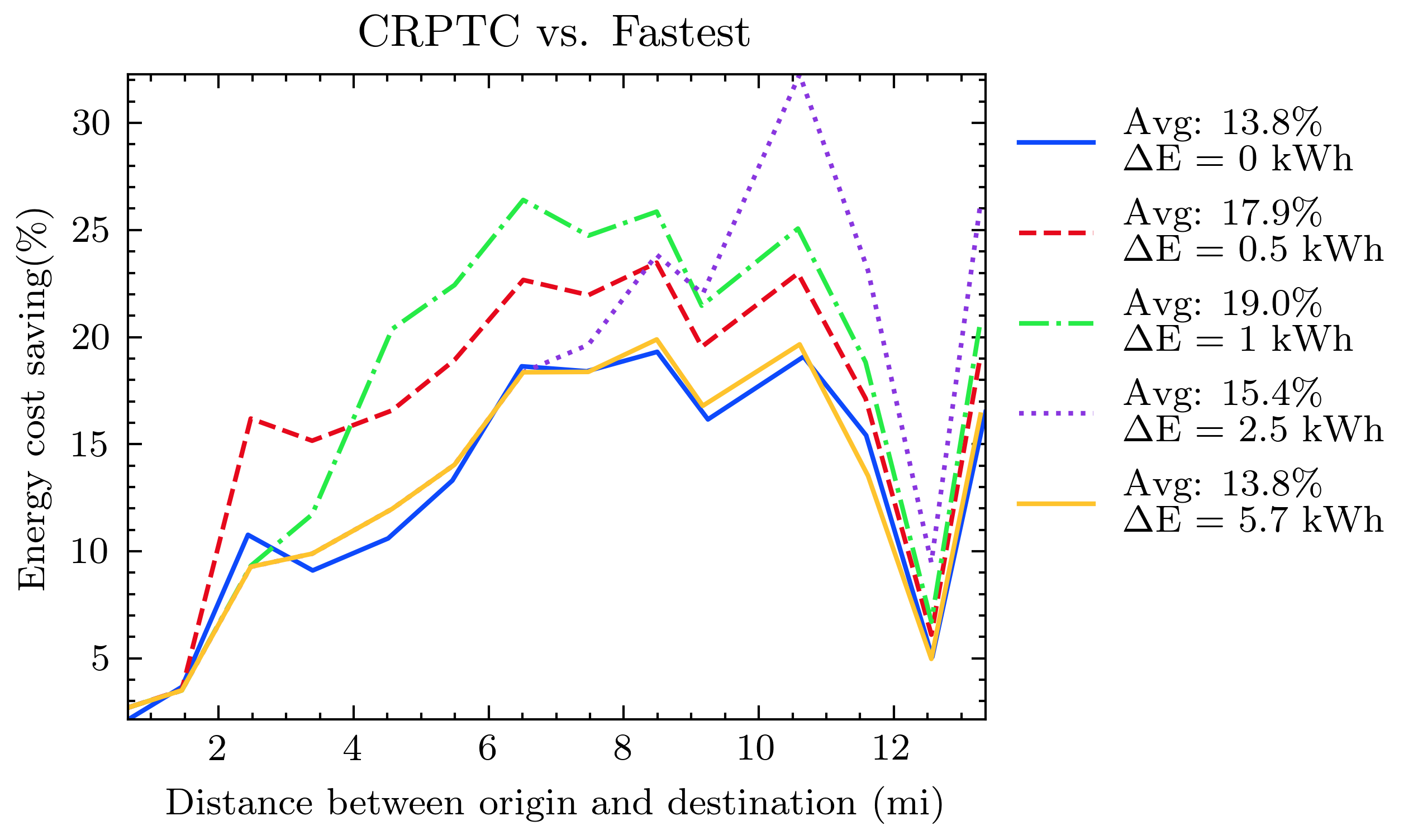}
\caption{}
\label{fig:cost saving CRPTC vs fastest}
\end{subfigure}\begin{subfigure}{.33\textwidth}
\centering
\includegraphics[width=0.98\linewidth]{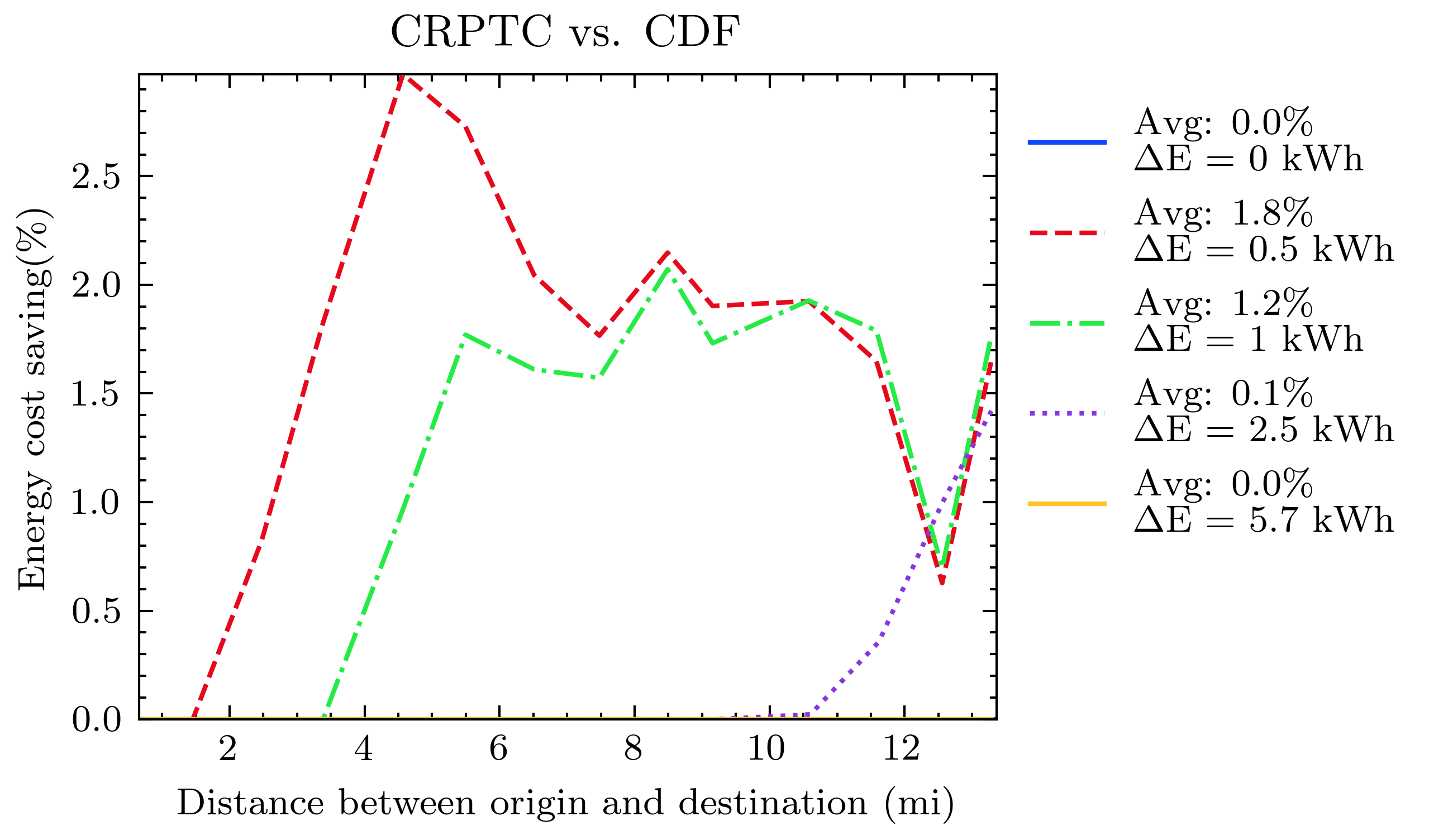}
\caption{}
\label{fig:cost saving CRPTC vs CDF}
\end{subfigure}\begin{subfigure}{.33\textwidth}
\centering
\includegraphics[width=0.98\linewidth]{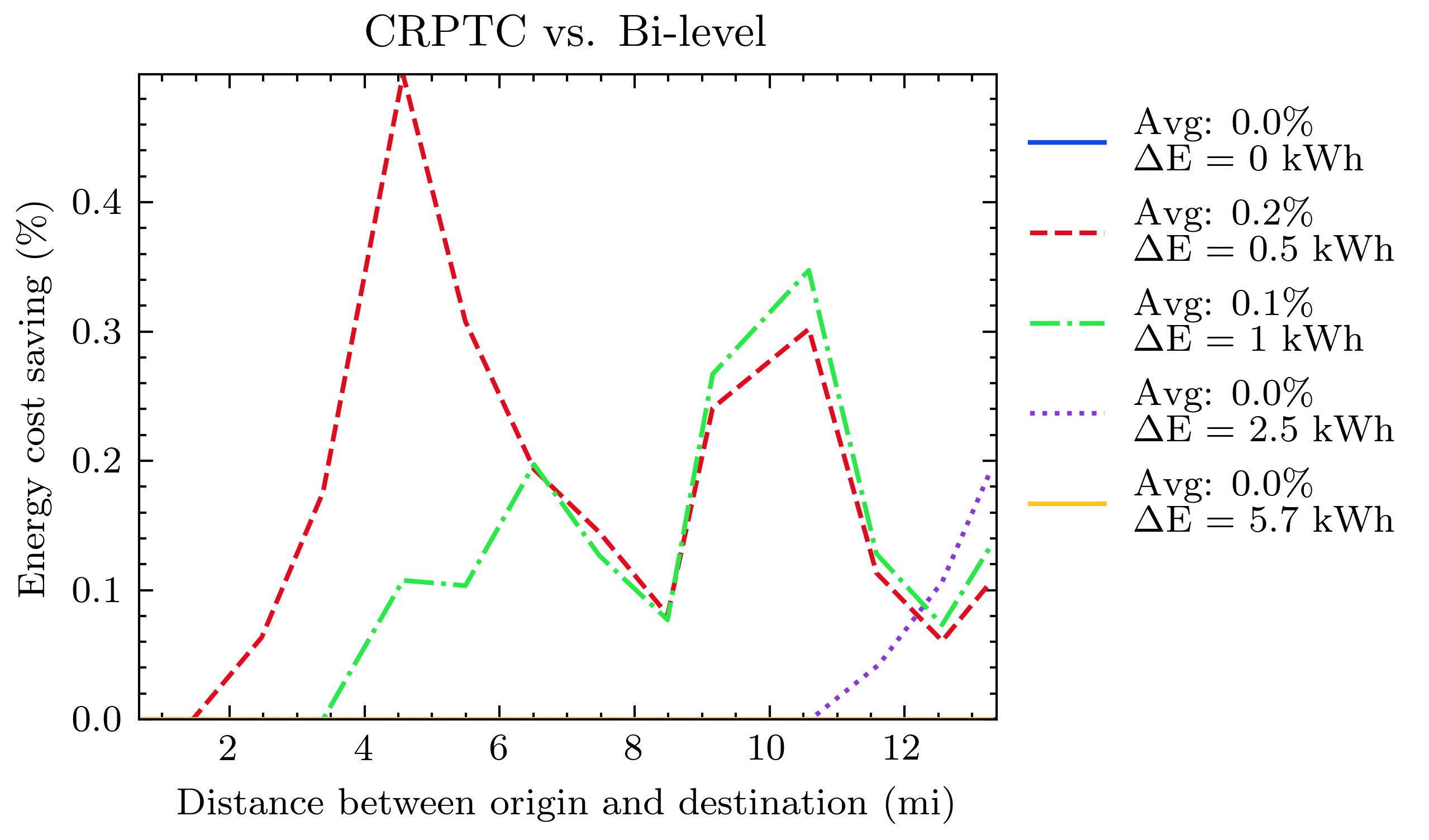}
\caption{}
\label{fig:cost saving CRPTC vs bilevel}
\end{subfigure}
\begin{subfigure}{.33\textwidth}
\centering
\includegraphics[width=0.98\linewidth]{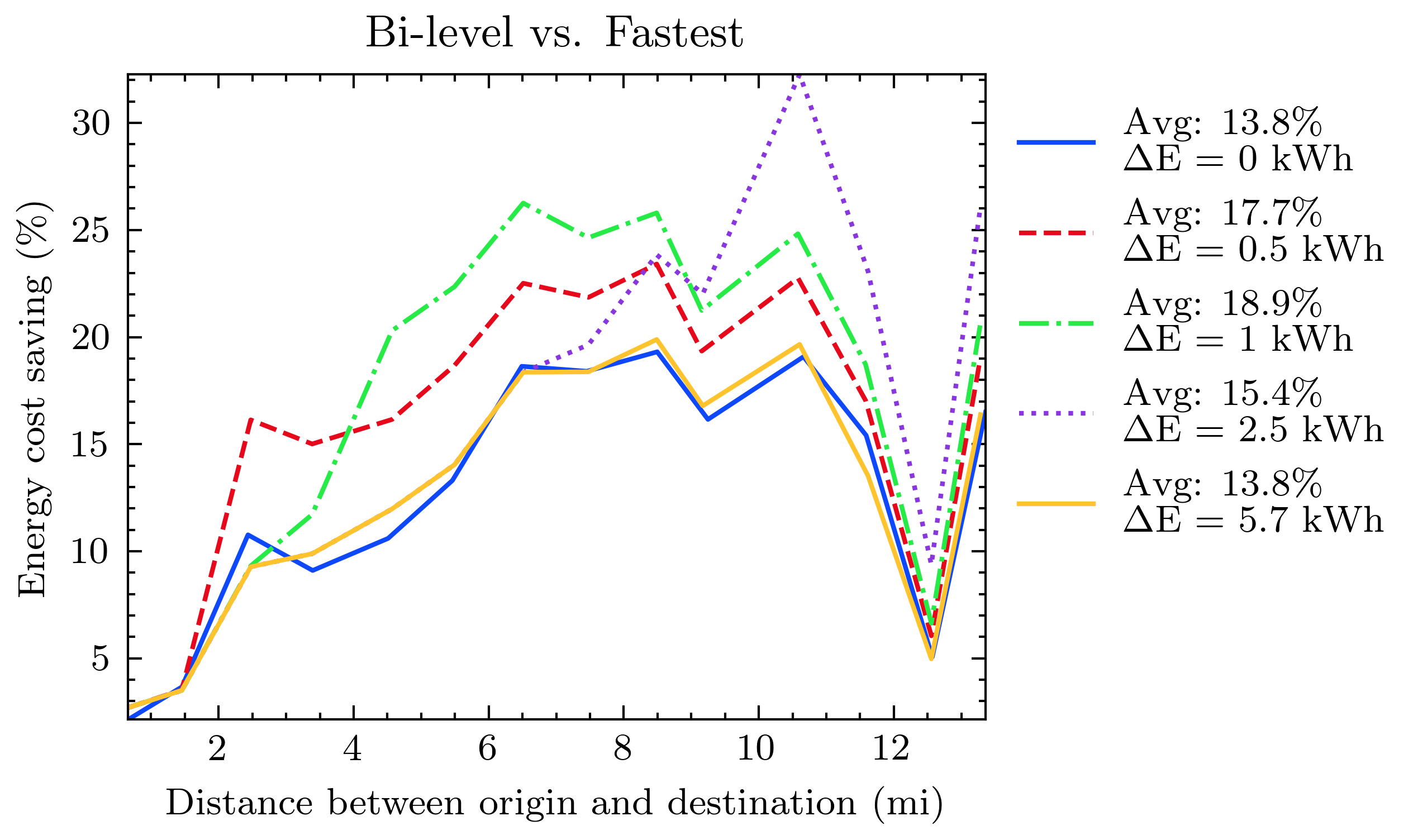}
\caption{}
\label{fig:cost saving bilevel vs fastest}
\end{subfigure}\begin{subfigure}{.33\textwidth}
\centering
\includegraphics[width=0.98\linewidth]{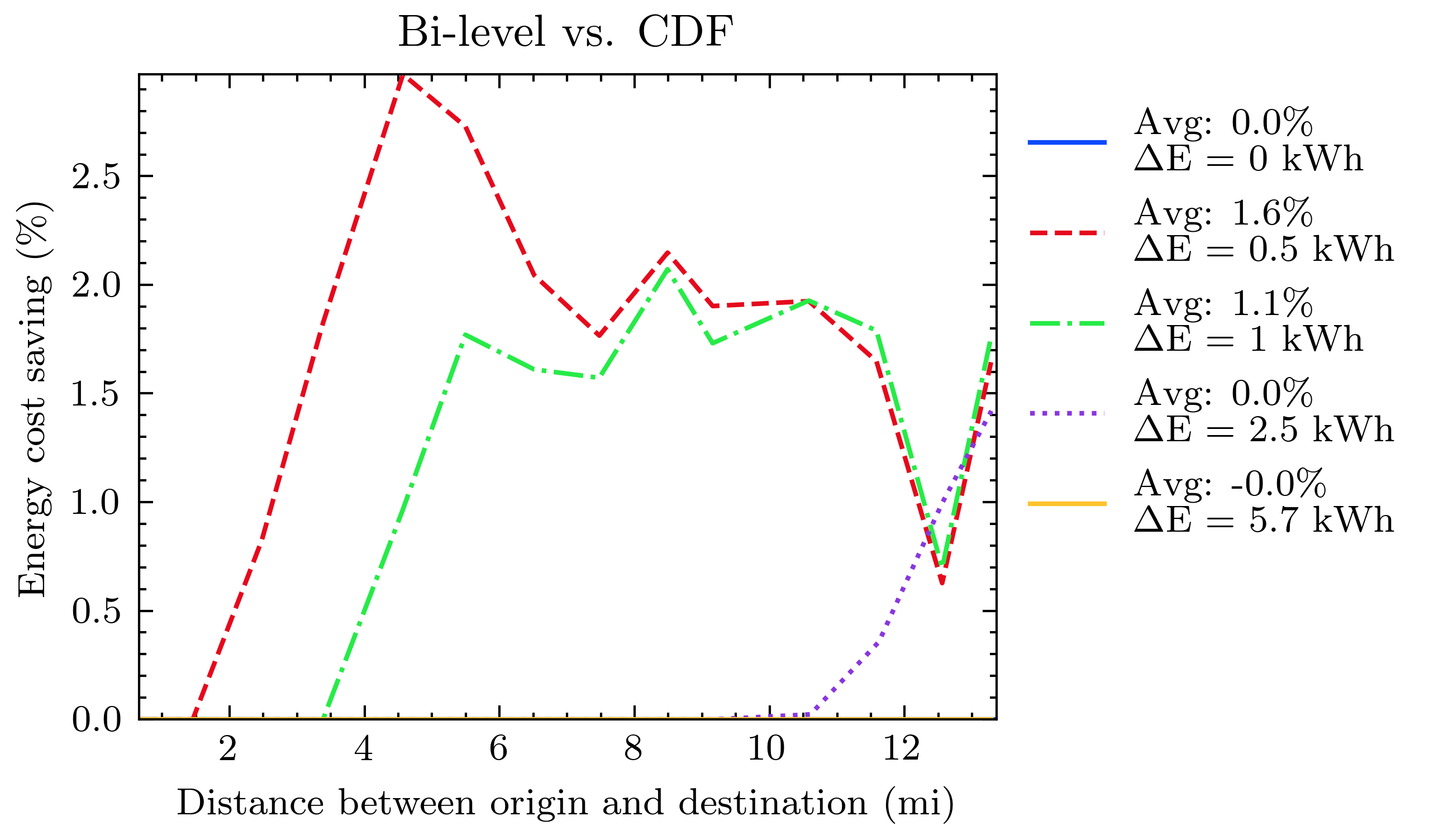}
\caption{}
\label{fig:cost saving bilevel vs CDF}
\end{subfigure}\begin{subfigure}{.33\textwidth}
\centering
\includegraphics[width=0.98\linewidth]{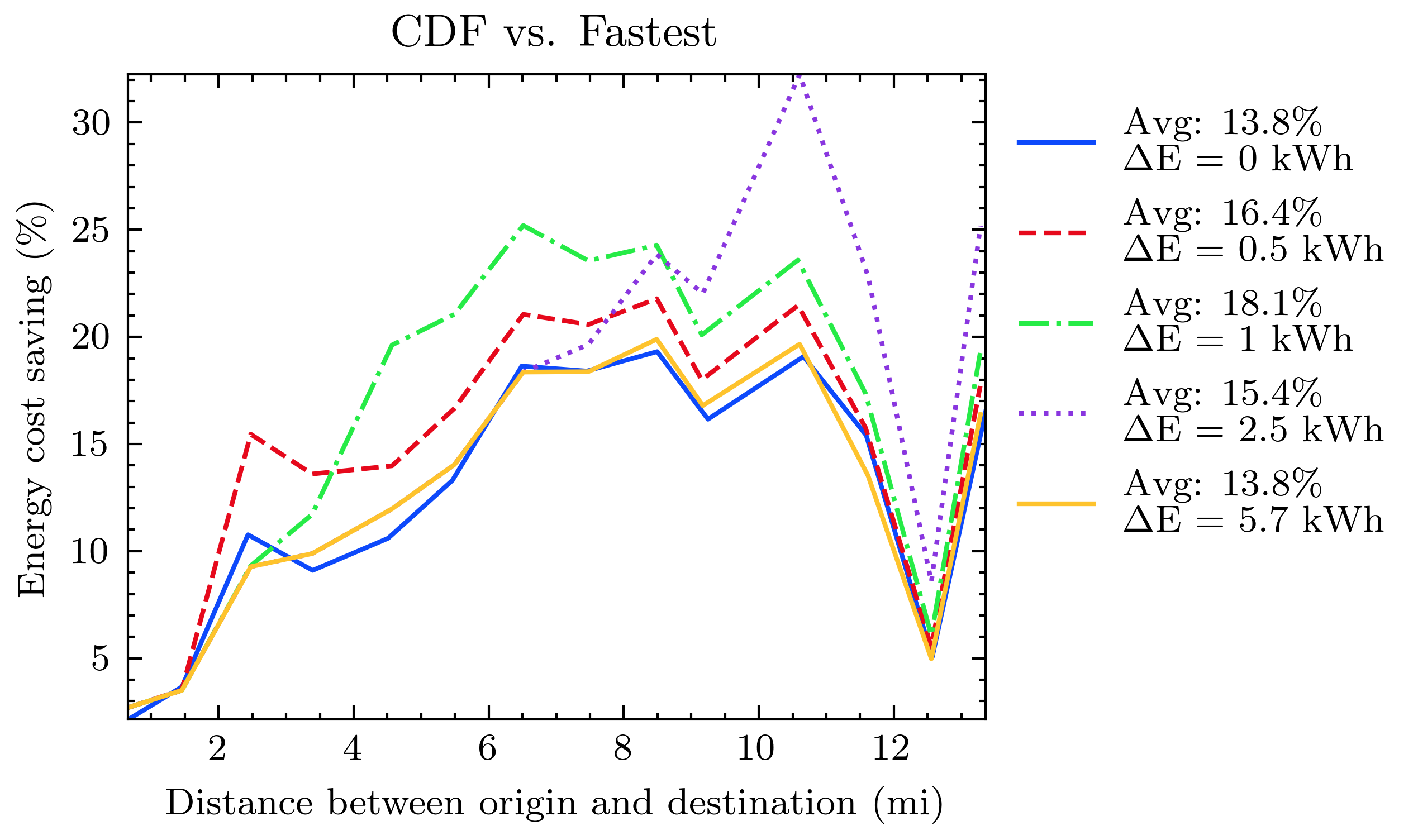}
\caption{}
\label{fig:cost saving CDF vs fastest}
\end{subfigure}
\caption{ Average energy cost savings for the selected O-D pairs based on the shortest distance between them}
\label{fig:Energy savings based on distance}%
\end{figure*}
\begin{figure*}[ptb]
\centering
\begin{subfigure}{.33\textwidth}
\centering
\includegraphics[width=0.98\linewidth]{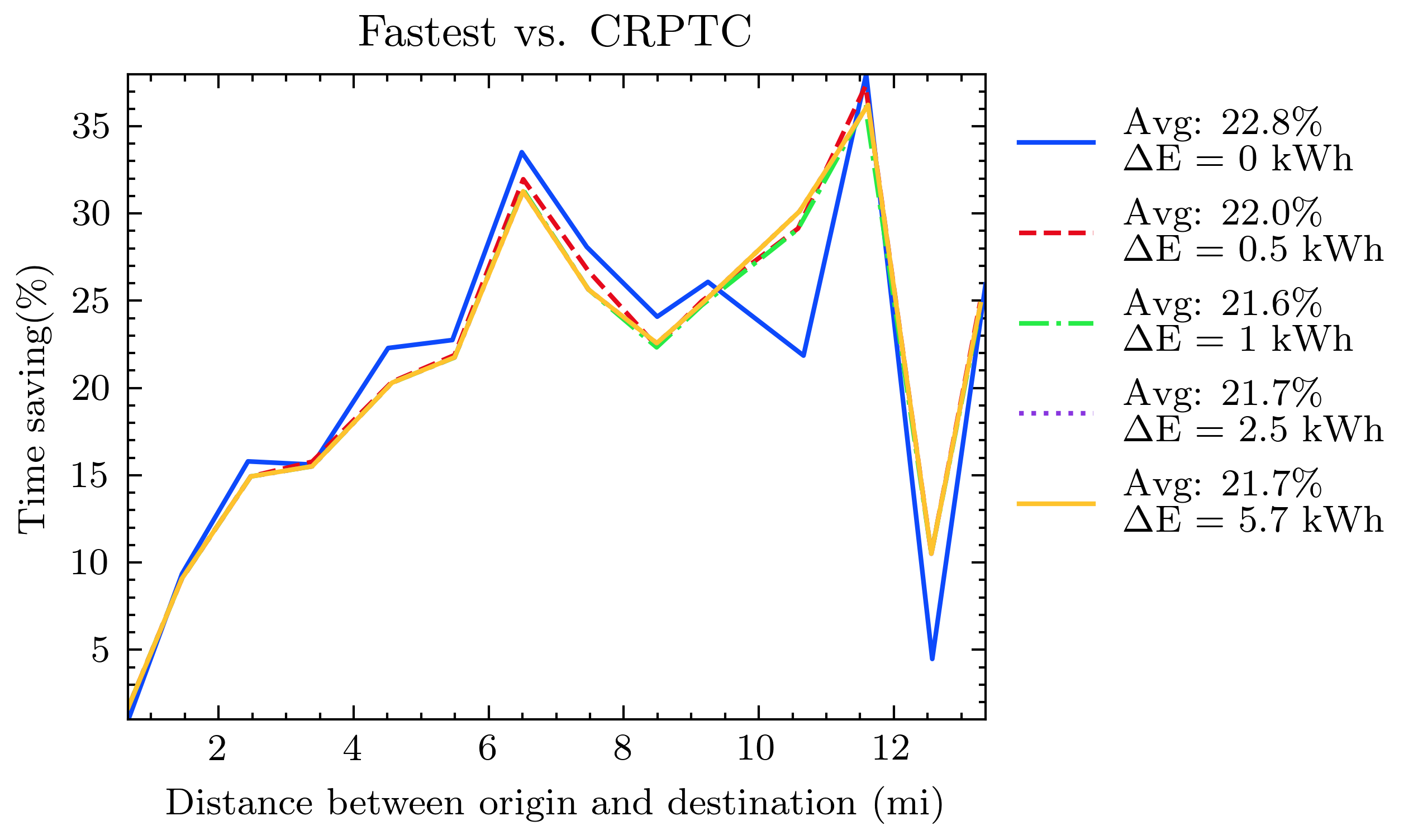}
\caption{}
\label{fig:time saving fastest vs CRPTC}
\end{subfigure}\begin{subfigure}{.33\textwidth}
\centering
\includegraphics[width=0.98\linewidth]{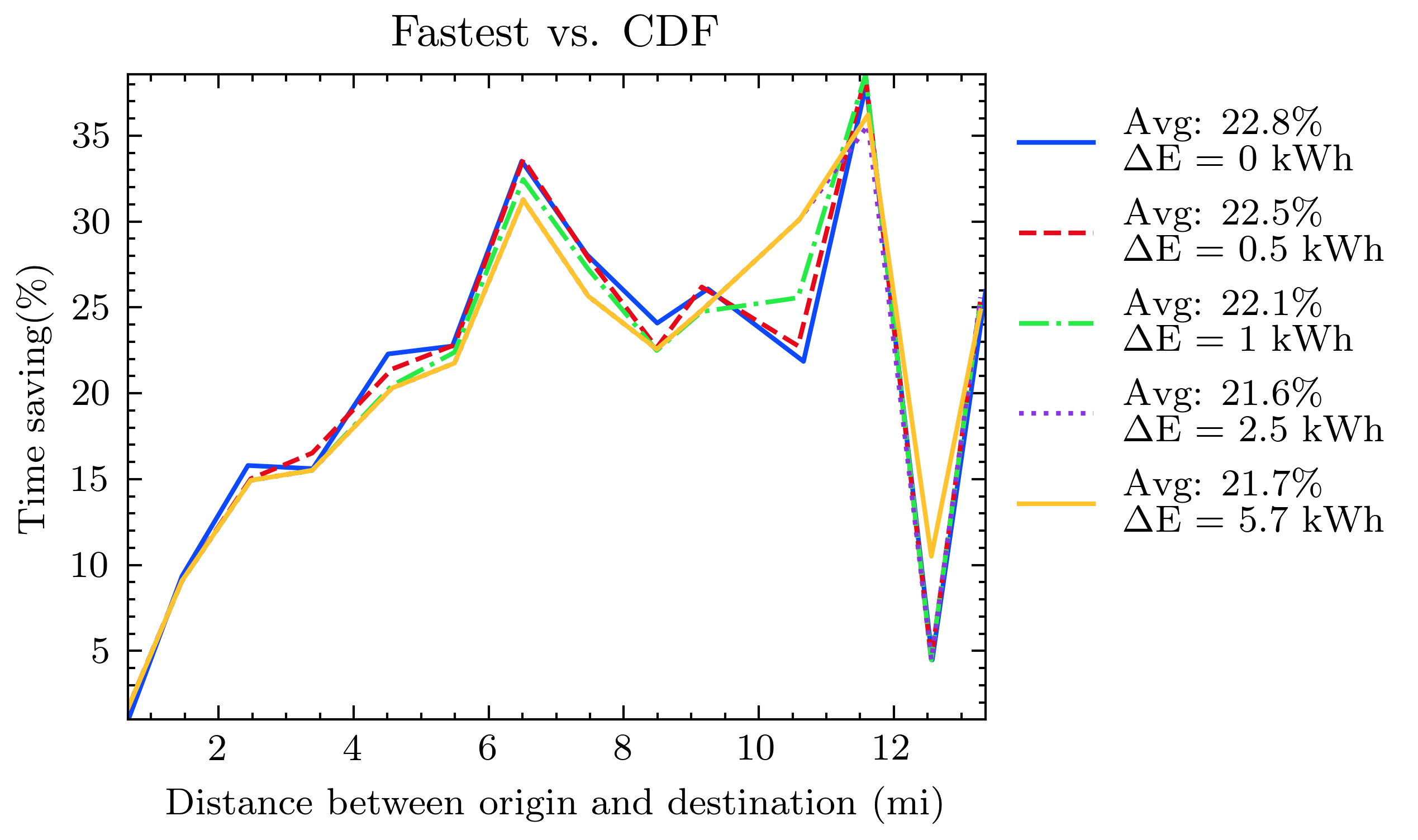}
\caption{}
\label{fig:time saving fastest vs CDF}
\end{subfigure}\begin{subfigure}{.33\textwidth}
\centering
\includegraphics[width=0.98\linewidth]{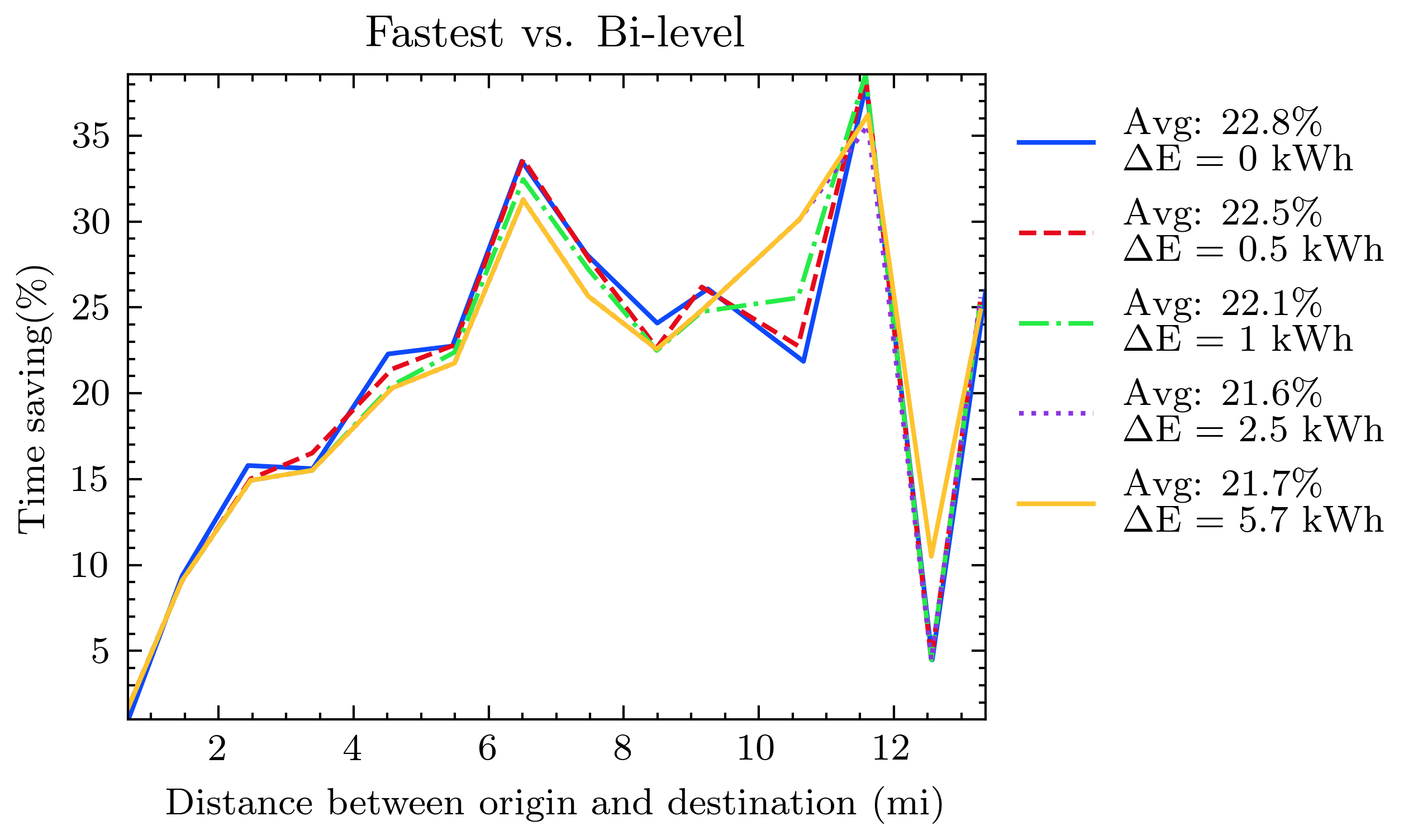}
\caption{}
\label{fig:time saving CDF vs CRPTC}
\end{subfigure}
\caption{ Average travel time savings for the selected O-D pairs based on the
shortest distance between them}%
\label{fig:Time savings based on distance}
\end{figure*}
\begin{figure*}[ptb]
\centering
\begin{subfigure}{.33\textwidth}
\centering
\includegraphics[width=0.98\linewidth]{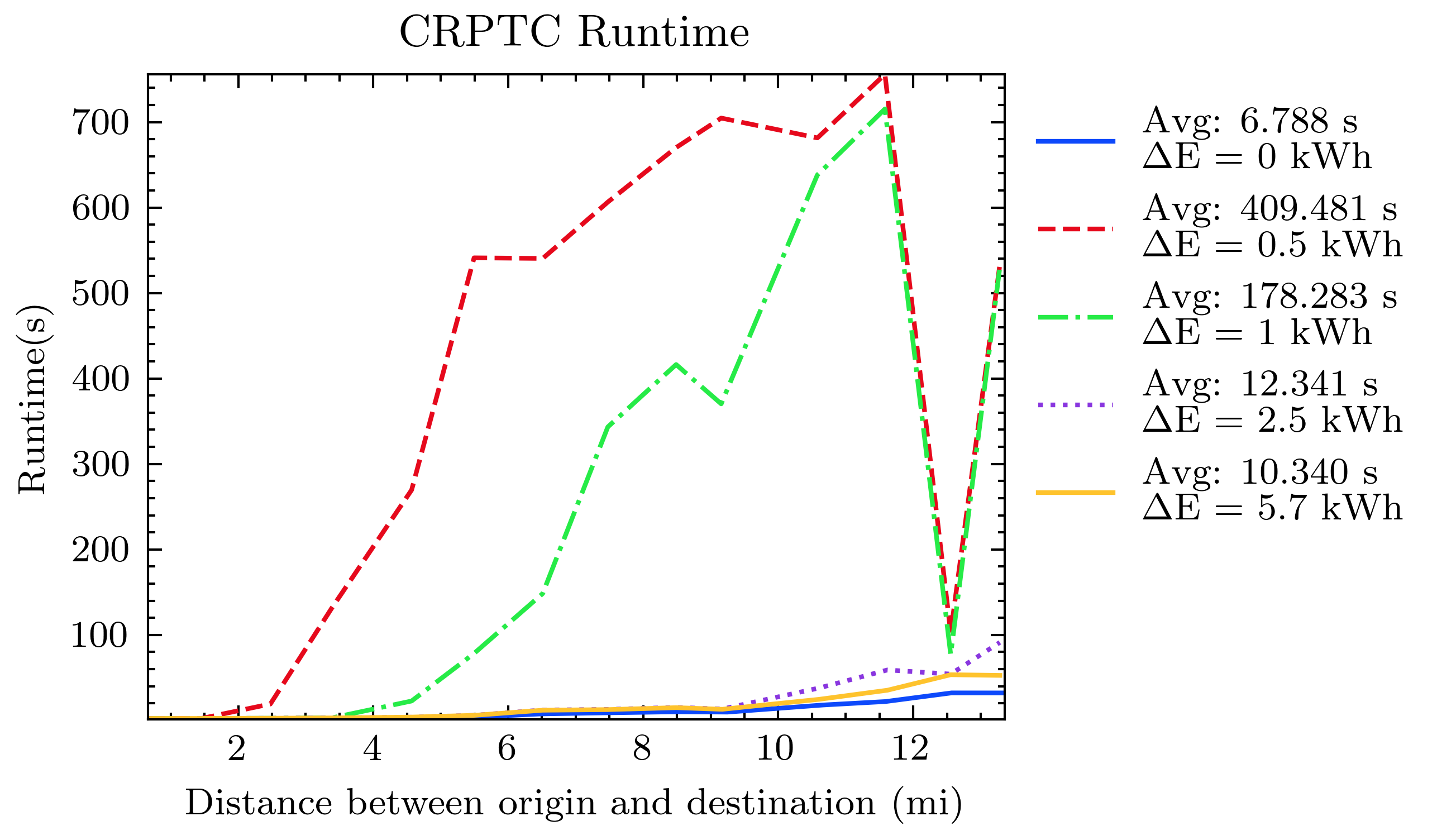}
\caption{}
\label{fig:runtime CRPTC - dist}
\end{subfigure}\begin{subfigure}{.33\textwidth}
\centering
\includegraphics[width=0.98\linewidth]{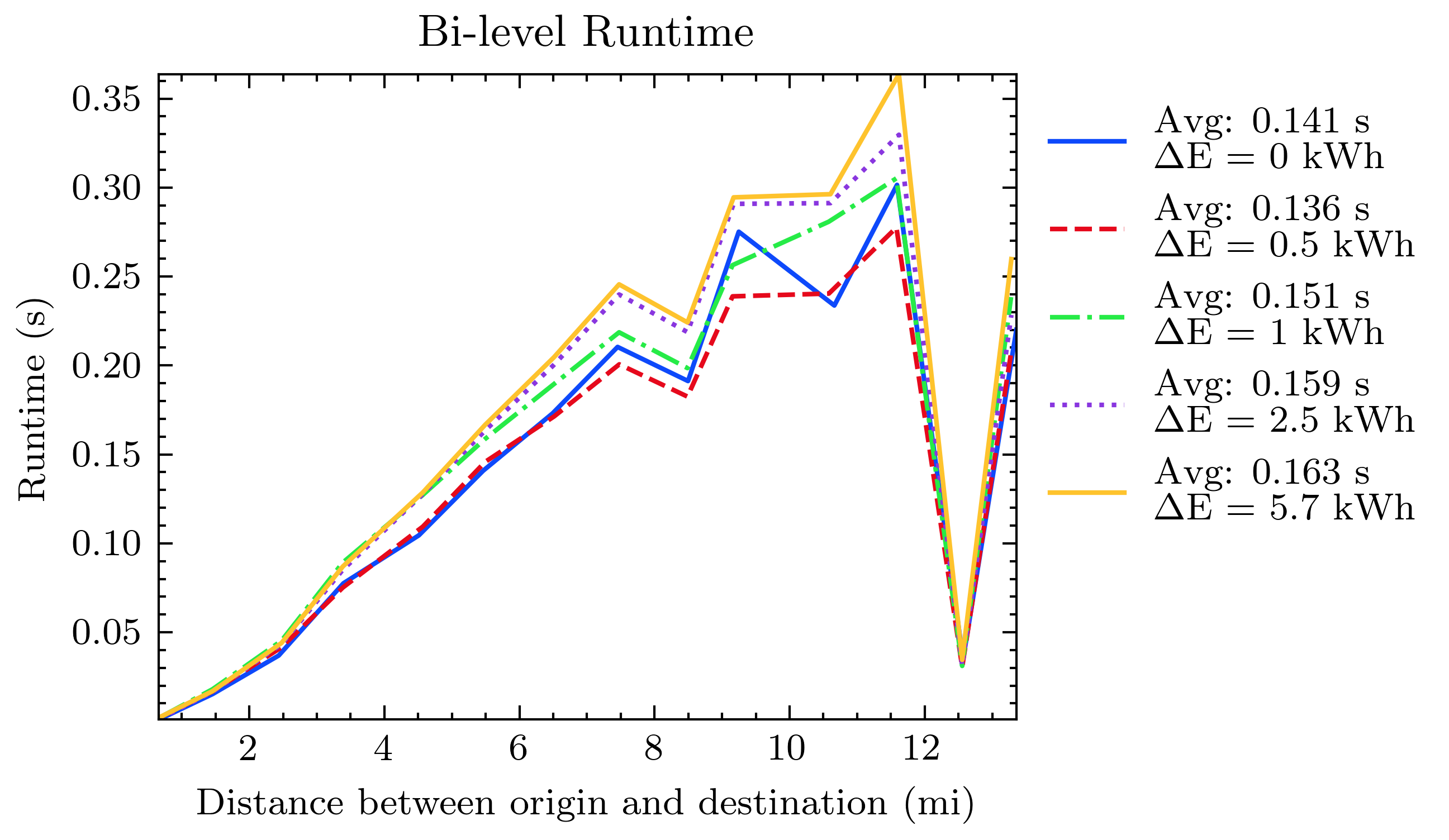}
\caption{}
\label{fig:runtime Bilevel - dist}
\end{subfigure}\begin{subfigure}{.33\textwidth}
\centering
\includegraphics[width=0.98\linewidth]{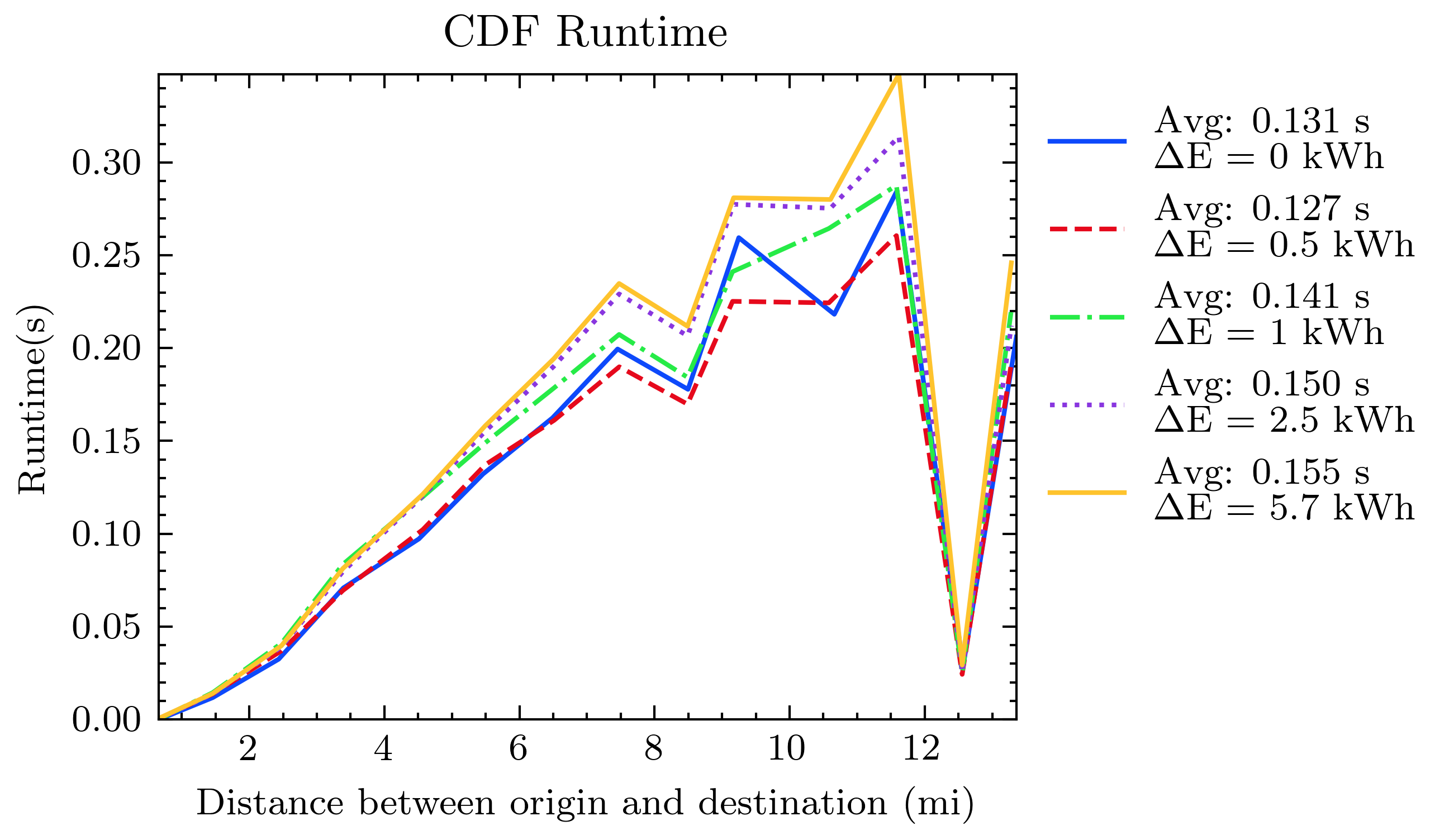}
\caption{}
\label{fig:runtime CDF - dist}
\end{subfigure}
\caption{Average runtime distribution of eco-routing algorithms based on the
shortest distance between O-D pairs}%
\label{fig:runtime dist}
\end{figure*}

In summary, the average energy and travel time savings of different routing
strategies are shown in Tables \ref{tab: energy savings} and
\ref{tab: time savings} respectively. It can be seen that the CRPTC algorithm
always has the best performance with energy savings of up to 19\% compared to
the fastest route. The bi-level eco-routing approach is within 0.2\% of the
CRPTC approach; even though it does not provide the global optimal solution,
its much faster execution time makes it particularly attractive. Note that
when the car has enough energy to travel the entire path with electricity, all
the eco-routing algorithms give us the same results. This is because in this
case the optimal PT control strategy is to use electricity (CD mode) on all
links since it is cheaper than using gas; as a result, CDF and CRPTC give us
the same solutions, as does the bi-level optimization approach. By the same
token, when we do not let the car deplete electricity along the route ($\Delta
E=0$), the optimal solutions of CRPTC, CDF, and bi-level become the same,
since the car can only use gas (CS mode) to travel through the links. In
general, the CRPTC eco-routing algorithm has the best performance when the
distance between origin and destination is relatively high ($>30$ miles),
since in those cases we have more options for choosing when and where to use
the CD mode on the path and solving the combined problem can give us better
results. When we solved the eco-routing problem for the Eastern Massachusetts
highway sub-network \cite{houshmand_eco-routing_2018} the distances were more
than 30 miles with the CRPTC approach outperforming CDF by an average of more
than 2.1\%. As expected, there is a trade-off between energy and time savings,
and time-optimal routes can be more than 20\% faster than the energy-optimal
routes. In order to consider time when solving eco-routing problems, we can
solve problem (\ref{prob: CRPTC with time}) or (\ref{prob: CDF with time})
with different $\alpha$ values to establish a desired balance between time and energy.

\begin{table}[t]
\caption{Average energy cost savings of the proposed eco-routing algorithms ($\Delta E$ values are in kWh)}
\label{tab: energy savings}
\centering
\resizebox{1\columnwidth}{!}{\renewcommand{\arraystretch}{1.15}
\begin{tabular}{c|c|c|c|c|c|c|}
\cline{2-7}
& \multicolumn{6}{c|}{Energy Cost Saving (\%)}                                                                                                                                                                                                                                                                                                                                                          \\ \cline{2-7}
& \begin{tabular}[c]{@{}c@{}}CRPTC \\ vs. \\ Fastest\end{tabular} & \begin{tabular}[c]{@{}c@{}}Bi-level \\ vs. \\ Fastest\end{tabular} & \begin{tabular}[c]{@{}c@{}}CDF \\ vs. \\ Fastest\end{tabular} & \begin{tabular}[c]{@{}c@{}}CRPTC\\  vs.\\  CDF\end{tabular} & \begin{tabular}[c]{@{}c@{}}CRPTC \\ vs.\\  Bilevel\end{tabular} & \begin{tabular}[c]{@{}c@{}}Bi-level\\  vs. \\ CDF\end{tabular} \\ \hline
\multicolumn{1}{|c|}{$\Delta E = 0$}    & 13.8                                                            & 13.8                                                               & 13.8                                                          & 0                                                           & 0                                                               & 0                                                              \\ \hline
\multicolumn{1}{|c|}{$\Delta E = 0.5$}  & 17.9                                                            & 17.7                                                               & 16.4                                                          & 1.8                                                         & 0.2                                                             & 1.6                                                            \\ \hline
\multicolumn{1}{|c|}{$\Delta E = 1$}    & 19                                                              & 18.9                                                               & 18.1                                                          & 1.2                                                         & 0.1                                                             & 1.1                                                            \\ \hline
\multicolumn{1}{|c|}{$\Delta E = 2.5$}  & 15.4                                                            & 15.4                                                               & 15.4                                                          & 0.1                                                         & 0                                                               & 0                                                              \\ \hline
\multicolumn{1}{|c|}{$ \Delta E = 5.7$} & 13.8                                                            & 13.8                                                               & 13.8                                                          & 0                                                           & 0                                                               & 0                                                              \\ \hline
\end{tabular}}\end{table}
\begin{table}[h]
\caption{Average travel time savings results under different routing scenarios ($\Delta E$ values are in kWh)}
\label{tab: time savings}
\centering
\resizebox{1\columnwidth}{!}{\renewcommand{\arraystretch}{1}
\begin{tabular}{c|c|c|c|c|c|}
\cline{2-6}
& \multicolumn{5}{c|}{Travel Time Saving (\%)}                                                                                                                                                                                                                                                                                    \\ \cline{2-6}
& \begin{tabular}[c]{@{}c@{}}Fastest \\ vs.\\  CRPTC\end{tabular} & \begin{tabular}[c]{@{}c@{}}Fastest\\  vs.\\  CDF\end{tabular} & \begin{tabular}[c]{@{}c@{}}Fastest\\ vs.\\ Bi-level\end{tabular} & \begin{tabular}[c]{@{}c@{}}CRPTC\\ vs.\\ CDF\end{tabular} & \begin{tabular}[c]{@{}c@{}}CRPTC\\ vs.\\ Bi-level\end{tabular} \\ \hline
\multicolumn{1}{|c|}{$\Delta E = 0$}    & 22.8                                                            & 22.8                                                          & 22.8                                                             & 0                                                         & 0                                                              \\ \hline
\multicolumn{1}{|c|}{$\Delta E = 0.5$}  & 22                                                              & 22.5                                                          & 22.5                                                             & 1                                                         & 1                                                              \\ \hline
\multicolumn{1}{|c|}{$\Delta E = 1$}    & 21.6                                                            & 22.1                                                          & 22.1                                                             & 0.8                                                       & 0.8                                                            \\ \hline
\multicolumn{1}{|c|}{$\Delta E = 2.5$}  & 21.7                                                            & 21.6                                                          & 21.6                                                             & -0.1                                                      & -0.1                                                           \\ \hline
\multicolumn{1}{|c|}{$ \Delta E = 5.7$} & 21.7                                                            & 21.7                                                          & 21.7                                                             & 0                                                         & 0                                                              \\ \hline
\end{tabular}}\end{table}

\subsection{Algorithm Execution Time Comparison Results}

An important factor in assessing the performance of eco-routing algorithms,
aside from their energy improvement, is their execution time (runtime). It is
essential that an algorithm can compute the eco-route quickly and is able to
re-calculate the energy-optimal route in case of sudden changes in traffic
patterns in the network. As a result, we calculated the execution time of the
three proposed eco-routing algorithms and reported the corresponding averages
in Table \ref{tab: Runtimes}. We have also included the runtime distribution box-plot
of each algorithm in Fig. \ref{fig:runtime box} and showed the execution time
dependency on the distance between O-D pairs in Fig. \ref{fig:runtime dist}.
All these algorithms have been coded in Python 3.7.1 and executed on a desktop
computer with a 4.2GHZ Core i7 CPU and 16 GB of RAM. We used
Gurobi\cite{gurobi} as the MILP solver in this setting. As we can see in Table
\ref{tab: Runtimes}, the CRPTC algorithm runtime is on average 10.34s when
the battery is full ($\Delta E=5.7kWh$). This is an interesting observation since
CRPTC is a MILP problem which is NP-hard. However, as we decrease the
allowable energy depletion from the battery, we see that the runtime of the
CRPTC algorithm starts to increase. This is due to the fact that we are
imposing a tighter constraint (\ref{eqn: CRPTC energy constraint}) to problem
(\ref{prob: CRPTC MILP}) which forces the optimizer to explore more options in
seeking the optimal solution. Recalling that the CRPTC problem has two sets of
decision variables, the routing decision vector $\mathbf{x}$ and PT control
strategy vector $\mathbf{y}$, when the car has enough energy to travel the
entire route with electricity, the optimal PT control strategy is to set
$y_{ij}=1$ for all $(i,j)\in\mathcal{A}$, and the optimizer can find the
optimal solution easily. As we decrease the allowable $\Delta E$ value, we
increase the search space for the optimization problem, consequently the
runtime increases.

As expected, both CDF and bi-level eco-routing algorithms have near real-time
execution times ($\sim150ms$). This is due to the fact that both of these
algorithms use a modified version of Dijkstra's algorithm to find the
energy-optimal paths (Algorithm \ref{algo: CDF-Dijkstra}) which has a
complexity of $O(n\log n)$ where \textit{n} is the number of nodes in the graph.

Another interesting observation is that as the shortest distance between O-D
pairs increases, the runtime of the eco-routing algorithms typically increases
(Fig. \ref{fig:runtime dist}). In particular, CRPTC's runtime is the most
sensitive to the distance between O-D pairs and the runtime can increase more
than 700 times as the distance between O-D pairs increases. Considering the
difference between the runtimes of CRPTC and bi-level, while their respective
performance is virtually indistinguishable (Table \ref{tab: energy savings}),
the use of the bi-level eco-routing algorithm is practically attractive in
urban settings where the O-D pairs are relatively close to each other.

\begin{table}[h]
\caption{Eco-routing algorithms runtime comparison }
\label{tab: Runtimes}
\centering
\resizebox{0.7\columnwidth}{!}{\renewcommand{\arraystretch}{1.15}
\begin{tabular}{c|c|c|c|}
\cline{2-4}
& \multicolumn{3}{c|}{Runtime (s)} \\ \cline{1-4}
\multicolumn{1}{|c|}{$\Delta E (kWh)$}  & CRPTC     & Bi-level    & CDF    \\ \hline
\multicolumn{1}{|c|}{$\Delta E = 0$}    & 6.79      & 0.14        & 0.13   \\ \hline
\multicolumn{1}{|c|}{$\Delta E = 0.5$}  & 409.48    & 0.14        & 0.13   \\ \hline
\multicolumn{1}{|c|}{$\Delta E = 1$}    & 178.28    & 0.15        & 0.14   \\ \hline
\multicolumn{1}{|c|}{$\Delta E = 2.5$}  & 12.34     & 0.16        & 0.15   \\ \hline
\multicolumn{1}{|c|}{$ \Delta E = 5.7$} & 10.34     & 0.16        & 0.15   \\ \hline
\end{tabular}}\end{table}

\begin{figure}[ptb]
\centering
\begin{subfigure}{.24\textwidth}
\centering
\includegraphics[width=0.98\linewidth]{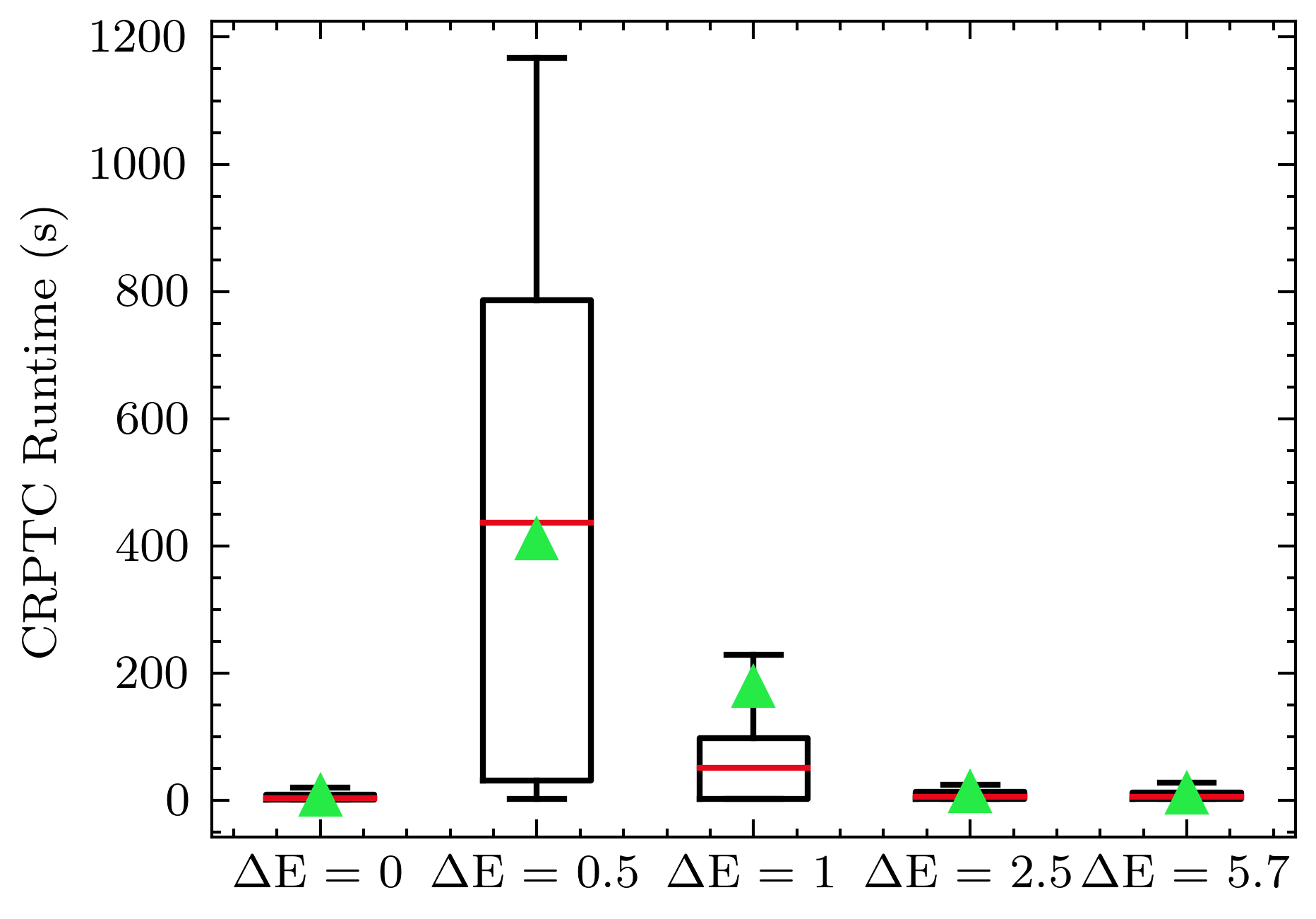}
\caption{}
\label{fig:runtime CRPTC - box}
\end{subfigure}\begin{subfigure}{.24\textwidth}
\centering
\includegraphics[width=0.98\linewidth]{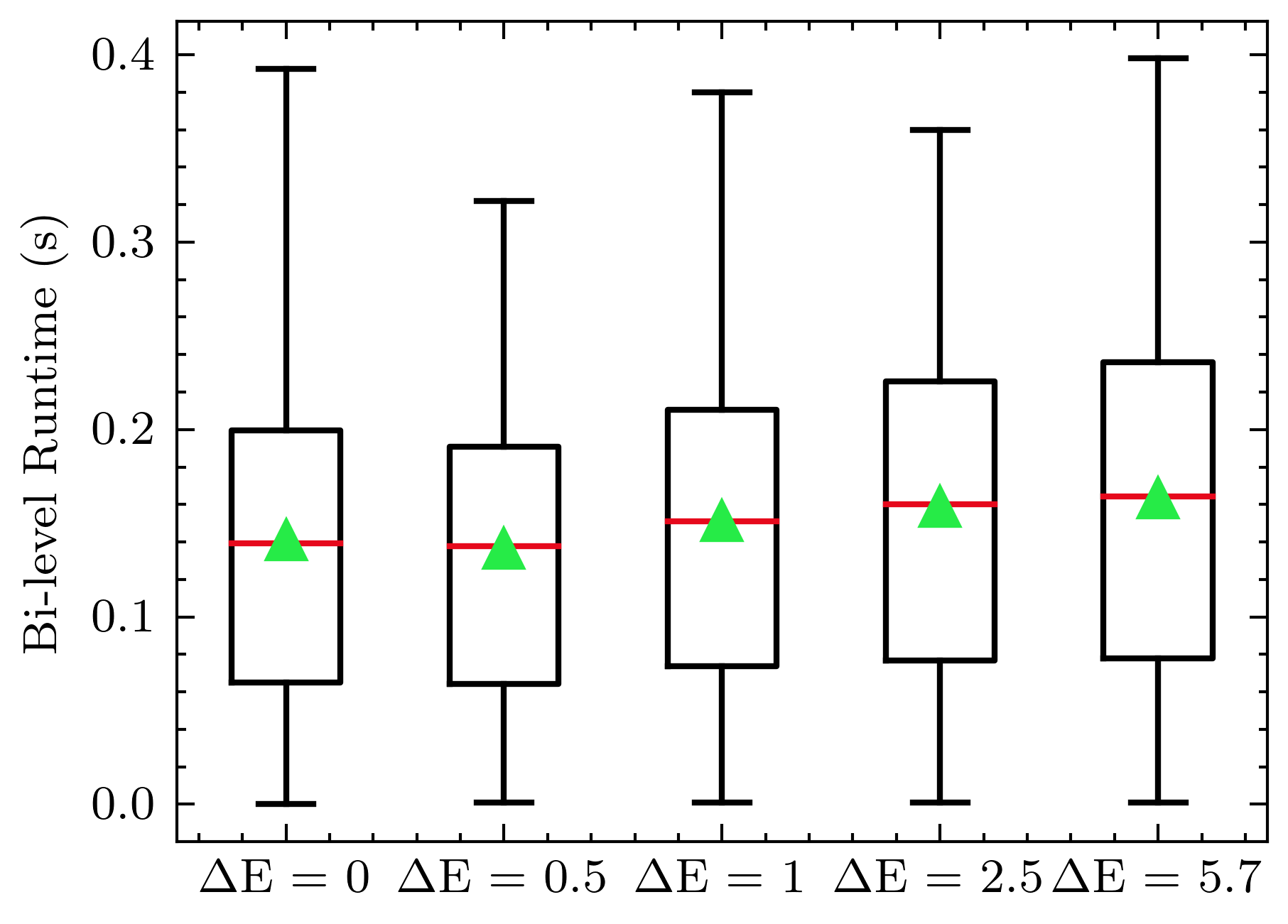}
\caption{}
\label{fig:runtime Bilevel - box}
\end{subfigure}
\begin{subfigure}
{.24\textwidth}
\centering
\includegraphics[width=0.98\linewidth]{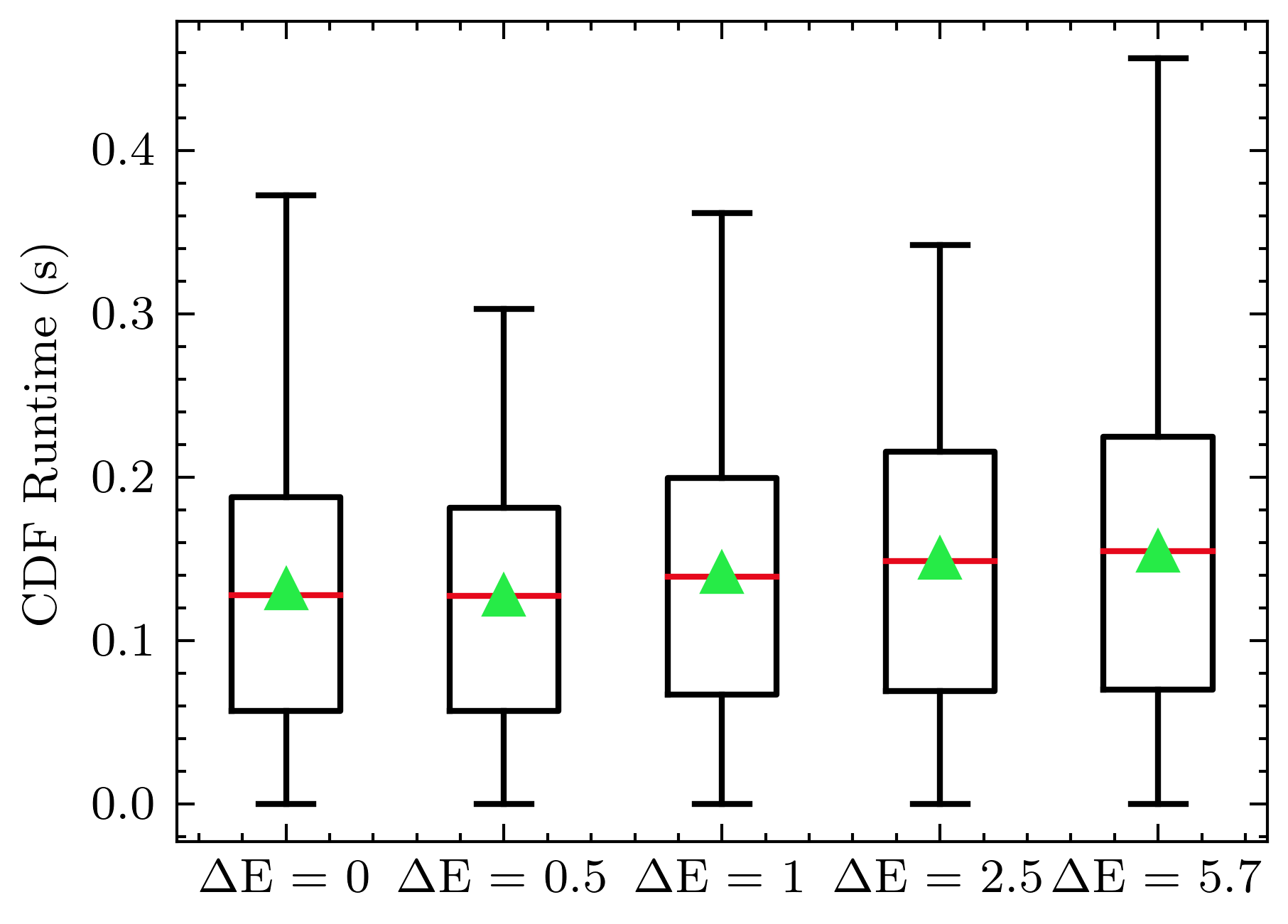}
\caption{}
\label{fig:runtime CDF - box}
\end{subfigure}
\caption{ Average runtime distribution of eco-routing algorithms}%
\label{fig:runtime box}
\end{figure}


\section{Validation Using Simulation Models}\label{sec:validation}

Throughout this paper we have used a simplified energy model to estimate the
energy consumption of PHEVs only based on the traffic intensity on each
road-link. In this section, we investigate the accuracy of this energy model
by using the traffic simulator Simulation of Urban MObility (SUMO)
\cite{krajzewicz_recent_2012} along with a modified version of the
Vehicle-Engine SIMulation (VESIM) model \cite{malikopoulos_simulation_2006}
which is a high fidelity energy modelling software tool calibrated for the
Audi A3 e-tron in Simulink. We use an Audi A3 e-tron since our proposed algorithms have been extensively tested on this vehicle (Fig. \ref{fig:Audi}) at the University of Michigan's M-City and also using Chassis Dyno by Bosch. We start by briefly reviewing our simulation
modeling frameworks in SUMO and VESIM, and then explaining how using them
allows us to validate the accuracy of our eco-routing algorithm.

\begin{figure}[h]
\centering
\includegraphics[width=0.4\textwidth]{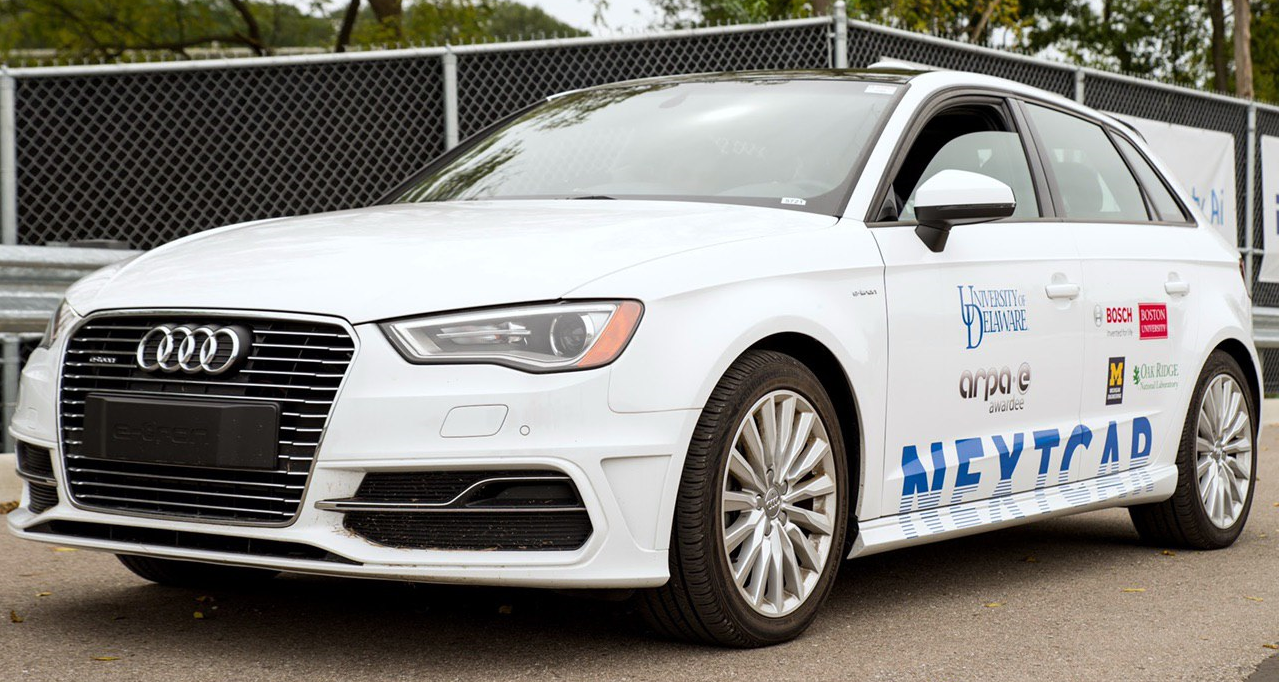}\caption{The Audi A3 e-tron which was used in this study}%
\label{fig:Audi}
\end{figure}

\textbf{SUMO}: We use the SUMO to evaluate the performance of the eco-routing
algorithm. SUMO is an open source traffic simulation package which can generate speed trajectories
of each individual vehicle. In order to have realistic traffic scenarios, we
use the calibrated SUMO model for the Ann Arbor network which we previously
built in \cite{huang2019eco}. To briefly summarize, the Ann Arbor traffic
model consists of 11,265 road segments (links) and 8,660 traffic junctions.
There are 327 traffic lights and 11,857 stop signs, all of which are in
accordance with real-world information. The travel demands used in the SUMO
model are generated according to a calibrated POLARIS model
\cite{auld2016polaris} which is an agent-based mesoscopic traffic simulation
package developed by the Argonne National Lab. The focus of the POLARIS model
is on realistic generation of travel demands based on travel activities of
individual agents within each household using ADAPTS (Agent-based Dynamic
Activity Planning and Travel Scheduling) \cite{auld2012activity}. The travel
demands are calibrated with the dataset from the Safety Pilot Model Deployment
(SPMD) \cite{bezzina2014safety} with records of 321,945 trips in Ann Arbor
between 2013 and 2015 with 2,800 passenger cars, trucks, and buses. The route
choices of individual vehicles in SUMO are then calibrated again to ensure
that the average speed of each road in the simulation match the observed
average speeds from the SPMD dataset. \begin{figure}[h]
\centering
\includegraphics[width=0.4\textwidth]{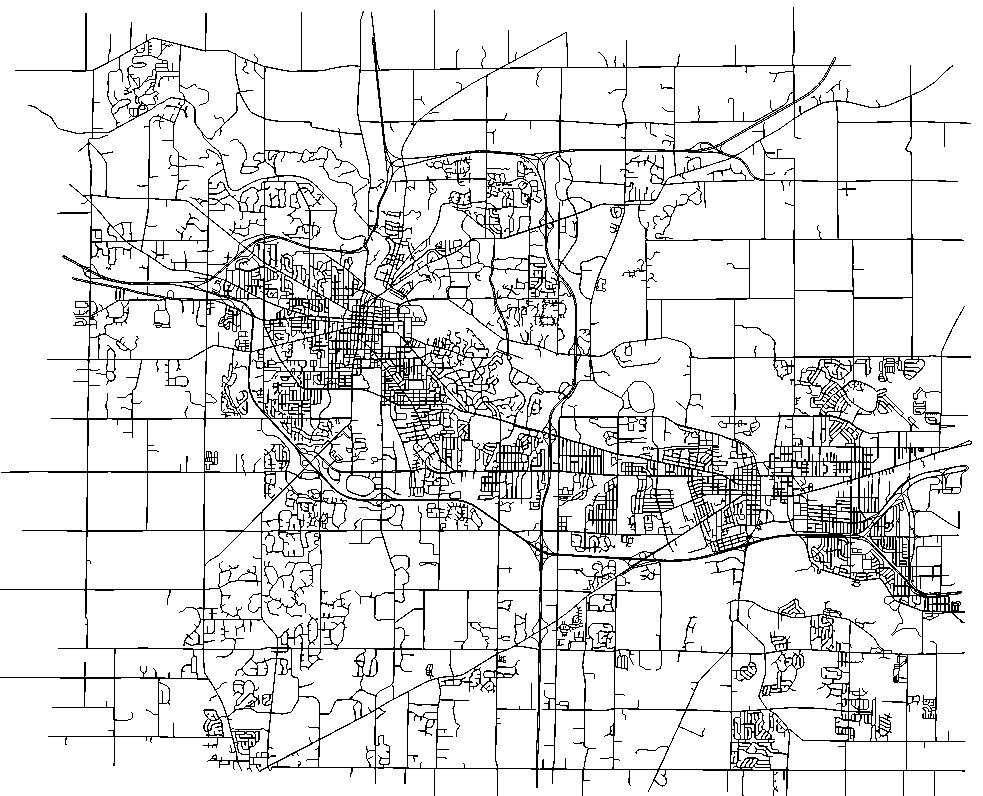}\caption{Ann
Arbor's traffic network in SUMO}%
\end{figure}

\textbf{VESIM}: This is a Simulink based power-train modeling framework to
calculate energy costs for any given speed profile. Our VESIM model is
calibrated for an Audi A3 e-tron which is a PHEV, and the engine and electric
motor efficiency maps are modified to match that of the Audi. A schematic of the VESIM model in Simulink is shown
in Fig. \ref{fig:VESIM}. \begin{figure}[h]
\centering
\includegraphics[width=0.4\textwidth]{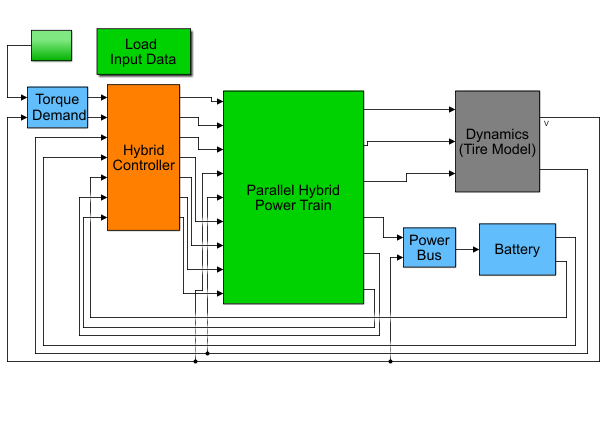}\caption{{VESIM
Model for modelling the power-train of the plug-in hybrid electric Audi A3
e-tron\cite{mahbub2019concurrent}}}%
\label{fig:VESIM}
\end{figure}

\subsection{Eco-routing performance validation}

We consider the time-optimal path as the baseline for measuring the
performance of eco-routing algorithms. In order to validate the effectiveness
of our eco-routing algorithm we need to compare the energy cost of travelling
through the time-optimal path to that of the energy-optimal path. Since it is
very difficult and costly to conduct such experiments in a real-world setting,
we use computer simulation software to perform this task. The procedure of
using SUMO and VESIM to validate the energy saving results is as follows:
\end{subequations}
\begin{enumerate}
\item Use the calibrated SUMO model of Ann Arbor to simulate traffic in the network.

\item Choose an O-D pair and find the eco-route and fastest route between the
selected origin and destination.

\item Send two vehicles in SUMO to follow the fastest route and eco-route and
collect their respective speed trajectories.

\item Import the speed trajectories from SUMO to VESIM and calculate the
energy costs of both the fastest route and the eco-route.

\item Compare the energy costs of the eco-route and fastest route and report
the energy savings.
\end{enumerate}

\begin{figure}[h]
\centering
\includegraphics[width=0.25\textwidth]{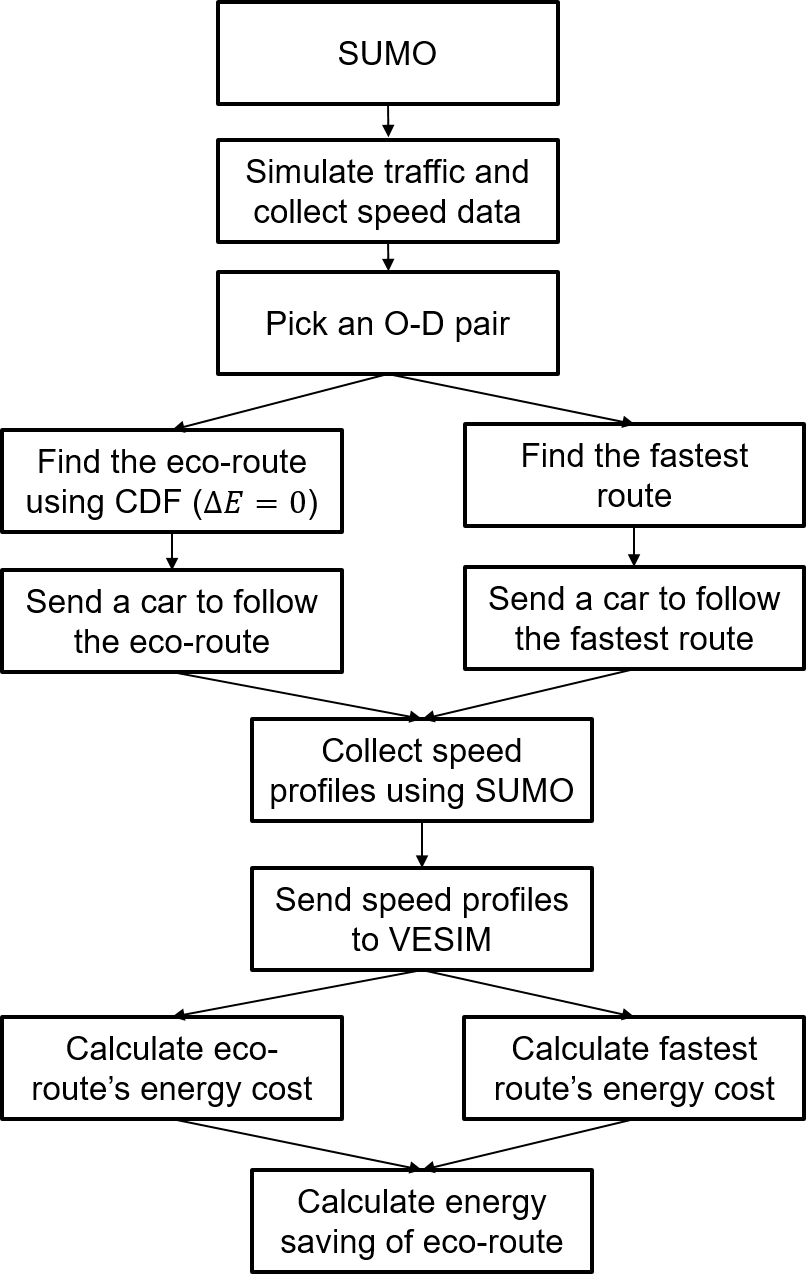}\caption{{Procedure
for calculating energy costs using SUMO and VESIM}}%
\label{fig:SUMO+VESIM Framework}%
\end{figure}

\begin{figure}[h]
\centering
\includegraphics[width=0.48\textwidth]{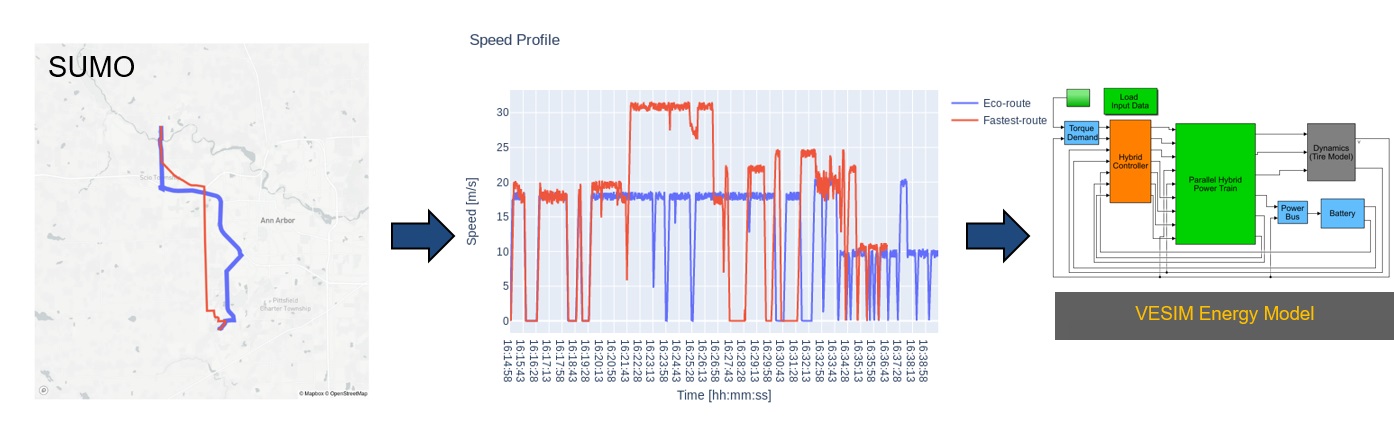}\caption{Procedure
of collecting speed profiles from SUMO and sending them to VESIM to validate
the performance of the eco-routing algorithm}%
\label{fig:SUMO+VESIM example}%
\end{figure}

Since our CRPTC eco-routing algorithm (\ref{prob: CRPTC MILP}) simultaneously
finds both the energy-optimal route and the optimal switching strategy between
CD and CS modes, accurate energy results require us to consider the optimal PT
control strategy commands from CRPTC while finding energy values using VESIM.
So far, we have not incorporated the CRPTC control decisions for the PT
controller into VESIM; as a result, to simplify this process, we consider a
hybrid electric vehicle (HEV) instead of a PHEV, which only operates in CS
mode (setting $\Delta E=0$ while solving the eco-routing problem). Hence, we
use Algorithm (\ref{algo: CDF-Dijkstra}) with $E_{0}=0$ to find the eco-route,
and in VESIM we choose the CS mode and let the SOC change $\pm5\%$ throughout
the route (Figs. \ref{fig:SUMO+VESIM Framework} and
\ref{fig:SUMO+VESIM example}). 

\subsubsection{Results Validation}

Using the aforementioned framework (Fig. \ref{fig:SUMO+VESIM Framework}) we
randomly select 2200 O-D pairs in the Ann Arbor SUMO network and calculate
eco-routes and fastest routes for each of these pairs. We then import the
collected speed profiles from SUMO into VESIM and calculate the energy costs
savings. We repeat this procedure under two different traffic conditions in
the SUMO network: medium and high traffic. The results for the medium traffic
network are summarized as follows:

\begin{itemize}
\item The average actual energy saving of eco-routes vs. fastest routes
calculated by VESIM over the selected O-D pairs  is \textbf{12.59\%}.

\item The average predicted energy saving using our simplified energy model
(\ref{eqn: CDF cost}) is \textbf{12.52\%}.

\item In 80.38\% of cases, we correctly predicted the energy optimal routes and
the average savings of these cases is 18.46\%.

\item In 19.62\% of cases, the predicted eco-routes consumed more energy than
the fastest route, which we refer to as \textquotedblleft false
positives\textquotedblright. The average energy loss of false positive cases
is -11.46\%.
\end{itemize}

\begin{table}[ptb]
\caption{VESIM eco-routing validation results for the 2200 O-D pairs selected
in SUMO under medium and high traffic network traffic }%
\label{tab:VESIM validation results}
\centering
\resizebox{0.99\columnwidth}{!}{\renewcommand{\arraystretch}{1.2}
\begin{tabular}{|c|c|c|c|c|c|}
\hline
& \begin{tabular}[c]{@{}c@{}}Actual average \\ energy savings \\ (VESIM Results)\end{tabular} & \begin{tabular}[c]{@{}c@{}}True Positive \\ rate\end{tabular} & \begin{tabular}[c]{@{}c@{}}True Positive\\ energy savings\end{tabular} & \begin{tabular}[c]{@{}c@{}}False Positive\\ rate\end{tabular} & \begin{tabular}[c]{@{}c@{}}False Positive\\ energy saving\end{tabular} \\ \hline
Medium Traffic & 12.59\%                                                                                     & 80.38\%                                                       & 18.46\%                                                                & 19.62\%                                                       & -11.46\%                                                               \\ \hline
High Traffic   & 5.35\%                                                                                      & 69.00 \%                                                      & 17.47\%                                                                & 31.00\%                                                       & -21.64\%                                                               \\ \hline
\end{tabular}}\end{table}As we can see in the medium traffic network, the
expected actual energy saving of eco-route vs. fastest route is almost the
same as the predicted energy saving using our simplified energy model.
However, the energy model does not always predict the eco-route correctly and
sometimes finds a route that in reality consumes more energy than the fastest
route (Table \ref{tab:VESIM validation results}). One of the reasons is that
we classified links based on their traffic intensity into only three
categories: low, medium, and high traffic intensity links. We also assign
standard speed profiles (Tables \ref{tab: Drive cycle assignment} and
\ref{tab: conversion factors}) to each of these links and use their average
energy consumption as the predicted energy consumption over the link. In order
to increase the accuracy of the energy model, we can increase the number of
link categories and assign more suitable speed profiles to each link based on
different factors such as free flow speed, location of the link
(urban/highway), road grade, traffic lights/stop signs, etc. Moreover, we
update the average speed data in SUMO every 15 minutes. However, since the Ann
Arbor network includes many traffic lights and there is often high traffic in
the network, traffic conditions may significantly deviate within 15 minutes.
This behavior is more evident in the high traffic network case (Table
\ref{tab:VESIM validation results}). As a result, we may need to extend our
analysis to dynamic re-routing whenever traffic changes occur in the network,
expecting to improve the accuracy of our results. 

\begin{figure}[h]
\centering
\begin{subfigure}{.24\textwidth}
\centering
\includegraphics[width=0.98\linewidth]{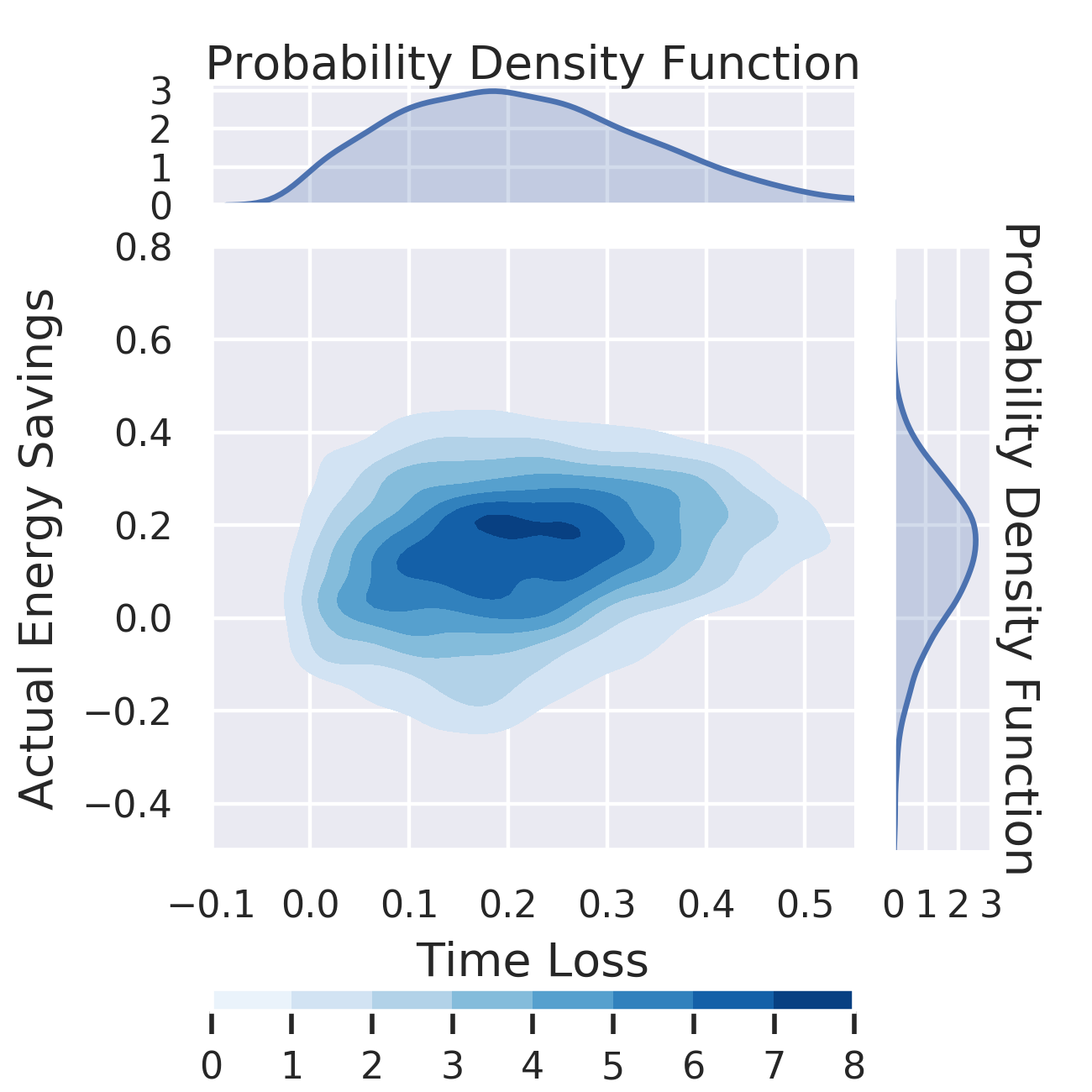}
\caption{Medium traffic network}
\label{fig:UI dl}
\end{subfigure}\begin{subfigure}{.24\textwidth}
\centering
\includegraphics[width=0.98\linewidth]{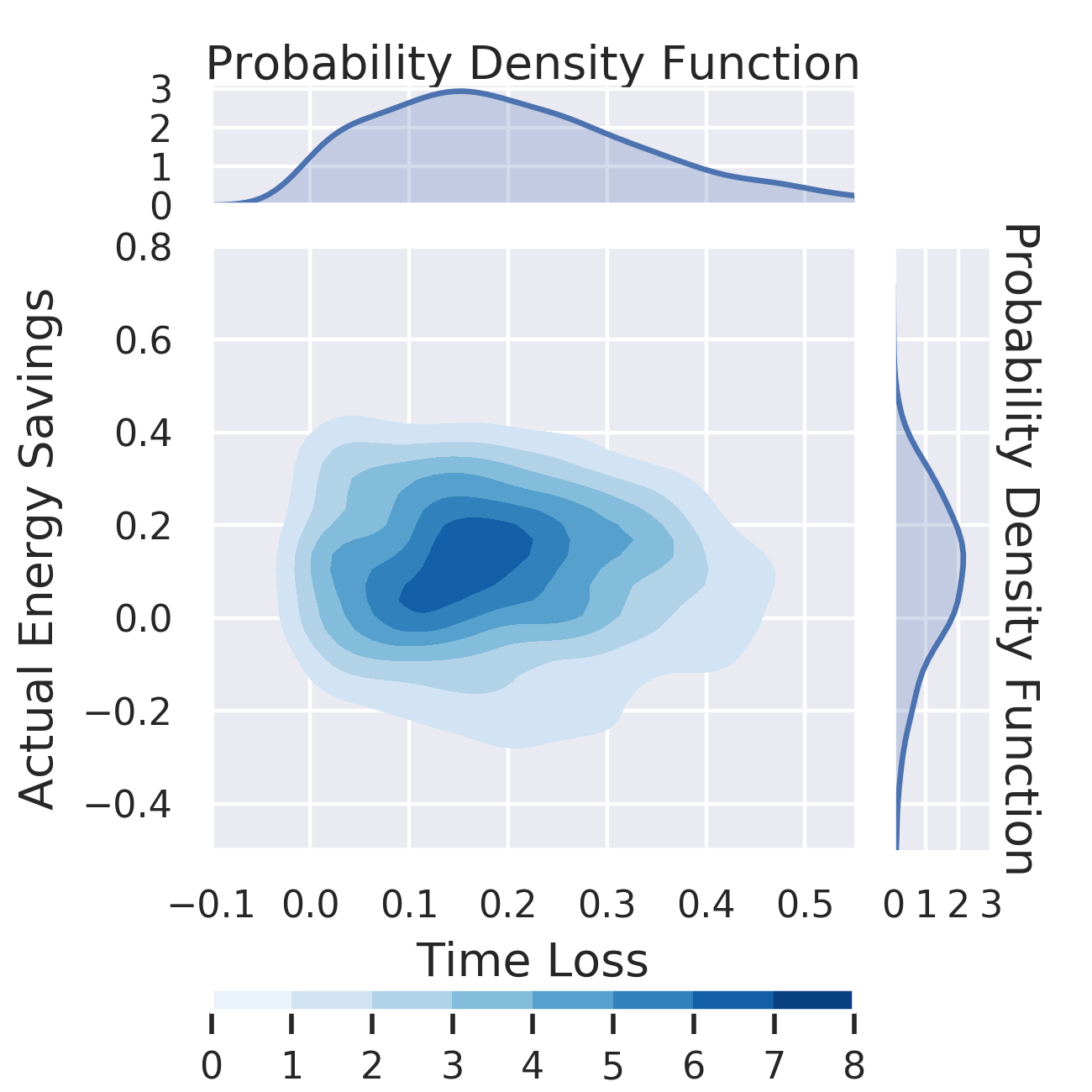}
\caption{High traffic network.}
\label{fig:UI vis}
\end{subfigure}
\caption{ SUMO+VESIM validation results: Trade-off between energy and time
saving in eco-routes (colors represent the probability density values)}%
\label{fig:SUMO Kernel plots}%
\end{figure}

\begin{figure}[h]
\centering
\begin{subfigure}{.24\textwidth}
\centering
\includegraphics[width=0.98\linewidth]{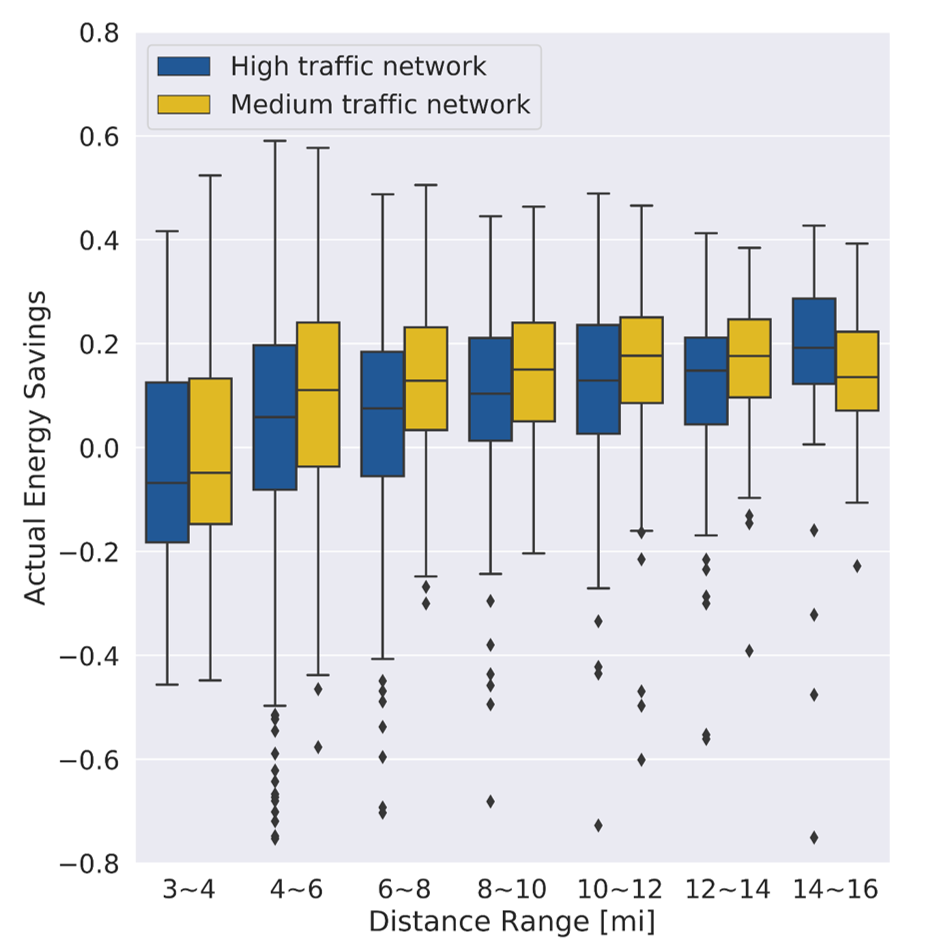}
\caption{Energy saving}
\label{fig:UI dl}
\end{subfigure}\begin{subfigure}{.24\textwidth}
\centering
\includegraphics[width=0.98\linewidth]{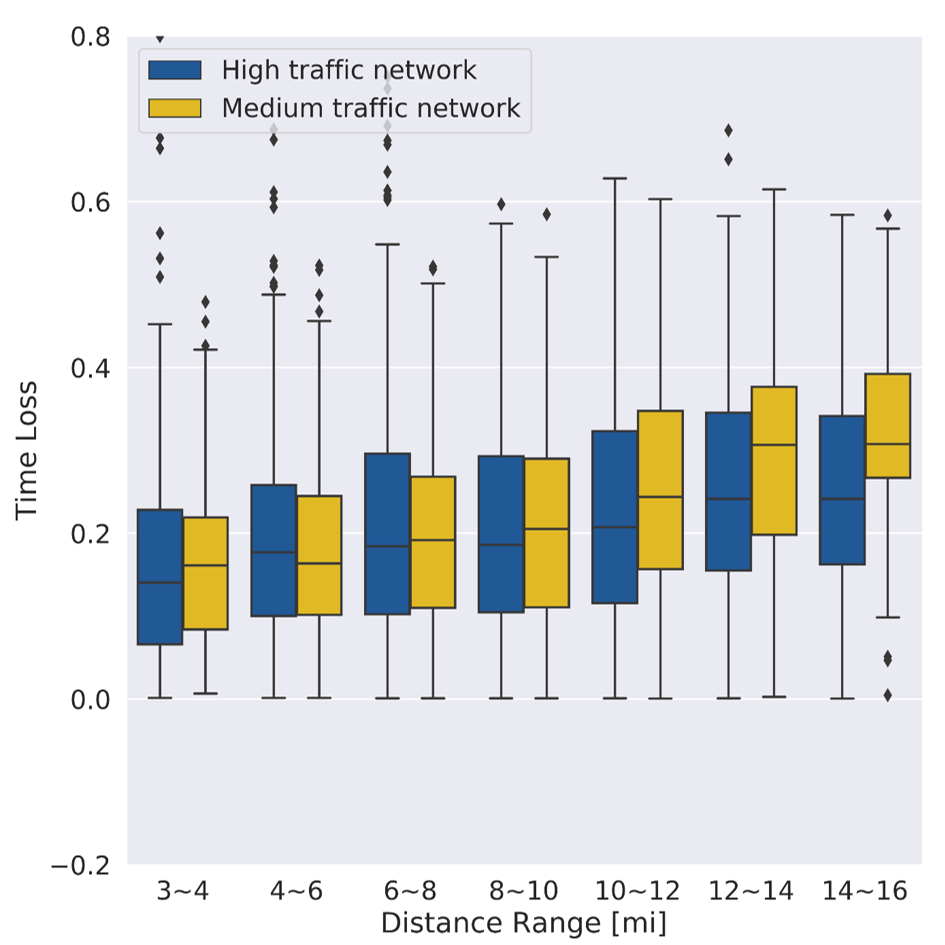}
\caption{Time loss}
\label{fig:UI vis}
\end{subfigure}
\caption{ SUMO+VESIM validation results: Energy saving and time loss of
eco-route vs. fastest route as a function of distance between O-D pairs.}%
\label{fig:validation energy and time}%
\end{figure}

The relationship between energy saving and time loss is shown in Fig.
\ref{fig:SUMO Kernel plots} where the colored contours show the probability
density values. We use kernel density estimation to estimate the probability density function. A Kernel density estimator can be viewed as a special case of the Gaussian mixture model with the weight of each component set to $n^{-1}$.   Moreover, Fig. \ref{fig:validation energy and time} shows the
distribution of energy savings of the eco-route and time savings of the
fastest route as a function of the shortest distance between O-D pairs.



\section{Conclusions and Future Work}

\label{sec: conclusions}

We have proposed two methods to solve the minimum-energy cost problem for a
single vehicle routing: the CRPTC and bi-level optimization algorithms. The
proposed methods are capable of finding both the optimal path and the optimal
switching strategy between CD and CS modes on each link, and can be
implemented in real time. An open source framework for downloading traffic
data was also developed for this work and was used to implement the
eco-routing algorithms on a large urban transportation network. The
performance of the eco-routing algorithm was validated by using SUMO and VESIM
and the results show energy savings of more than 12\%. We have also shown that
there is a trade-off between energy saving and time saving.

So far, we have not considered dynamically updating routing decisions at
network nodes to account for sudden changes in traffic conditions (e.g., due
to accidents). We have also limited our analysis to a single vehicle scenario
with a known origin and destination. As a next step, we will consider
connectivity among vehicles and determine the social optimum for the network
considering 100\% penetration rate of connected automated vehicle \cite{houshmand_penetration_2019}. Moreover, we plan to include multiple
vehicle architectures with different fuel consumption models and different
initial energies to the problem, as well as adding charging stations into the
network to let vehicles recharge their batteries if necessary.

\bibliographystyle{IEEEtran}
\begin{tiny}
\bibliography{NEXTCAR}
\end{tiny}

\appendix  \label{Appendix A}
{\tiny \begin{algorithm}[H]
\caption{{\bf CDF-Dijkstra} \label{algo: CDF-Dijkstra}}
\begin{algorithmic}[]
\Procedure {CDF}{$G$, origin, destination, $E_0$}
\State create node set $Q$
\ForAll {\(v \in \mathcal{N}\)}
\State cost[origin] $\leftarrow \infty$
\State prev[v] $\leftarrow$ UNDEFINED
\State add $v$ to $Q$
\EndFor
\State $cost[origin] \leftarrow 0$
\State $E[origin] \leftarrow E_0$
\While {$Q$ is not empty}
\State $u \leftarrow$ vertex in $Q$ with min cost[$u$]
\State remove $u$ from $Q$
\If{$u = $ destination}
\State $S \leftarrow$ empty sequence
\State $u \leftarrow$ destination
\If{prev$[u]$ is defined \textbf{or} $u =  origin$}
\While {$u$ is defined}
\State insert $u$ at the beginning of $S$
\State $u \leftarrow$ prev$[u]$
\EndWhile
\EndIf
\State break the while loop
\EndIf
\ForAll { neighbor $v$ of $u$}
\If{energy$[u] \geq \frac{dist[u,v]}{\mu_{CD}[u,v]}$}
\State cost$[u,v] \leftarrow C_{ele}\frac{dist[u,v]}{\mu_{CD}[u,v]}$
\State $E_{temp} \leftarrow E[u] - \frac{dist[u,v]}{\mu_{CD}[u,v]}$
\Else
\State{
\begin{equation*}
\begin{array}{llllll}
& & & & cost[u,v] \leftarrow  & C_{gas}\frac{dist[u,v]-\mu_{CD}[u,v]E[u]}{\mu_{CS}[u,v]}\\
& & & & & + C_{ele}E[u]
\end{array}
\end{equation*}}
\State $E_{temp} \leftarrow 0$
\EndIf
\State alt $\leftarrow$ cost$[u]$ + cost$(u,v)$
\If{alt $<$ cost$[v]$}
\State cost$[v] \leftarrow$ alt
\State $E[v] \leftarrow E_{temp}$
\State prev$[v] \leftarrow u$
\EndIf
\EndFor
\EndWhile
\State \textbf{return} $cost[ \text{ }], S[\text{ } ]$
\EndProcedure
\end{algorithmic}
\end{algorithm}
}

\end{document}